\RequirePackage{amsthm}
\documentclass[pdflatex,sn-mathphys-num]{sn-jnl}% Math and Physical Sciences Numbered Reference Style 
\usepackage{graphicx}%
\usepackage{multirow}%
\usepackage{amsmath,amssymb,amsfonts}%
\usepackage{amsthm}%
\usepackage{mathrsfs}%
\usepackage[title]{appendix}%
\usepackage{xcolor}%
\usepackage{textcomp}%
\usepackage{manyfoot}%
\usepackage{booktabs}%
\usepackage{algorithm}%
\usepackage{algorithmicx}%
\usepackage{algpseudocode}%
\usepackage{listings}%
\usepackage{float}
\usepackage{caption}
\usepackage{lipsum}
\usepackage{lmodern}
\usepackage{anyfontsize}
\usepackage{lscape}
\usepackage{anyfontsize}
\theoremstyle{thmstyleone}%
\newtheorem{theorem}{Theorem}%  meant for continuous numbers
\newtheorem{corollary}{Corollary}

\theoremstyle{thmstyletwo}%

\theoremstyle{thmstylethree}%

\usepackage{paralist}
\usepackage{color,soul}
\usepackage{subcaption}
\usepackage{latexsym}
\usepackage{flexisym}
\usepackage{comment}
\usepackage{lipsum}
\usepackage{mathtools}

\usepackage{caption}
\usepackage{graphicx}

\usepackage{xcolor}
\usepackage{rotating}    
\usepackage{makecell}   
\usepackage{tabularray}
\UseTblrLibrary{booktabs, counter, varwidth}
\usepackage{ragged2e}

\DeclarePairedDelimiter\floor{\lfloor}{\rfloor}

\newcommand\Mycomb[2][^n]{\prescript{#1\mkern-0.5mu}{}C_{#2}}

\newcommand{\black}{\color{black}}

\newcommand{\insertpar}[1]{\noalign{\vskip6pt{\nipar #1}\vskip6pt}}

\def\nipar{\par\noindent\ignorespaces }

\begin{document}

\title[Branching Adaptive Surrogate Search
Optimization (BASSO)]{Branching Adaptive Surrogate Search
Optimization (BASSO)}

%%=============================================================%%
%% GivenName	-> \fnm{Joergen W.}
%% Particle	-> \spfx{van der} -> surname prefix
%% FamilyName	-> \sur{Ploeg}
%% Suffix	-> \sfx{IV}
%% \author*[1,2]{\fnm{Joergen W.} \spfx{van der} \sur{Ploeg} 
%%  \sfx{IV}}\email{iauthor@gmail.com}
%%=============================================================%%

\author*[1]{\fnm{Pariyakorn} \sur{Maneekul}}\email{parim@uw.edu}

\author[1]{\fnm{Zelda B.} \sur{Zabinsky}}\email{zelda@uw.edu}
%\equalcont{These authors contributed equally to this work.}

\author[2]{\fnm{Giulia} \sur{ Pedrielli}}\email{giulia.pedrielli@asu.edu}
%\equalcont{These authors contributed equally to this work.}

\affil[1]{\orgdiv{Industrial and Systems Engineering Department}, \orgname{University of Washington}, \orgaddress{\city{Seattle}, \postcode{98195-2650}, \state{WA}, \country{USA}}}

\affil[2]{\orgdiv{School of Computing and Augmented Intelligence}, \orgname{Arizona State University}, \orgaddress{\city{Tempe}, \postcode{85281}, \state{AZ}, \country{USA}}}

%%==================================%%
%% Sample for unstructured abstract %%
%%==================================%%

\abstract{Global optimization of black-box functions is challenging in high dimensions. We introduce a conceptual adaptive random search framework, Branching Adaptive Surrogate Search Optimization (BASSO), that combines partitioning and surrogate modeling for subregion sampling. 
%The approach is inspired by the desirable complexity of theoretical adaptive random search approaches that are effective for global optimization problems, where the expected  time to reach a target region increases linearly with dimension. 
We present a finite-time analysis of BASSO, and establish conditions under which it is theoretically possible to scale to high dimensions. While we do not expect that any implementation will achieve the theoretical ideal, we experiment with several BASSO variations and discuss implications on narrowing the gap between theory and implementation.
Numerical results on test problems are presented.}

\keywords{Black-box global optimization,
Derivative-free optimization, 
Adaptive random search, 
Surrogate model, 
Gaussian process, 
Boltzmann distribution}

%%\pacs[JEL Classification]{D8, H51}

%%\pacs[MSC Classification]{35A01, 65L10, 65L12, 65L20, 65L70}

\maketitle

\section{Introduction}\label{sec:Intro}
Optimization of a black-box function is a classical problem studied in several disciplines and can be found in many real-world applications including experimentation, engineering tasks, industrial design, economics, biology, model estimation in statistics and machine learning. These applications have objective functions that are typically non-linear and non-convex in nature, expensive to evaluate and not known in closed form. New technologies in sensing, computing, and storage with cloud and edge solutions have resulted in increasingly large numbers of parameters (i.e., decision variables) with increasing complex interactions representing one of the major challenges for black-box optimization. An example of such complexity arises in machine learning applications with models with a huge number of hyper-parameters to be learned from increasingly large data sets. %This adds to the complexity of these problems, which are , requiring the use of black-box techniques. 
%In fact, black-box global optimization methods present important challenges in scaling to high dimensions and large data samples, challenges that typically will occur at the same time. 

In order to develop scalable algorithms, it is important to gain insight into the impact of high dimensions on the sampling distribution. 
 The impact of high dimensions on computation has been studied through finite-time analyses of several stochastic adaptive search methods, including, Pure Adaptive Search (PAS)~\cite{Zabinsky1992, Zabinsky1995, zabinsky2003stochastic} and Hesitant Adaptive Search (HAS)~\cite{Bulger1998,Wood2001,HASELinzZabinsky}. {\black In particular, under certain conditions, the expected number of PAS (and HAS) function evaluations required to sample below a specified objective
function value increases only {\it linearly} in the dimension of the input when optimizing a function without noise. We refer to this as the ``linearity result." 
While PAS and HAS {\black are not directly implementa{ble, the theoretical analyses contribute important insights} to global optimization. 

The introduction of Annealing Adaptive Search (AAS)~\cite{Shen2005, Shen2007, WoodZabinskyChapter}
was an attempt to narrow the gap between theory and implementation.  AAS was created as a conceptual form of simulated annealing {\black where points are sampled over the original domain according to the Boltzmann distribution parameterized by temperature. The desired linearity result can be achieved by} the conceptual AAS with a derived cooling schedule for the temperature parameter~\cite{Shen2007}. 
However, it is still impractical to sample  exactly from a Boltzmann distribution, so an implementation is still illusive.
%{\black The primary challenge in implementing PAS, HAS, and AAS, is identifying practical sampling distributions that satisfy the assumptions needed for the ideal performance.} 

Another situation where theory was used to aid implementation was in \cite{Zabinsky2010}. Since cooling the temperature in simulated annealing too quickly may result in premature convergence, the HAS theory was adapted to model the behavior of simulated annealing and similar stochastic optimization algorithms. The combination of theory and observations are used to determine a stopping and restarting strategy that resulted in an
implementable algorithm, Dynamic Multistart Sequential Search (DMSS).

Many optimization algorithms have been developed for black-box optimization that search the whole domain for the global optimum, as reviewed in \cite{Audet2020, AudetHare, RiosSahinidis}. 
{\black Surrogate models of the function to optimize, especially Gaussian processes, have been used widely to determine the next point in the domain to evaluate~\cite{rasmussen2006gaussian,gramacy2020surrogates,frazier2018tutorial}.}
Partition-based algorithms, as in \cite{ShiNestedPartition,zabinsky2019PBnB}, {\black create smaller subregions in order to prioritize the exploration of promising areas} as a way to dynamically update the sampling distribution.  Partitioning combined with surrogate models, as in \cite{DIRECT,BIRECT,DDSBB,pedrielli2023part,snobfit2008,BAM}, have been shown to be successful in many applications. 
The prior analyses of PAS, HAS, and AAS are not immediately applicable to surrogate modeling and partition-based algorithms.
The question remains, 
{\it is it possible for surrogate modeling and partition-based algorithms to achieve the ideal linearity result?  And if so, what is needed to make these algorithms scalable to high dimensions?}

To address this, we introduce Branching Adaptive Surrogate Search Optimization (BASSO), which conceptualizes the use of branching and surrogate modeling.  We present a new finite-time analysis of BASSO and prove that the desired linearity result is achievable with two assumptions. 
%While Assumptions 1 and 2 are strong, they provide insight into what is needed in an algorithm to approximate the ideal linearity result. 
{\black Specifically, we establish conditions under which the expected number of function evaluations required to reach an $\epsilon$-optimal solution is bounded by an expression that is linear in dimension.}
While these conditions are impractical to satisfy completely in practice, we explore them with numerical implementations and discuss insights gained.
{\black  The analysis is valid when the domain is continuous and when the domain is} finite, such as a bounded integer lattice.
}}

The main contribution of this paper is the introduction of BASSO
with   
a finite-time analysis of performance,
providing a framework with partitioning schemes and surrogate modeling within subregions that are effective for black-box global optimization. 
While PAS, HAS, and AAS are not directly implementable,
a contribution of the BASSO framework is that it provides conditions that would make an algorithm scalable to high dimensions.  
The theoretical result allows researchers to make informed trade-offs between exploration and exploitation in designing algorithms.

Our second contribution is exploring several intuitive ways to implement BASSO, and evaluating which variations come closest to satisfying the conditions. The variations contrast exploitation with exploration to select subregions.  We also consider two surrogates, a Gaussian process and  a regularized quadratic regression, and compare to uniform sampling (i.e., no surrogate). In the numerical section, we discuss the advantages and disadvantages of these variations.  We also perform numerical experiments with three publicly available  algorithms to test performance in high dimensions.
 We do not expect that any implementation will achieve the theoretical ideal, but combining the theory with numerical experiments provides insights to push the boundary of tackling high dimensional  black-box problems.

BASSO is presented in Section~\ref{sec:BASSO}, and the finite-time analysis with important assumptions is in Section~\ref{sec:analysis}.
Numerical experiments testing the assumptions and scaling to high dimensions are presented in Section~\ref{sec:expts}.
The paper is concluded in Section~\ref{sec:conclusion}.

\section{Branching Adaptive Surrogate Search Optimization (BASSO)}
\label{sec:BASSO}

The BASSO algorithm uses 
%partitioning and surrogate modeling to direct the adaptive sampling distribution.
%Whereas PAS and HAS are not directly implementable, BASSO is computable and numerical results are presented later.   
%Specifically, BASSO uses 
the information obtained on sampled points and their observed objective function values to partition the domain, construct surrogate models, and update a sampling procedure that focuses computational effort on promising subregions.

%We use the following notation throughout.  
The  optimization problem we consider is  
\begin{equation} \tag{P}
\label{eq:optproblem}
\underset{x \in S}{\text{min}} \ f(x) %\nonumber
\end{equation}
where $f(x):S \rightarrow \mathbb{R}$, and $S \subset 
\mathbb{R}^{d}$.
Assume $S$ is compact. The feasible region $S$ may be a continuous space, or a finite space.
The decision variable $x$ is a vector in $d$ dimensions.
%, and the values may be integer- or real-valued. 
We assume a unique minimum, and denote the minimum value and the optimal
point in the domain, respectively, as:
\[
y_* = \underset{x \in S}{\text{min}} \ f(x) \quad \textrm{and} \quad x_* = \underset{x \in S}{\text{arg min}} \ f(x).
\]
%Similarly, denote the maximum value and a maximal point as:
%\[
%y^* = \underset{x \in S}{\text{max}} \ f(x) \quad \textrm{and} \quad x^* = \underset{x \in S}{\text{arg max}} \ f(x).
%\] 
Similarly, denote the maximum value as
$
y^* = \underset{x \in S}{\text{max}} \ f(x)$.

The adaptive sampling scheme in BASSO  to generate a new sample point on the $k$th iteration has three main components: 
\begin{enumerate}
\item[(i)]  Branch a collection of subregions.
%in $\Sigma^C_{k}$. 
\item[(ii)]  Evaluate the adaptive subregion probabilities 
$\tilde{p}_i$ for current subregions $\sigma_i$
%  Select a subregion using the adaptive subregion probabilities 
%$\tilde{p}_i(\tilde{y}_k^*)$ for subregions $\sigma_i$ in $\Sigma^C_{k}$
and select a promising subregion for sampling. 
\item[(iii)] Generate a point within the selected subregion according to the uniform distribution
or a surrogate model.
\end{enumerate}

\begin{figure}[ht!]
    \centering
    \includegraphics[width=0.6\textwidth]{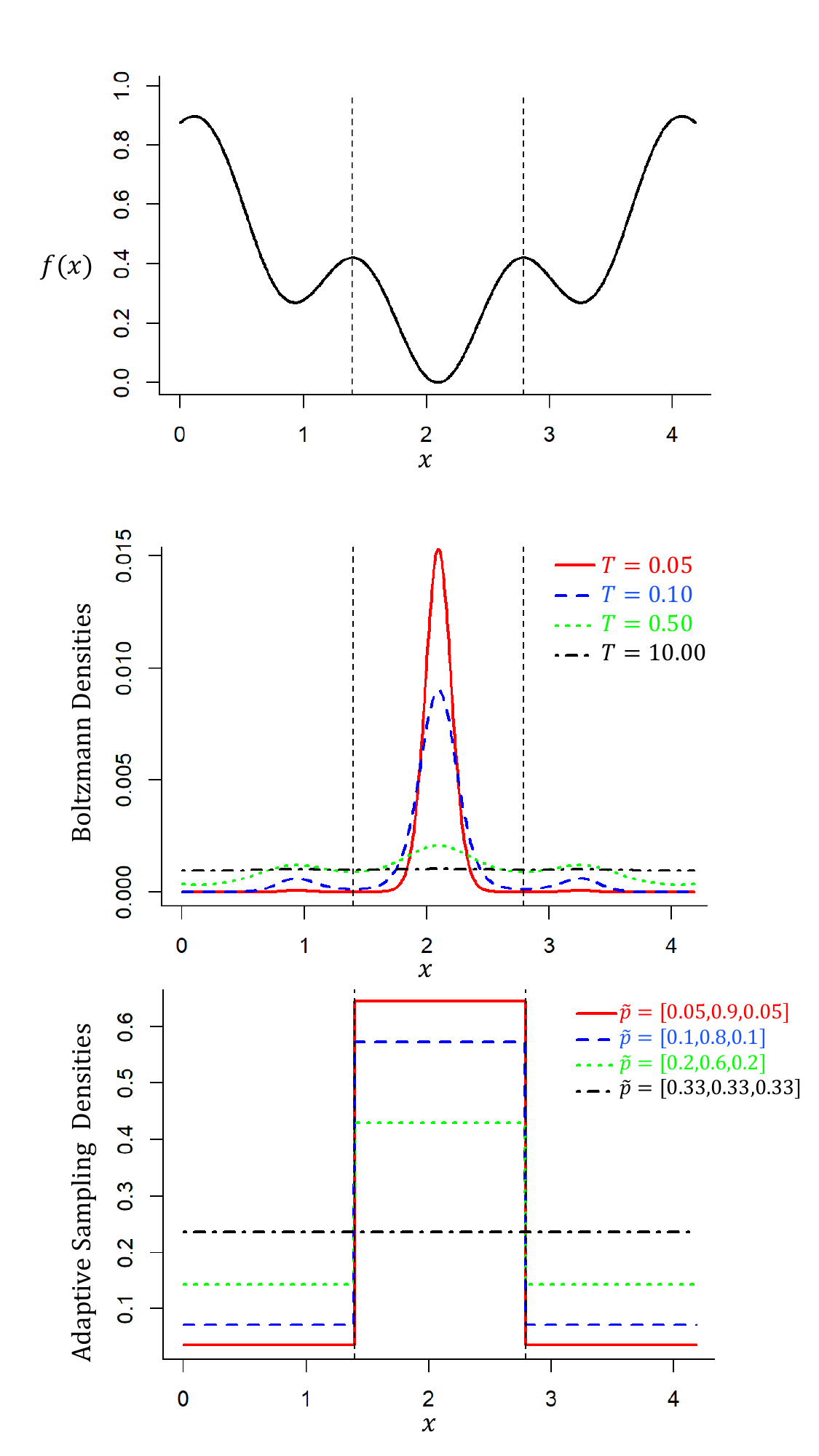}
    \caption{Sampling from a Boltzmann distribution parameterized by temperature induces the densities to focus more on the optimum as the temperature decreases. In BASSO, the adaptive subregion probability parameter $\tilde{p}$ impacts the sampling distribution.  
    %We hope that, as BASSO progresses,  $\tilde{p}$ will gradually focus on the subregions that contain the optimum.
    }
    \label{fig:boltzmann}
\end{figure}

A motivation for BASSO's adaptive sampling scheme is the effectiveness of a sequence of Boltzmann distributions, as in AAS.
Figure~\ref{fig:boltzmann} illustrates the concept of AAS and BASSO on a one-dimensional objective function (shown in the top panel) to be minimized. The adaptive sampling procedure in AAS samples from a Boltzmann distribution where the temperature parameter $T$ is lowered when an improving point is found. The middle panel in Figure~\ref{fig:boltzmann} illustrates how the Boltzmann densities increasingly focus on the optimum as temperature decreases.

In BASSO, the adaptive subregion probability parameter 
%of selecting a subregion aids in weighing promising regions, and impacting  the sampling distribution. The adaptive probability parameter 
$\tilde{p}$ is a vector that captures the probability of selecting a subregion for further sampling. The bottom panel in Figure~\ref{fig:boltzmann} illustrates how adapting $\tilde{p}$ can focus the sampling distribution on the promising subregions. In the example, the domain is partitioned into three subregions, and initially $\tilde{p}$ is roughly $[1/3, 1/3, 1/3]$, representing a uniform distribution.  
 When the adaptive subregion probability parameter is updated to $\tilde{p}=[0.2, 0.6, 0.2]$, and then to $\tilde{p}=[0.1, 0.8, 0.1]$, there is a higher probability of sampling from the middle subregion with the global optimum. We draw an analogy to sampling from a Boltzmann distribution parameterized by temperature, to sampling from subregions parameterized by $\tilde{p}$.  
 %If $\tilde{p}$ satisfies Assumption 1, 
 %including \eqref{ratioassumption}, 
 As BASSO progresses,  $\tilde{p}$ is adapted to   focus on the subregion that contains the optimum.  
In the analysis of BASSO, Assumption 1 gives a condition on how  $\tilde{p}$ is updated that essentially ensures that  $\tilde{p}$ is consistent with the contours of the objective function and does not focus on subregions that do not contain the global optimum. In practice, this condition is difficult to always satisfy, as we show in the numerical experiments.

\subsection{BASSO Framework}
The BASSO framework subdivides the feasible region into smaller subregions, using sequential partitioning. Let $\Sigma^C_{k}$  be the collection of subregions that are produced by branching on the $k$th iteration, $\Sigma_k^C = \{\sigma_1, \ldots, \sigma_{m_k}\}$, where $m_k$ is the number of subregions on the $k$th iteration.  We assume that the branching strategy is such that the union of all subregions equals $S,$ and that the subregions are mutually exclusive. The subregions do not need to be of equal size.
We define the size of subregion $\sigma_i$, denoted $v_{\sigma_i}$, as the $d$-dimensional volume of $\sigma_i$ when $S$ is real-valued, and a count of points in $\sigma_i$ when $S$ is finite.
 We slightly abuse notation when the domain is finite or has mixed integer- and real-valued variables; technical details can be found in \cite{Zabinsky1995} and  \cite{Wood2001}.
 
We let the adaptive subregion probability $\tilde{p}_i(y)$
 %$\tilde{p}_i(\tilde{y}_k^*)$ 
 be the probability that subregion $\sigma_i$
%for $\sigma_i  \in \Sigma^C_{k}, 
is selected given an objective function value $y$, $y_* < y \leq y^*$.
%the incumbent objective function value is $\tilde{y}_k^*$ on iteration $k$. 
 BASSO keeps track of the best objective function value observed in each subregion $\sigma_i$, denoted  ${y}_i^*$, and 
  the overall best objective function value, denoted $\tilde{y}_k^*$ on iteration $k$, which we call the incumbent value. BASSO uses this incumbent value  to parameterize the adaptive subregion probabilities, i.e., $\tilde{p}_i(\tilde{y}_k^*)$.  In the analysis, we do not limit the parameter $y$ in $\tilde{p}_i(y)$ to be the incumbent value, because we want to analyze the impact of the parameter on the adaptive subregion probability.  It becomes important in establishing a condition in Assumption~1 on the adaptive subregion probabilities to represent the promising subregions.
 %, i.e., $\tilde{p}_i(\tilde{y}_k^*)$.  

% A new sample point on the $k$th iteration is generated with the adaptive sampling scheme. 
%The adaptive sampling scheme in BASSO  to generate a new sample point on the $k$th iteration has three main components: 
%\begin{enumerate}
%\item[(i)]  Branch a collection of subregions in $\Sigma^C_{k}$. 
%\item[(ii)]  Evaluate the adaptive subregion probabilities 
%  Select a subregion using the adaptive subregion probabilities 
%$\tilde{p}_i(\tilde{y}_k^*)$ for subregions $\sigma_i$ in $\Sigma^C_{k}$
%to select a promising subregion for sampling. 
%\item[(iii)] Generate a point within the selected subregion according to a surrogate model or the uniform distribution.
%\end{enumerate}
 The steps of the framework follow.\\ 
% (with pseudo-code in Appendix~\ref{sec:BASSO Pseudocode})

%\textcolor{red}{I think that Fig. 1 would fit best here.  I (Zelda) will move the figure and discussion, and rework the beginning of this section.}

 \noindent
 \textbf{Branching Adaptive Surrogate Search Optimization (BASSO). }
 
 \textbf{Step 0: Initialize.}  
Set $\sigma_1=S$, $\Sigma^C_{0}=\{\sigma_1\}$, $m_0 = 1$, and iteration counter $k =0$. 
Sample a point $X$  from the uniform distribution over $S$. Evaluate $f(X)$ and set 
${y}^*_1 =  f(X)$. Also, set  $\tilde{y}_0^*=f(X)$ to keep track of the incumbent function value observed. Set $\tilde{p}_1(\tilde{y}_0^*) = 1$.

\textbf{Step 1: Sample a new point.
%according to adaptive sampling probabilities.
}
Choose a subregion $\sigma_i$ from $\Sigma_{k}^C$ using the adaptive subregion probabilities $\tilde{p}_i(\tilde{y}_k^*)$ for $i = 1,\ldots,m_{k}$.  
%Generate a point $X_{k}$ on the selected subregion and set $Y_{k}=f(X_{k})$.
Generate a point $X$ on the selected subregion either uniformly or using surrogate modeling, and evaluate $f(X)$. Update the best objective function value $y_i^*$ in the selected subregion $\sigma_i$.
 
If the new objective function value $f(X)$ is less than the incumbent value $\tilde{y}^*_{k} $, 
%i.e., if $y_i^* < \tilde{y}^*_{k}$, 
update the incumbent function value 
%$Y_{k+1}=\tilde{y}^*_{k+1} = f(X)$ 
$\tilde{y}^*_{k+1} = f(X)$ and proceed to Step~2. Otherwise, repeat Step~1 for a maximum number of function evaluations before proceeding to Step~2.

\textbf{Step 2: Branch subregions.}
Branch a collection of  subregions in $\Sigma_k^C$ according to a branching strategy. Update  the new number of subregions $m_{k+1}$, and update the collection of subregions, $\Sigma_{k+1}^C = \{\sigma_1,\ldots, \sigma_{m_{k+1}}\}$. 

\textbf{Step 3: Update adaptive subregion probabilities.}
Update $\tilde{p}_i(\tilde{y}_{k+1}^*)$ for $i=1, \ldots, m_{k+1}$.
%$A$, $B$ or $C$.\\

\textbf{Step 4: Evaluate stopping criterion.} If a stopping criterion is met, stop.  Otherwise increment the iteration counter $k\leftarrow k+1$ and go to Step 1.\\

\begin{figure}[htbp!]
     \centering
     \begin{subfigure}[b]{0.45\textwidth}
         \centering
         \includegraphics[width=\textwidth]{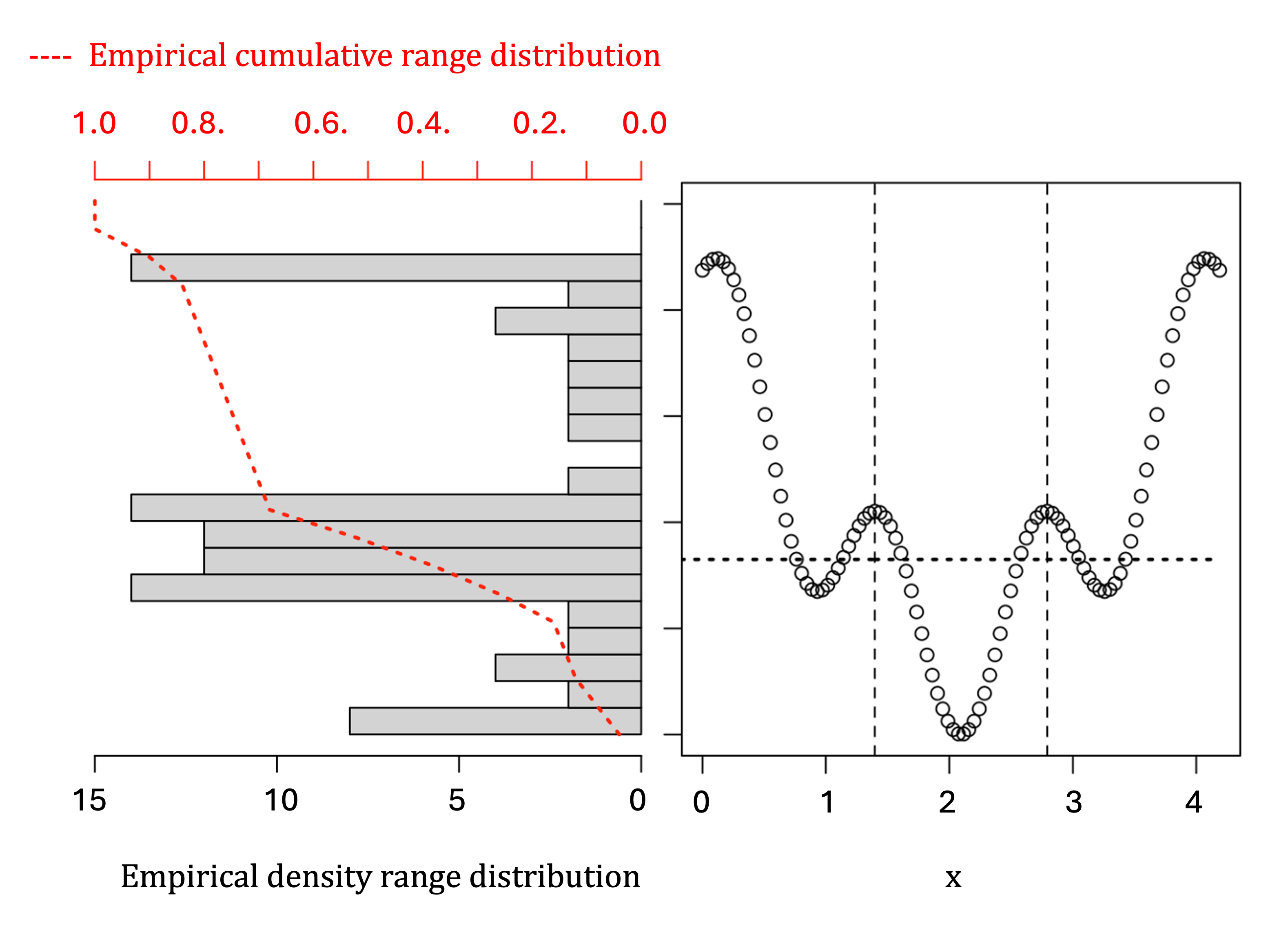}
         \caption{Uniform sampling density}
         \label{fig:rangedist1}
     \end{subfigure}
     \begin{subfigure}[b]{0.45\textwidth}
         \centering
         \includegraphics[width=\textwidth]{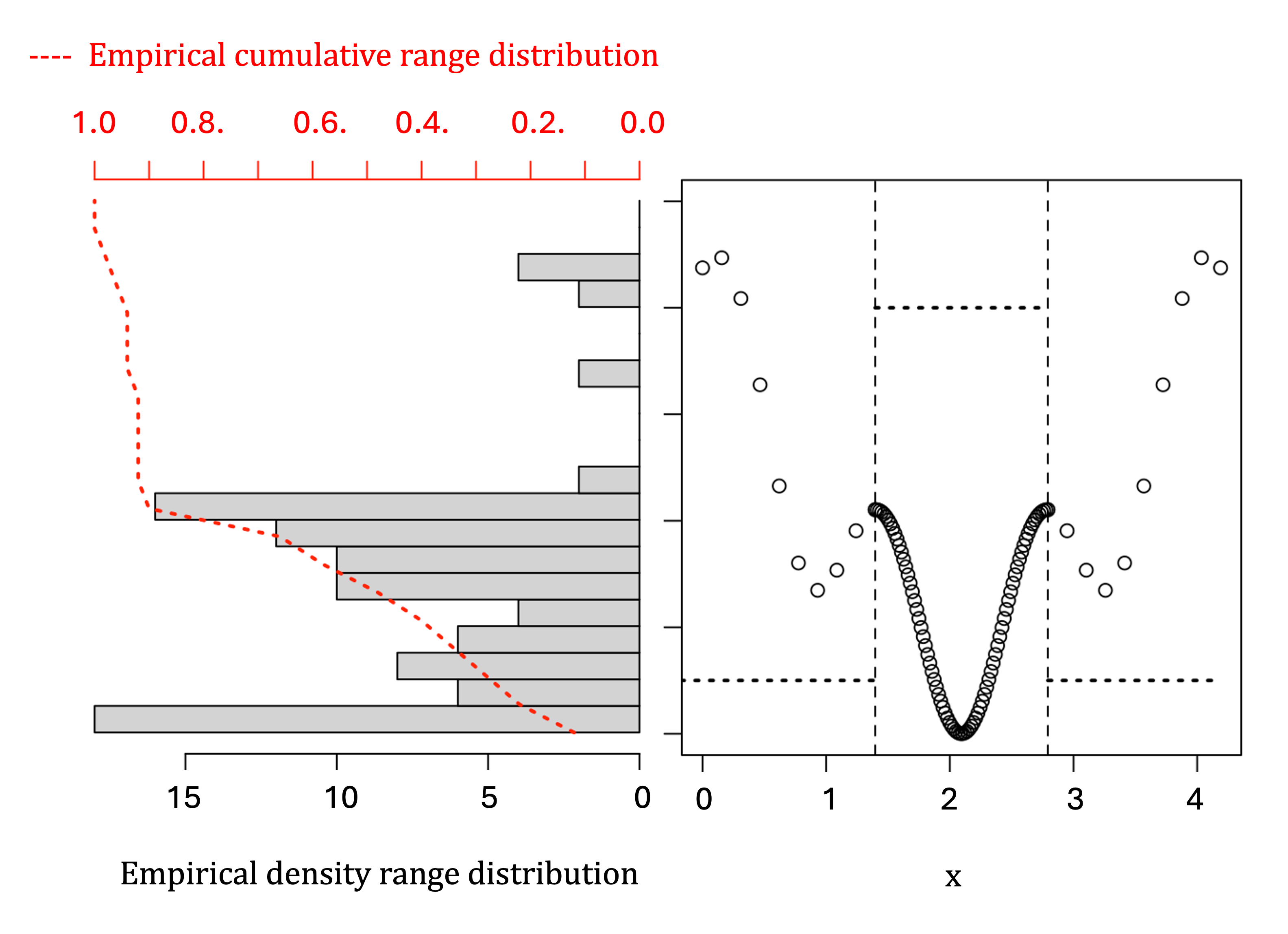}
         \caption{Adaptive sampling density}
         \label{fig:rangedist2}
     \end{subfigure}
        \caption{Empirical density range distribution (gray histogram) and empirical cumulative range distribution (dashed red line) for  100 points sampled independently, where the domain is partitioned into three subregions (i.e., intervals). In (a), the adaptive subregion probability is $\tilde{p}=[0.33,0.33,0.33]$, representing a uniform distribution.  
In (b),  the adaptive subregion probability $\tilde{p}=[0.1,0.8,0.1]$.  Once a subregion is selected, a point is sampled uniformly on that interval.}
        \label{fig:rangedist}
\end{figure}

Figure~\ref{fig:rangedist} illustrates the impact of the adaptive sampling distribution on sampling points with good objective function values. In Figure~\ref{fig:rangedist}, the domain $S$ is partitioned into three subregions, and 100 points are generated independently.  For each point, a subregion is selected using the adaptive subregion probabilities, and then a point is generated uniformly within the selected subregion.  
In  Figure~\ref{fig:rangedist1},
 the  
 adaptive subregion probabilities are $\left[1/3, 1/3, 1/3\right]$, whereas the adaptive subregion probabilities in  Figure~\ref{fig:rangedist2} are  $\left[0.1, 0.8, 0.1\right]$ for the three subregions.

The empirical density range distribution  is the number of objective function values that fall in an interval, as illustrated by the vertical histogram in Figure~\ref{fig:rangedist}.
 %The vertical histogram in Figure~\ref{fig:rangedist} shows the empirical density range distribution of 
 This provides a way to visualize the impact of the adaptive sampling distribution.  The red dashed curve illustrates the empirical cumulative range distribution, which is the probability the objective function value is below $y$ when a point is sampled according to a sampling distribution. Notice that the range distribution in Figure~\ref{fig:rangedist2} puts more weight on the lower objective function values. The impact of adaptive sampling on the range distribution is important for the finite-time analysis of BASSO.

 \subsection{BASSO Variations}
 \label{sec:BASSOvars}

%{\red{
To implement BASSO, it is important to specify: %\begin{enumerate*}
%\item[(i)] 
(i) the branching scheme, 
%\item[(ii)]   
(ii) the adaptive subregion probabilities 
$\tilde{p}_i(\tilde{y}_k^*)$ for subregions $\sigma_i$ in $\Sigma^C_{k}$,
and 
%\item[(iii)] 
(iii) the method for generating a point within the selected subregion (surrogate modeling or the uniform distribution).
%\end{enumerate*} 
%}}

%\usepackage[inline]{enumitem} and \begin{enumerate*}\end{enumerate*} (with asterisks) should do your work.

The implementation variations demonstrate the interplay between branching strategy, the adaptive subregion probabilities, and the surrogate model used to generate a point on the selected subregion. The variations follow.
Additional implementation details for the variations used in the numerical experiments are in Section~\ref{sec:expts}.
\\

\noindent \textbf{(i) Branching strategy}

\noindent
Create a collection of subregions to branch.
%, by ranking the subregions by their best observed function values. 
Add   $\lceil 0.1m_k\rceil$ subregions (where $m_k$ is the number of subregions in $\Sigma^C_k$  on the $k$th iteration) ranked by their best observed function values to the collection, and add additional  $\lceil 0.1m_k\rceil$ of the {\bf remaining} subregions  ranked by  the largest volume. We branch two types of subregions, those with good observed objective function values and those with large volumes, as an attempt to balance exploitation with exploration.

We subdivide the subregions in the collection  coordinate-wise by dividing the largest dimension in half.
If any newly created subregion has 
fewer than two observations in it,
%less than two observations in it,
sample one or two points uniformly in the subregion so that each subregion has at least two observed points. \\

\noindent \textbf{(ii) Adaptive subregion probabilities, $\tilde{p}$}

\noindent
Four variations of adaptive subregion probabilities are implemented.The first variation (a) is greedy in the sense that it exploits the subregions with the best observed function value. The second variation (b) favors exploration as it uses an estimate of the variance in the probability of selecting subregions.  The third variation (c) attempts to balance exploration and exploitation by estimating the distribution of the range on a subregion with a one-dimensional Gaussian process. The fourth variation (d) relaxes the greedy variation by estimating confidence intervals around the incumbent value and each subregion's best observed function value. \\

\begin{itemize}
\item[\textbf{a: Observed best function value.}]
Choose  subregion $\sigma_i$ from the subregions in $\Sigma^C_k$  on the $k$th iteration with probability  $\tilde{p}_i(\tilde{y}^*_k)$, where
\begin{equation}
   \tilde{p}_i(\tilde{y}^*_k)= 
  \left(\left(y^*_i-\tilde{y}^*_k+1\right)
   \sum_{j=1, \ldots, m_k}\ \left(\frac{1}{y_j^*-\tilde{y}^*_k+1}\right)
   \right)^{-1}. 
\label{eq:ptildeincumbent}
\end{equation}
%{\red{Does not seem to satisfy Assumption 1.}}
%{\blue{I will try with an example.  THANK YOU, Pete!
%Consider 3 subregions with $y_i^*$ = 10, 5, 1, and incumbent $\tilde{y}^*_k=1$.
%Then $\mathcolor{orange}{\sum_{j=1,2,3}\ \left(\frac{1}{y_j^*-\tilde{y}^*_k+1}\right) = \frac{1}{10-1+1} + \frac{1}{5-1+1}+\frac{1}{1-1+1} = \frac{1}{10} + \frac{1}{5}+\frac{1}{1} = \frac{13}{10}}$.
%\begin{align*}
%    \tilde{p}_1(\tilde{y}^*_k)&= \left( (y^*_1-\tilde{y}^*_k+1)\mathcolor{orange}{\sum_{j=1,2,3}\ \left(\frac{1}{y_j^*-\tilde{y}^*_k+1}\right)}\right)^{-1}\\
%    &=\left( (10-1+1)\mathcolor{orange}{\frac{13}{10}}\right)^{-1}\\
%    &=1/13
%\end{align*}

%\begin{align*}
%    \tilde{p}_2(\tilde{y}^*_k)&= \left( (y^*_1-\tilde{y}^*_k+1)\mathcolor{orange}{\sum_{j=1,2,3}\ \left(\frac{1}{y_j^*-\tilde{y}^*_k+1}\right)}\right)^{-1}\\
%    &=\left( (5-1+1)\mathcolor{orange}{\frac{13}{10}}\right)^{-1}\\
%    &=2/13
%\end{align*}

%\begin{align*}
%    \tilde{p}_3(\tilde{y}^*_k)&= \left( (y^*_1-\tilde{y}^*_k+1)\mathcolor{orange}{\sum_{j=1,2,3}\ \left(\frac{1}{y_j^*-\tilde{y}^*_k+1}\right)}\right)^{-1}\\
%    &=\left( (1-1+1)\mathcolor{orange}{\frac{13}{10}}\right)^{-1}\\
%    &=10/13
%\end{align*}

\item[\textbf{b: Sample variance.}]
Choose  subregion $\sigma_i$ from the subregions in $\Sigma^C_k$  on the $k$th iteration with probability  $\tilde{p}_i$, where
\begin{equation}
   \tilde{p}_i= \left\{
   \begin{array}{l l}
     v_{\sigma_i}/v_{\Sigma^C_{k}} & {\rm if} \  k =1\\
   \ \  & \\
   \frac{\left(s_{k}(\sigma_i)\right)^2}{\sum_{j=1, \ldots, m_k 
   %\vert \vert\Sigma^C_{k}\vert \vert
   }\left(s_{k}(\sigma_i)\right)^2}\ 
    & {\rm if} \  k > 1,
   \end{array}
   \right.  \label{eq:ptildesamplevariance}
\end{equation}
and 
$(s_{k}
(\sigma_i))^2$ is the sample variance in  subregion $\sigma_i$, that is,
\begin{equation}(s_{k}(\sigma_i))^2 = 
\left(\frac{1}{N^{i}_k-1}\right) \sum_{j=1}^{N^i_k} (y_j - \overline{f_i})^2,
\label{eq:samplevariance}
\end{equation}
and $y_j$ is the function value of the $j$th previously sampled point in subregion $\sigma_i$,
and
$\overline{f_i}$ is the sample mean of the $N^i_k$ observed objective function values in subregion $\sigma_i$, i.e.,
$\overline{f_i}=(1/N^i_k)\sum_{j=1,\ldots,N^i_k} y_j$.
\\
%{\red{Does not seem to satisfy Assumption 1.}}

%\end{itemize}
%For use in Step~5, set $s^{max}_{iter}(\sigma_i)=s_{iter}(\sigma_i)$.
%Then, uniformly generate a sample point within the chosen subregion $\sigma_i$.
%Update the number of points that have been sampled in $\sigma_i$ as  $N^i$. 

% \item[\textbf{c. Gaussian process on the domain}]
% We use, with the overall incumbent function value $\tilde{y}_{k}^*$ as its argument, that is,
% \[
% ACQ_i(\tilde{y}_{k}^*) = \int_{x \in \sigma_i}\left(\frac{\tilde{y}^*_{k}- \hat{g}_{i,k}(x)}{\hat{s}_{i,k}(x)}\right) 
% + \hat{s}_{i,k}(x)\Phi\left(
% \left(\frac{\tilde{y}^*_{k} - \hat{g}_{i,k}(x)}{\hat{s}_{i,k}(x)}\right) 
% \right)dx,
% \]
% which uses the distance to the incumbent objective function value combined with uncertainty.
% Then, the probability of selecting subregion $\sigma_i$ is,
% \begin{equation}
% \label{eq:piB}
% \tilde{p}_i(\tilde{y}_{k}^*) = \frac{ACQ_i(\tilde{y}^*_{k})}{\sum_{i=1}^{m_{k}} ACQ_i(\tilde{y}^*_{k})}.
% \end{equation}
% \end{itemize}

\item[\textbf{c. Gaussian process on the range distribution.}]
For each of the subregions $\sigma_i$ in $\Sigma^C_{k}$, $i=1, \ldots, m_{k}$, we construct a Gaussian process 
 on the empirical cumulative range distribution using the objective function values of the $N^i_k$ points that have been sampled in $\sigma_i$ on the $k$th iteration. 
Specifically, we
use the $N^i_k$ observed objective function values in subregion $\sigma_i$, $y_j$, $j=1, \ldots, N^i_k$ and  the associated empirical cumulative function values $F(y_j)$, using
%(e.g., cumulative histogram for each $y_{j}$),
$\mathcal{D}= \{ {y}_j, F({y}_j) \}^{N^i_k}_{j=1}$ 
as the input space. 
The empirical cumulative function  $F(y_j)$ is the count of observed function values less than or equal to $y_j$, divided by $N^i_k$.

The Gaussian process on the range distribution provides a predictor function on subregion $\sigma_i$ as an  approximation of  the cumulative range distribution on $\sigma_i$, where the mean of the Gaussian process is denoted  $\hat{F}_{i,k}(y)$.
Note that the Gaussian process predictor is a one-dimensional function (as illustrated in Figure~\ref{fig:rangedist}), and is easily computed.

The adaptive subregion probability of selecting subregion $i$ is,
\begin{equation}
\label{eq:piC}
\tilde{p}_i(\tilde{y}^*_k + \delta) = \frac{\hat{F}_{i,k}(\tilde{y}^*_{k}+\delta)}{\Sigma_{i=1}^{m_{k} }\hat{F}_{i,k}(\tilde{y}^*_{k}+\delta)}.    
\end{equation}
where $\tilde{y}^*_{k}+\delta$ is the value of the fifth lowest observed objective function value. \\

%{\red{Seems to SATISFY Assumption 1.}}

\item[\textbf{d. Confidence bounds.}]
\label{eq:pid}
Choose  subregion $\sigma_i$ from the subregions in $\Sigma^C_k$ on the $k$th iteration  with probability  $\tilde{p}_i(\tilde{y}_k^*)$, where

\begin{equation}
   \tilde{p}_i(\tilde{y}^*_k)=  \left\{
   \begin{array}{l l}
    0  & {\rm if} \  LB_i \geq 
    \widetilde{UB}_k  
   \\
   \ \  & \\
   \frac{\widetilde{UB}_k - LB_i}{\sum_{i =1,\ldots,m_k}\widetilde{UB}_k - LB_i}\ \ 
    & {\rm otherwise} \  
   \end{array}
   \right.  \label{eq:ptildeUBLB}
\end{equation}
and $LB_i=y^*_i - s_k(\sigma_i)$, $\widetilde{UB}_k=\tilde{y}^*_k+s_k(\tilde{\sigma})$, $\tilde{\sigma}$ is the subregion containing a point with the incumbent value$\tilde{y}^*_k$, and $s_k(\cdot)$ is the square root of the sample variance as in \eqref{eq:samplevariance}.
%{\red{Seems to SATISFY Assumption 1.}}
\end{itemize}

\vspace{0.25in}
\noindent \textbf{(iii) Generate a point on the selected subregion}

\noindent
Three variations of generating a point within a subregion are implemented.
The first variation (A) has no surrogate model and simply samples uniformly on the the selected subregion.  The second variation (B) uses a Gaussian process as a surrogate model on the selected subregion, and the third variation (C) uses a regularized quadratic regression as the surrogate model.\\

%{\red{Why doesn't A, B, C indent like above for a, b, c, d?}}
\begin{itemize}
\item[\textbf{A: Uniform.}]
Generate a point $X$ uniformly distributed on the selected subregion.\\

\item[\textbf{B: Maximize expected improvement for a Gaussian process on the domain.}]
For each of the subregions $\sigma_i$ in $\Sigma^C_{k}$, $i=1, \ldots, m_{k}$, we construct a Gaussian process 
 on the domain
% For the selected subregion $\sigma_i$ in $\Sigma^C_{k}$, $i=1, \ldots, m_{k}$, construct a Gaussian process 
 using the $N^i_k$ points that have been sampled in $\sigma_i$ on the $k$th iteration and their associated function values, using $\mathcal{D}= \{ {x}_j, f({x}_j) \}^{N^i_k}_{j=1}$ as the input space.
%The Gaussian process provides the predictor function, {\red{mean function,}} denoted   $\hat{g}_{i,k}(x)$, and model variance,
%denoted
 $\hat{s}^2_{i,k}(x)$,
 for $x\in \sigma_i$. 
The estimated mean function of the Gaussian process that approximates the objective function value using the sample points in  subregion $\sigma_i$ on iteration $k$ is denoted $\hat{g}_{i,k}(x)$ and model variance is
denoted
 $\hat{s}^2_{i,k}(x)$,
 for $x\in \sigma_i$. 
 
The expected improvement function for subregion $\sigma_i$
uses the best objective function value observed in that subregion $y^*_{i}$ to capture a better point, and a 
 standard Normal cumulative  distribution $\Phi$ to capture uncertainty,
\begin{equation}
EI_i(x) =\left(\frac{{y}^*_{i} - \hat{g}_{i,k}(x)}{\hat{s}_{i,k}(x)}\right) 
+ \hat{s}_{i,k}(x)\ \Phi
\left(\frac{{y}^*_{i} - \hat{g}_{i,k}(x)}{\hat{s}_{i,k}(x)} 
\right).
\label{eq:EIx}
\end{equation}

Generate a point $X$ in subregion $\sigma_i$ 
%in Step~1, 
that maximizes the expected improvement function over $x\in \sigma_i$.  
The maximization of $EI_i(x)$ is implemented numerically by performing a grid search on the subregion. This implementation of Variation B is computationally intensive for high dimensions.  In the numerical experiments, we executed Variation B on problems in dimensions 20 and 50, and not on the higher dimensions.
\\

\item[\textbf{C: Minimize a regularized quadratic regression on the domain.}]
For each of the subregions $\sigma_i$ in $\Sigma^C_{k}$, $i=1, \ldots, m_{k}$, we construct a regularized quadratic regression 
 on the domain
%
%For the  selected subregion $\sigma_i$ in $\Sigma^C_{k}$, $i=1, \ldots, m_{k}$, construct a regularized quadratic regression 
using the  $N^i_k$ points that have been sampled in $\sigma_i$ on the $k$th iteration and their associated function values, using $\mathcal{D}= \{ {x}_j, f({x}_j) \}^{N^i_k}_{j=1}$ as the input space.
  The regularized quadratic regression function for $x\in \sigma_i$ is 
 \begin{equation}\label{eqn:fquadZZ}
\hat{q}^i_k(x) = \hat{\beta}_0 + \sum_{\ell=1}^{d} \hat{\beta}_\ell x_{\ell} + \sum_{\ell=1}^{d} \sum_{m=\ell}^{d} \hat{\beta}_{\ell m} x_\ell x_m
\end{equation}

\noindent 
  where the $\ell$th component (decision variable) of a  point $x$ is denoted  $x_{\ell}$.
The  coefficients of the quadratic regression ($\hat{\beta}_0$ is the intercept term,
$\hat{\beta}_\ell$ are the coefficients of the linear terms, and $\hat{\beta}_{\ell m}$ are the coefficients of the quadratic terms, $\ell, m=1,\ldots,d$)
are determined by minimizing the loss function 
given in \eqref{eqn:qreglossZZ}.

The loss function
uses the $N^i_k$ sampled points $x_j \in \sigma_i$ and their objective function values $y_j$, for $j=1,\ldots,N^i_k$, as 
%An additional subscript is used to  denote the $\ell$th component (decision variable) of a sample point, i.e., $x_{j,\ell}$. 

\begin{equation}\label{eqn:qreglossZZ}
\mathcal{L}(\hat{\beta}) = \frac{1}{2N^i_k}  \sum_{j=1}^{N^i_k} (\hat{q}^i_k (x_{j})- y_j)^2 + \lambda \left( \sum_{l=1}^{d} \lvert\hat{\beta_\ell}\rvert + \sum_{\ell=1}^{d} \sum_{m=\ell}^{d} \lvert\hat{\beta}_{\ell m}\rvert \right),
\end{equation}
where
$\hat{\beta}$ is the vector of coefficients, including $\hat{\beta}_0$, $\hat{\beta}_\ell$, and $\hat{\beta}_{\ell m}$.
The first term in the loss function is the squared difference between the quadratic model and the observed value.  The second term in the loss function is a penalty term 
weighted by a 
hyperparameter $\lambda$, that is a regularization parameter controlling the amount of regularization applied.  
%In the numerical experiments, $\lambda=1$.
In the numerical experiments,  we update $\lambda$ based on observations. For subregions with more than 50 points, we determine the regularizer $\lambda$ by performing a 5-fold cross-validation on that subregion and use the value of $\lambda$ that minimizes the cross-validation prediction error.  For subregions with 50 or fewer points, we set $\lambda=1$.

Generate a point $X$ in subregion $\sigma_i$ that minimizes the regularized quadratic regression $\hat{q}^i_{k}(x)$ over $x\in \sigma_i$.
In the numerical experiments, we use a trust region algorithm 
\cite{trustregion2000}
to approximately minimize $\hat{q}^i_{k}(x)$ over $\sigma_i$.
\end{itemize}

%\vspace{0.2in}
%In the numerical experiments, we explore the performance of 12 BASSO variations: Aa, Ab, Ac, Ad, Ba, Bb, Bc, Bd, Ca, Cb, Cc, and Cd.

%Generate a sample point in the chosen subregion, $x\in \sigma_i$, that 
%minimizes 
%where  $\hat{q}^i_{k}(x)$ is a regularized quadratic regression over $\sigma_i$, as in 
%\eqref{eqn:fquad}, using the coefficients from minimizing the loss function with an additional penalty term, as in \eqref{eqn:qregloss}.

\section{Analysis}
\label{sec:analysis}
The BASSO analysis begins by assuming that sample points are generated according to a uniform distribution on selected subregions. %sampling within a subregion is performed according to a uniform distribution.  
We show in Theorem~\ref{Theorem 1} that the objective function values of BASSO with uniform sampling within a subregion stochastically dominate those of a special case of HAS, HAS\textprime, under Assumption~1 on the adaptive subregion probabilities $\tilde{p}_i$ (Appendix~\ref{sec:HAS} provides an overview of HAS). We then relax the assumption of uniform sampling within the subregion and analyze the case of surrogate-driven sampling. % allow surrogate sampling on a subregion, and show in 
Theorem~\ref{Theorem 2} shows that the objective function values of BASSO with surrogate sampling stochastically dominate those of BASSO with uniform sampling on a subregion, under Assumption~2 on the surrogate model. 
   
Using prior analyses of PAS and HAS~\cite{Zabinsky1992,Zabinsky1995,Bulger1998,Wood2001}, we then bound the expected number of BASSO function evaluations needed to achieve an $\epsilon$-optimal solution.  Corollary~\ref{corollary1} provides the desired linearity result for functions satisfying the Lipschitz condition, and Corollary~\ref{corollary2} provides the desired linearity result for functions on a finite domain, specifically a $d$-dimensional integer lattice. 
 
%The special case of HAS, i.e., 
%The special case 
HAS\textprime \ 
 uses the uniform distribution as its underlying  sampling distribution, hence, given a value $z$, the \textit{bettering probability}, i.e., the probability to sample a location with value below $z$, can be specified as:%is: %.  
%$\delta$ 
%The bettering probability of HAS\textprime \  is specified as
\begin{equation}
\label{eq:betteringprob}
b(z)=\sum_{i=1,\ldots,m_k} \frac{\tilde{p}_i(z)}{v_{\sigma_i}}
v_{\sigma_i}(z)
\end{equation}
where $y_* < z \leq y^*$, 
given  subregions $\sigma_i \in \Sigma^C_k$,  $m_k$, and probabilities $\tilde{p}_i(z)$ for $i=1, \ldots, m_k$ with $y_* < z \leq y^*$.
As earlier, 
 $v_{\sigma_i}$ is the size of subregion  $\sigma_i$ 
 and 
 $v_{\sigma_i}(z)$ is the size of the set $\{x \in \sigma_i: f(x) \leq z\}$ 
for $y_* < z \leq y^*$ 
 (i.e.,  size is $d$-dimensional volume when $S$ is real-valued, or a count when $S$ is finite). 
 The bettering probability $b(z)$ depends on the adaptive subregion probability parameter
$\tilde{p}_i(z)$ for $i=1, \ldots, m_k$ with $y_* < z \leq y^*$, which drives the 
 the sampling distribution.

 Let the random variable %$Y_k^{BASSO{\rm unif}}$ 
 $Y_k^{{\rm \it BASSOunif}}$ be the incumbent function value $\tilde{y}_k^* $
 %from $BASSOunif$ 
 on the $k$th iteration when BASSO generates a point uniformly distributed within a selected subregion.  
 %We start by assuming that, once a subregion is selected with the adaptive subregion probabilities, the method of generating a sample point in the subregion uses a uniform distribution. 
Let $Y_k^{HAS'}$ be the incumbent function value from HAS\textprime \   on the $k$th iteration. 

The parameter $z$ in $\tilde{p}_i(z)$ is the incumbent function value in BASSO.  However, in the  analysis,  we do not limit the parameter $y$ in $\tilde{p}_i(y)$ to be the incumbent value because we need to characterize the conditional probability of improvement as the parameter $y$ in $\tilde{p}_i(y)$  changes.
This is important in Assumption 1 where we need to account for the amount of improvement as the incumbent value improves. 
%
%we do not limit the parameter $y$ in $\tilde{p}_i(y)$ to be the incumbent value, because we want to analyze the impact of the parameter on the adaptive subregion probability.  This is important in Assumption 1 where we want to characterize the conditional probability of improvement as the parameter $y$ in $\tilde{p}_i(y)$  changes.

We note that the probability of generating a point with objective function value no greater than $y$  for BASSO, given the incumbent  value $z$,
%$z=\tilde{y}^*_k$,  
adaptive subregion probabilities $\tilde{p}_i(z)$, when the uniform distribution is used on the selected subregion, can be expressed as
\begin{align*}
P(Y^{{\rm \it BASSOunif}}_{k+1}\leq y |Y^{{\rm \it BASSOunif}}_{k}=z) &= \sum_{i=1,\ldots,m_k} \frac{\tilde{p}_i(z)}{v_{\sigma_i}}
v_{\sigma_i}(y)
\end{align*}
for $y_* <y \leq z \leq y^*$.
An understanding of this expression is that each subregion has a probability $\tilde{p}_i(z)$ of being selected, and once selected, a point is sampled uniformly on the subregion with a probability of achieving a function value less than or equal to $y$ (i.e.,  $v_{\sigma_i}(y)/v_{\sigma_i}$).
Another interpretation of this expression is viewing the sampling density as $\tilde{p}_i(z)/v_{\sigma_i}$.  

The challenge in partitioning and using an adaptive subregion probability is focusing on appropriate subregions. 
We make the following assumptions regarding the adaptive subregion parameters
$\tilde{p}_i(z)$ for $i=1, \ldots, m_k$ with $y_* < z \leq y^*$.\\

\noindent{\bf Assumption 1:  Adaptive subregion probabilities $\tilde{p}_i(z)$}
\begin{enumerate}
    \item  $0 \leq \tilde{p}_i(z)\leq 1$ for $i=1, \ldots, m_k$ with $y_* < z \leq y^*$
    \item $\sum_{i=1,\ldots,m_k} \tilde{p}_i(z) = 1$
    \item Given two values $t_z$ and $t_{z'}$ such that $y_* < t_z \leq t_{z'}\leq y^*$, the following holds, %ratio is non-increasing in $t_z$, that is,
\begin{equation}\label{ratioassumption} \frac{\sum_{i=1,\ldots,m_k} \left(\frac{\tilde{p}_i(t_z)}{v_{\sigma_i}} \right)
{v_{\sigma_i}(y)}}{\sum_{i=1,\ldots,m_k} \left(\frac{\tilde{p}_i(t_z)}{v_{\sigma_i}} \right)
{v_{\sigma_i}(z)}} 
\geq
\frac{\sum_{i=1,\ldots,m_k} \left(\frac{\tilde{p}_i(t_{z'})}{v_{\sigma_i}} \right)
{v_{\sigma_i}(y)}}{\sum_{i=1,\ldots,m_k} \left(\frac{\tilde{p}_i(t_{z'})}{v_{\sigma_i}} \right)
{v_{\sigma_i}(z)}}
\end{equation} 
for $y_* < y \leq z\leq y^*$.
\end{enumerate}

\vspace{.2in}
The first two assumptions ensure that the $\tilde{p}_i(z)$ form a proper probability. The third assumption 
is that the ratio (in the left-hand-side of \eqref{ratioassumption}) is non-increasing in $t_z$.
This
ensures that $\tilde{p}_i(z)$ focuses the sampling distribution on the regions with good objective function values. The ratio with $t_z=z$ can be interpreted as the conditional probability,
\begin{align*}
\label{eq:condprob}
&P(Y^{{\rm \it BASSOunif}}_{k+1}\leq y |Y^{{\rm \it BASSOunif}}_{k+1}\leq z,
Y^{{\rm \it BASSOunif}}_k=z
%\tilde{p}_i(z)
) \\
&\qquad \qquad = \frac{\sum_{i=1,\ldots,m_k} \left(\frac{\tilde{p}_i(z)}{v_{\sigma_i}} \right)
{v_{\sigma_i}(y)}}{\sum_{i=1,\ldots,m_k} \left(\frac{\tilde{p}_i(z)}{v_{\sigma_i}} \right)
{v_{\sigma_i}(z)}}    
\end{align*} 
for $y_* < y \leq z\leq y^*$, when uniform sampling is used on the selected subregion.  
The conditional probability captures, not only the probability of improving, but the amount of improvement to a function value $y$ where $y\leq z$.
Thus, the sampling distribution parameterized by $\tilde{p}_i(z)$ can improve the conditional probability by increasing the chance of selecting ``good'' subregions.
This is key to a good branching algorithm, although, as illustrated in the numerical experiments, it is challenging to devise a practical implementation of $\tilde{p}_i(z)$ that satisfies this assumption.

We now show that the objective function values of BASSO with uniform sampling within subregions stochastically dominate those of HAS\textprime.
\begin{theorem}\label{Theorem 1} 
Given the optimization problem \eqref{eq:optproblem} over a continuous/finite domain $S$ with Assumption~1 holding, the  function values of BASSO with uniform sampling within subregion stochastically dominate those of $HAS'$, i.e., 
\[P(Y^{{\rm \it BASSOunif}}_k \leq y) \geq P(Y^{HAS
'}_k \leq y)   \]
for  $k = 0,1,\ldots$,  and $y_* < y\leq y^*$.
\end{theorem}
\noindent \textbf{Proof:} 
%The proof of
%Theorem 1 relies on Lemma 30 in \cite{Shen2005}, which is included in Appendix~\ref{sec:Lemma30} for completeness.
%%
%To prove that $Y^{\rm \it BASSOunif}_k$ stochastically dominates $Y^{HAS'}_k$, for  $k = 0,1,2,\ldots$,  we  show that the two sequences 
%%$Y^{BASSOunif}_k$, for $k = 0,1,2,\ldots$ and $Y^{HAS'}_k$ for $k = 0,1,2,\ldots$ 
%satisfy the three conditions in the Lemma in Appendix~\ref{sec:Lemma30}.
%
To prove that $Y^{\rm \it BASSOunif}_k$ stochastically dominates $Y^{HAS'}_k$, for  $k = 0,1,2,\ldots$,  we  show that the two sequences 
%$Y^{BASSOunif}_k$, for $k = 0,1,2,\ldots$ and $Y^{HAS'}_k$ for $k = 0,1,2,\ldots$ 
satisfy three conditions as stated in Lemma 30 in \cite{Shen2005}, and included in Appendix~\ref{sec:Lemma30} for completeness.
At iteration $k$, we are given the subregions $\sigma_i$ in $\Sigma^C_{k}$ for $i=1,\ldots, m_k$, the incumbent value $Y_k=\tilde{y}^*_{k}$, and the adaptive subregion probabilities $\tilde{p}_i(\tilde{y}^*_{k})$. Let $z=\tilde{y}^*_{k}$ for ease of notation.
%that satisfy Assumptions (1)-(3).

To show Condition 1, we need to prove that  
\[
P(Y^{{\rm \it BASSOunif}}_{k+1}\leq y |Y^{{\rm \it BASSOunif}}_k\ =z) \geq  P(Y^{HAS'}_{k+1}\leq y |Y^{HAS'}_k\ =z)
\]
for any fixed $k$, 
$k = 0,1,2,\ldots$,  and $y$, $z$ values, for $y_* \leq y, z \leq y^*$. We are interested in the case when $y < z$, because, when
$y \geq z$, the conditional probabilities are both equal to 1, and so satisfy the inequality.

 For $y < z$, using conditional probabilities,
we have,
\begin{align*}
&P(Y^{{\rm \it BASSOunif}}_{k+1}\leq y |Y^{{\rm \it BASSOunif}}_k=z) 
\\
&=P(Y^{{\rm \it BASSOunif}}_{k+1}\leq y | Y^{{\rm \it BASSOunif}}_{k+1}\leq z,Y^{{\rm \it BASSOunif}}_k=z) \\
& \ \ \ \ \cdot
P(Y^{{\rm \it BASSOunif}}_{k+1}\leq z | Y^{{\rm \it BASSOunif}}_k=z)
\\
&= \frac{\sum_{i=1,\ldots,m_k} \frac{\tilde{p}_i(z)}{v_{\sigma_i}}
v_{\sigma_i}(y)}{\sum_{i=1,\ldots,m_k} \frac{\tilde{p}_i(z)}{v_{\sigma_i}}
v_{\sigma_i}(z)} \cdot \sum_{i=1,\ldots,m_k} \frac{\tilde{p}_i(z)}{v_{\sigma_i}}
v_{\sigma_i}(z)\\
\insertpar{and by  Equation \eqref{ratioassumption} 
in Assumption~1 with $z \leq y^*$, }
& \geq  \frac{\sum_{i=1,\ldots,m_k} \frac{\tilde{p}_i(y^*)}{v_{\sigma_i}}
v_{\sigma_i}(y)}{\sum_{i=1,\ldots,m_k} \frac{\tilde{p}_i(y^*)}{v_{\sigma_i}}
v_{\sigma_i}(z)} \cdot \sum_{i=1,\ldots,m_k} \frac{\tilde{p}_i(z)}{v_{\sigma_i}}
v_{\sigma_i}(z) \\
\insertpar{and because $\tilde{p}_i(y^*)$ represents a uniform distribution on $S$, i.e., $\tilde{p}_i(y^*)=v_{\sigma_i} / v_S$,}
&= \frac{v_S(y)}{v_S(z)}
\cdot \sum_{i=1,\ldots,m_k} \frac{\tilde{p}_i(z)}{v_{\sigma_i}}
v_{\sigma_i}(z) \\
&= \frac{v_S(y)}{v_S(z)}
\cdot b(z)
= P(Y^{HAS'}_{k+1}\leq y |Y^{HAS'}_{k}=z). 
\end{align*}

To show Condition 2, we need to prove that
$P(Y^{\rm \it BASSOunif}_{k+1}\leq y |Y^{\rm \it BASSOunif}_k =z)$ is non-increasing in $z$, that is, we need to prove 
\begin{small}
\begin{equation}
\label{eq:Condition2}
P(Y^{\rm \it BASSOunif}_{k+1}\leq y |Y^{\rm \it BASSOunif}_k=z)
\geq 
P(Y^{\rm \it BASSOunif}_{k+1}\leq y |Y^{\rm \it BASSOunif}_k=z') 
\end{equation}
\end{small}
for any fixed
$k$,
$k = 0,1,2,\ldots$, and any $y$ and $z$ values, $y_* < y \leq z \leq z' \leq y^*$. 

If $y \geq z$, 
$P(Y^{\rm \it BASSOunif}_{k+1}\leq y |Y^{\rm \it BASSOunif}_k\ =z)
=1$,
hence the inequality in \eqref{eq:Condition2} is satisfied.
For $y < z$,
 we have,
\begin{align*}
P(Y^{\rm \it BASSOunif}_{k+1}\leq y |Y^{\rm \it BASSOunif}_{k}=z) &= \sum_{i=1,\ldots,m_k} \frac{\tilde{p}_i(z)}{v_{\sigma_i}}
v_{\sigma_i}(y),
\end{align*}
and noting that 
$\sum_{i=1,\ldots,m_k} \frac{\tilde{p}_i(z)}{v_{\sigma_i}}
v_{\sigma_i}(y^*) = 1$, and by Equation \eqref{ratioassumption} 
in Assumption~1, since $z \leq z'$,
we have,
\begin{align*}
P(Y^{\rm \it BASSOunif}_{k+1}\leq y |Y^{\rm \it BASSOunif}_{k}=z) 
&\geq \frac{\sum_{i=1,\ldots,m_k} \frac{\tilde{p}_i(z')}{v_{\sigma_i}}
v_{\sigma_i}(y)}
{\sum_{i=1,\ldots,m_k} \frac{\tilde{p}_i(z')}{v_{\sigma_i}}
v_{\sigma_i}(y^*)}
\end{align*}
and again noting that $\sum_{i=1,\ldots,m_k} \frac{\tilde{p}_i(z')}{v_{\sigma_i}}
v_{\sigma_i}(y^*) = 1$, we satisfy the condition in \eqref{eq:Condition2}.

For $k = 0$, both ${\rm \it BASSOunif}$ and $HAS\textprime$  sample according to a uniform distribution over $S$, hence the third condition of the Lemma in Appendix~\ref{sec:Lemma30} is satisfied, and we have
%We now can apply Lemma 30 and get 
$P(Y^{\rm \it BASSOunif}_k\leq y) \geq  P(Y^{HAS'}_k\leq y)$ .
\hfill \qed
\\

Let the random variable $Y_k^{\rm \it BASSOsurr}$ be  the incumbent function value $\tilde{y}_k^* $ from BASSO on the $k$th iteration where we generate a point within a selected subregion with surrogate modeling. We make some assumptions regarding efficacy of the surrogate.\\

\noindent{\bf Assumption 2: Sampling within a subregion}
\begin{enumerate}
\item Assume that generating a point $X^{\rm \it BASSOsurr}$ on subregion $\sigma_i$ with surrogate modeling stochastically dominates uniform sampling, i.e., 
\begin{align}\label{eqn:surgassumption}
&P(Y^{\rm \it BASSOsurr}_{k+1}\leq y |Y^{\rm \it BASSOsurr}_{k}=z, X^{\rm \it BASSOsurr} \in \sigma_i) \nonumber \\
&\geq  P(Y^{\rm \it BASSOunif}_{k+1}\leq y |Y^{\rm \it BASSOunif}_{k}=z, X^{\rm \it BASSOunif} \in \sigma_i)
%\\
%& =  P(Y^{BASSOunif}_{k+1}\leq y |Y^{BASSOunif}_{k}=z, X^{BASSOunif}_{k+1} \in \sigma_i) + c_i(y,z)\\
%&=  \frac{v_{\sigma_i}(y)}{v_{\sigma_i}}  + c_i(y,z),
\end{align}

\item Assume that a lower incumbent function value improves the performance of the surrogate, that is,
\begin{align}\label{eqn:IncumbentBASSOsurg}
&P(Y^{\rm \it BASSOsurr}_{k+1}\leq y |Y^{\rm \it BASSOsurr}_{k}=z, X^{\rm \it BASSOsurr} \in \sigma_i) \nonumber \\
&\geq P(Y^{\rm \it BASSOsurr}_{k+1}\leq y |Y^{\rm \it BASSOsurr}_{k}=z', X^{\rm \it BASSOsurr} \in \sigma_i)
\end{align}
for $y_*<y \leq z \leq z' \leq y^*$.
%\\
%& =  P(Y^{BASSOunif}_{k+1}\leq y |Y^{BASSOunif}_{k}=z, X^{BASSOunif}_{k+1} \in \sigma_i) + c_i(y,z)\\
%&=  \frac{v_{\sigma_i}(y)}{v_{\sigma_i}}  + c_i(y,z),
%\end{align*}
\end{enumerate}

\vspace{.2in}
Assumption~2.1 is that surrogate modeling does no worse than uniform sampling on a subregion. 
Notice that the incumbent function value $z$ impacts the adaptive selection probabilities 
$\tilde{p}_i(z)$ but does not impact  uniform sampling on a subregion.  For {\rm \it BASSOunif}, we have
\[P(Y^{\rm \it BASSOunif}_{k+1}\leq y |Y^{\rm \it BASSOunif}_{k}=z, X^{\rm \it BASSOunif} \in \sigma_i) = \frac{v_{\sigma_i}(y)}{v_{\sigma_i}}\]
for any $z$.  It is assumed that surrogate modeling improves upon uniform sampling within a subregion, so~\eqref{eqn:surgassumption} implies that there exists a constant $c_i(y,z) \geq 0$ such that 
\begin{small}
\begin{align}
\label{eq:assume_cyz}
&P(Y^{\rm \it BASSOsurr}_{k+1}\leq y |Y^{\rm \it BASSOsurr}_{k}=z, X^{\rm \it BASSOsurr} \in \sigma_i) 
 =  \frac{v_{\sigma_i}(y)}{v_{\sigma_i}}  + c_i(y,z)
\end{align}
\end{small}

Assumption 2.2 is that a lower incumbent function value improves (or does no worse) surrogate performance than a higher incumbent function value.  Combining \eqref{eqn:surgassumption} with \eqref{eqn:IncumbentBASSOsurg}, we have 
\begin{small}
\begin{align}
\label{eq:assume_dyzprime}
&P(Y^{\rm \it BASSOsurr}_{k+1}\leq y |Y^{\rm \it BASSOsurr}_{k}=z', X^{\rm \it BASSOsurr} \in \sigma_i) 
=  \frac{v_{\sigma_i}(y)}{v_{\sigma_i}}  + d_i(y,z')
\end{align}
\end{small}
%
%\begin{align}
%\label{eq:assume_dyzprime}
%&P(Y^{\rm \it BASSOsurr}_{k+1}\leq y |Y^{\rm \it BASSOsurr}_{k}=z', X^{\rm \it BASSOsurr} \in \sigma_i) \nonumber \\
%&\qquad
%=  \frac{v_{\sigma_i}(y)}{v_{\sigma_i}}  + d_i(y,z')
%\end{align}
where  $c_i(y,z) \geq d_i(y,z')$ for $y_*<y \leq z \leq z' \leq y^*$.

A challenge to implementation is for surrogate modeling to have sufficient coverage in a subregion to improve over uniform sampling on a subregion.  This may appear easy to satisfy, however, if the sampling guided by the surrogate is too exploitative, the amount of improvement to a value $y, y_*<y \leq z$ may be worse than uniform sampling.  This gives a condition on the balance between exploration and exploitation. The second part of the assumption is that, having improved the incumbent value, the surrogate model will increase its probability for improvement to a value $y$. Thus, the accuracy of the surrogate model and the sampling distribution associated with it must satisfy the assumption that includes both the probability of improvement and the {\it amount} of improvement.

We now show that the objective function values of {\rm \it BASSOsurr} stochastically dominate those of 
{\rm \it BASSOunif}.

\begin{theorem}\label{Theorem 2} 
Given the optimization problem \eqref{eq:optproblem} over a continuous/finite domain $S$ with Assumptions~1 and~2 holding, the  function values of BASSO with surrogate sampling model within subregions stochastically dominate those of BASSO with uniform sampling within subregions, i.e., 
\[P(Y^{\rm \it BASSOsurr}_k \leq y) \geq P(Y^{\rm \it BASSOunif}_k \leq y)\]
for $k = 0,1,\ldots$, and  $y_* < y\leq y^*$.
\end{theorem}
\noindent \textbf{Proof:} 
To prove that $Y^{\rm \it BASSOsurr}_k$ stochastically dominates $Y^{\rm \it BASSOunif}_k$ for $k = 0,1,2,\ldots$, we again show that the two sequences 
satisfy the three conditions as stated in Lemma 30 in \cite{Shen2005} (and included in Appendix~\ref{sec:Lemma30}).

To show Condition 1, we need to prove that  
\begin{small}
\[
P(Y^{\rm \it BASSOsurr}_{k+1}\leq y |Y^{\rm \it BASSOsurr}_k=z) \geq  P(Y^{\rm \it BASSOunif}_{k+1}\leq y |Y^{\rm \it BASSOunif}_k=z)
\]
\end{small}
for any fixed $k$, 
$k = 0,1,2,\ldots$,  and any $y$ and $z$ values, for $y_* \leq y, z \leq y^*$. 
We have,
\begin{align*}
&P(Y^{\rm \it BASSOsurr}_{k+1}\leq y |Y^{\rm \it BASSOsurr}_k =z) \\
&\ \  = \sum_{i=1,\ldots,m_k}\tilde{p}_i(z) \ P(Y^{\rm \it BASSOsurr}_{k+1}\leq y | Y^{\rm \it BASSOsurr}_k=z, X^{\rm \it BASSOsurr} \in \sigma_i) 
\\
\insertpar{and by  Equation \eqref{eqn:surgassumption} 
in Assumption~2,  }
&\ \  \geq \sum_{i=1,\ldots,m_k}\tilde{p}_i(z) \  P(Y^{\rm \it BASSOunif}_{k+1}\leq y | Y^{\rm \it BASSOunif}_k\ =z, X^{\rm \it BASSOunif} \in \sigma_i)\\
&\ \  = P(Y^{\rm \it BASSOunif}_{k+1}\leq y |Y^{\rm \it BASSOunif}_k =z).
\end{align*}

To show Condition 2, we need to prove that
$P(Y^{\rm \it BASSOsurr}_{k+1}\leq y |Y^{\rm \it BASSOsurr}_k =z)$ is non-increasing in $z$, that is, we need to prove 
\begin{small}
\begin{align}
\label{eq:Condition2again}
&P(Y^{\rm \it BASSOsurr}_{k+1}\leq y |Y^{\rm \it BASSOsurr}_k =z) \geq 
P(Y^{\rm \it BASSOsurr}_{k+1}\leq y |Y^{\rm \it BASSOsurr}_k =z') 
\end{align}
\end{small}  
for any fixed
$k$,
$k = 0,1,2,\ldots$, and any $y$ and $z$ values, for $y_* < y \leq z \leq z' \leq y^*$. 
If $y \geq z$, 
$P(Y^{\rm \it BASSOsurr}_{k+1}\leq y |Y^{\rm \it BASSOsurr}_k\ =z)
=1$,
hence the inequality in \eqref{eq:Condition2again} is satisfied.
For $y < z$,
 we have,
\begin{align*}
&P(Y^{\rm \it BASSOsurr}_{k+1}\leq y |Y^{\rm \it BASSOsurr}_{k}=z) \\
&\ \  = \sum_{i=1,\ldots,m_k}\tilde{p}_i(z) \ P(Y^{\rm \it BASSOsurr}_{k+1}\leq y | Y^{\rm \it BASSOsurr}_k=z, X^{\rm \it BASSOsurr} \in \sigma_i) 
\\
\insertpar{and by 
Assumption~2 and \eqref{eq:assume_cyz},  }
%
%&\ \  \geq \sum_{i=1,\ldots,m_k}\tilde{p}_i(z) \  P(Y^{\rm \it BASSOunif}_{k+1}\leq y | Y^{\rm \it BASSOunif}_k\ =z, X^{\rm \it BASSOunif} \in \sigma_i)\\
%&\geq 
%P(Y^{\rm \it BASSOunif}_{k+1}\leq y |Y^{\rm \it BASSOunif}_{k}=z)\\
& =\sum_{i=1,\ldots,m_k} \frac{\tilde{p}_i(z)}{v_{\sigma_i}}v_{\sigma_i}(y) + \sum_{i=1,\ldots,m_k} {\tilde{p}_i(z)} \  c_i(y,z).
\end{align*}
We note that 
$\sum_{i=1,\ldots,m_k} \frac{\tilde{p}_i(z)}{v_{\sigma_i}}
v_{\sigma_i}(y^*) = 1$.
Therefore, we have,
\begin{align*}
&P(Y^{\rm \it BASSOsurr}_{k+1}\leq y |Y^{\rm \it BASSOsurr}_{k}=z) \\
& \ \ = \frac{\sum_{i=1,\ldots,m_k} \frac{\tilde{p}_i(z)}{v_{\sigma_i} }
v_{\sigma_i}(y) + \sum_{i=1,\ldots,m_k} {\tilde{p}_i(z)} \ c_i(y,z)}
{\sum_{i=1,\ldots,m_k} \frac{\tilde{p}_i(z)}{v_{\sigma_i}}
v_{\sigma_i}(y^*)}\\
\insertpar{and by \eqref{ratioassumption} 
in Assumption~1, and by \eqref{eq:assume_dyzprime} in Assumption~2, since $z \leq z'$, we have}
& \ \ \geq \frac{\sum_{i=1,\ldots,m_k} \frac{\tilde{p}_i(z')}{v_{\sigma_i}}
v_{\sigma_i}(y) + \sum_{i=1,\ldots,m_k} {\tilde{p}_i(z')} \ d_i(y,z')}
{\sum_{i=1,\ldots,m_k} \frac{\tilde{p}_i(z')}{v_{\sigma_i}}
v_{\sigma_i}(y^*)}\\
\insertpar{and again noting that $\sum_{i=1,\ldots,m_k} \frac{\tilde{p}_i(z')}{v_{\sigma_i}}
v_{\sigma_i}(y^*) = 1$, we have }
&\ \ =
 P(Y^{\rm \it BASSOsurr}_{k+1}\leq y |Y^{BASSOsurg}_{k}=z').
\end{align*}
For $k = 0$, both {\rm \it BASSOsurr} and {\rm \it BASSOunif}  sample according to a uniform distribution over $S$, hence the third condition of the Lemma in Appendix~\ref{sec:Lemma30} is satisfied, and we have
%We now can apply Lemma 30 and get 
$P(Y^{\rm \it BASSOsurr}_k\leq y) \geq  P(Y^{\rm \it BASSOunif}_k\leq y)$.
\hfill \qed
\\

The following corollaries provide an upper bound on the expected number of BASSO function evaluations to get within $\epsilon$ of the global optimum that is linear in dimension.  The corollaries assume a constant lower bound $B$ on the bettering probability in \eqref{eq:betteringprob}. Corollary~\ref{corollary1} states conditions when the domain $S$ is continuous, and Corollary~\ref{corollary2} is stated for a $d$-dimensional integer lattice.
%\\

\begin{corollary} \label{corollary1} Assume Assumptions~1 and~2 hold, and $B$ is a positive lower bound on $b(z)$ in \eqref{eq:betteringprob} for $y_* + \epsilon 
 < z\leq y^*$. For a global optimization problem $(P)$, when $f(x)$  satisfies the Lipschitz condition with Lipschitz constant at most $L$ and  $S$  is a convex feasible region in $d$ dimensions with
a maximum diameter of at most $D$, 
then
%the expected number of $BAS$ points to get within $\epsilon$ of global optimum is linear in dimension 
%When $S$ in $(P)$ is a convex feasible region in $d$ dimensions with
%a diameter $D$, then 
the expected number of function evaluations for BASSO, using uniformly distributed points or surrogate modeling, to reach a value 
$y_* + \epsilon,\  \epsilon> 0$ , is bounded by,
\begin{equation}\label{eqn:corollary1_cont}
\mathbb{E}[N^{\rm \it BASSO} (y_* + \epsilon)] \leq 1+ d\cdot \left(\frac{1}{B} \right) ln\left( \frac{L\cdot D}{\epsilon}\right).
\end{equation}
%where $B = 
%b(y_*+\epsilon)
%{v_S(y_* +\epsilon)}/{v_S}$. 
%is a lower bound on bettering probability $b(t)$ for $y_*+ \epsilon \leq t \leq y^*$.
 \end{corollary}

\textbf{Proof of Corollary 1} 
%Here we present the proof for the case with a continuous domain in corollary~\ref{corollary1} for the expected number of $BAS$ points to get within $\epsilon$ of global optimum for continuous domain as in equation~\eqref{eqn:corollary1_cont}. Then,equation~\eqref{eqn:corollary1_discrete} will follows.
%
By stochastic dominance in Theorem~\ref{Theorem 1}, the expected number of iterations
to achieve a value within $y_*+\epsilon$ for BASSO is less than or equal to the number for $HAS\textprime$, i.e.,
\begin{equation} \label{eq:bas_has_prime_bound} E[N^{BASSO}(y_* +\epsilon)] \leq  
E[N^{HAS\textprime}(y_* +\epsilon)].
\end{equation}

%\begin{comment}
%Since the bettering probability for HAS\textprime is defined as  $b(z) =\sum_{i=1,\ldots,m_k} \frac{\tilde{p}_i(z)}{v_{\sigma_i}}
%v_{\sigma_i}(z)$, for  $y_* < z \leq y^*$, given  subregions $\sigma_i \in \Sigma^C_k$,  $m_k$, and probabilities $\tilde{p}_i(z)$ for $i=1, \ldots, m_k$ with $y_* < z \leq y^*$, using equation~\eqref{eq:has_bound_origin} and equation \eqref{eq:bas_has_prime_bound} we have,
%\end{comment}
%{\color{blue}
%In the continuous case, we show that a lower bound of the bettering probability $b(y)$ for $y_{*} +\epsilon < y < y^*$ is $\frac{v_s(y_*+\epsilon)}{v_S}$.
\noindent
Using \eqref{ratioassumption} in Assumption~1, with $t_z = y$ and $t_{z'} = y^*$, we have 
\[
\frac{\sum_{i=1,\ldots,m_k} \left(\frac{\tilde{p}_i(y)}{v_{\sigma_i}} \right)
{v_{\sigma_i}(y)}}{\sum_{i=1,\ldots,m_k} \left(\frac{\tilde{p}_i(y)}{v_{\sigma_i}} \right)
{v_{\sigma_i}(y^*)}} 
\geq
\frac{\sum_{i=1,\ldots,m_k} \left(\frac{\tilde{p}_i(y^*)}{v_{\sigma_i}} \right)
{v_{\sigma_i}(y)}}{\sum_{i=1,\ldots,m_k} \left(\frac{\tilde{p}_i(y^*)}{v_{\sigma_i}} \right)
{v_{\sigma_i}(y^*)}}.
\]
Note that the numerator on the left-hand side is $b(y)$ by definition in \eqref{eq:betteringprob}, and the   denominator is equal to 1,  because $v_{\sigma_i}(y^*) = v_{\sigma_i}$.
%Hence, 
%\begin{eqnarray*}
%\sum_{i=1,\ldots,m_k} \left(\frac{\tilde{p}_i(y)}{v_{\sigma_i}} \right)
%{v_{\sigma_i}(y^*)} \sum_{i=1,\ldots,m_k}\left(\frac{\tilde{p}_i(y)}{v_{\sigma_i}} \right)
%{v_{\sigma_i}} = \sum_{i=1,\ldots,m_k}{\tilde{p}_i(y^*)} = 1
%\end{eqnarray*}

 The right-hand side denominator is also equal to 1. For the right-hand side numerator,  $\tilde{p}_i(\tilde{y}_k^*)$ represents a uniform sampling probability, 
 %as in \eqref{eqn:piA}, 
 where
 $\tilde{p}_i(\tilde{y}_k^*) = v_{\sigma_i}/v_S$.  We then have 
 \begin{equation} \nonumber
 b(y) \geq \sum_{i=1,\ldots,m_k} \left(\frac{v_{\sigma_i}} {v_S}\frac{1}{v_{\sigma_i}} \right)
{v_{\sigma_i}(y)} = \frac{v_S(y)}{v_S}
 \end{equation}
 for $y_*+\epsilon \leq y \leq y^*$.
% A lower bound of $b(y)$  is $b(y_* +\epsilon)$, since $\frac{v_S(y)}{v_S} \geq \frac{v_S(y_* +\epsilon)}{v_S}$, for $y_* +\epsilon < y< y^*$.
%}
%Since $\frac{v_S(y)}{v_S} \geq \frac{v_S(y_* +\epsilon)}{v_S}$, for $y_* +\epsilon \leq y \leq y^*$, we have a lower bound on $b(y)$, i.e., $b(y) \geq  \frac{v_S(y_* +\epsilon)}{v_S}=B$. 
%
Since $HAS\textprime$ uses  uniform sampling with a bettering probability bounded below by $B$, from \cite{Bulger1998} we have
\begin{align}
 E[N^{HAS'}(y_* +\epsilon)] \leq   1 + \frac{1}{B} \  ln \left( \frac{v_S}{v_S(y+\epsilon)} \right).
\end{align}

Combined with the bounds on $\frac{v_S}{v_S(y_* + \epsilon)}$
in \cite{Bulger1998}, we have the desired bound in \eqref{eqn:corollary1_cont} that is  linear in dimension. \qed
%\begin{equation*}
%\mathbb{E}[N^{BASSO} (y_* + \epsilon)] \leq 1+ d \cdot \left(  \frac{1}{B} \right) ln\left( \frac{L\cdot D}{\epsilon}\right).
%\end{equation*}

\begin{corollary} \label{corollary2}
Assume Assumptions~1 and~2 hold, 
and $B$ is a positive 
 lower bound on $b(z)$ in \eqref{eq:betteringprob} for $y_* + \epsilon 
 < z\leq y^*$.  For a global optimization problem $(P)$, when $S$ in $(P)$ is a $d$-dimensional lattice $\{1,\ldots,k\}^d$, 
%with distinct objective function values, 
the expected number of function evaluations for BASSO, using uniformly distributed points or surrogate modeling, to reach the value $y_*$ is bounded by,
\begin{equation}\label{eqn:corollary1_discrete}
\mathbb{E}[N^{\rm \it BASSO} (y_*)] < 2+ d\cdot  ln(k).
\end{equation}
\end{corollary}

The proof for the finite case of Corollary 2 follows similarly, using \cite{Zabinsky1995}. 
%\hfill \qed

\section{Numerical Experiments}
\label{sec:expts}

We performed numerical experiments with 12 variations of BASSO (Aa, Ab, Ac, Ad, Ba, Bb, Bc, Bd, Ca, Cb, Cc, and Cd) to explore the impact of adaptive subregion probabilities and surrogate modeling on computational performance. These variations combine four adaptive subregion probability methods (a: the best observed function value, b: sample variance, c: Gaussian process on the range distribution, and d: confidence bound on the incumbent value), with three surrogate models (A: uniform, B: Gaussian process, and C: quadratic regression. 
 
In Section~\ref{sec:testAssumption1}, we performed experiments to evaluate whether the
adaptive subregion probabilities $\tilde{p}_i$ (a, b, c, and d) satisfy Assumption~1.  The numerical experiments illustrate the impact of exploitation versus exploration in selecting subregions.  In Section~\ref{sec:TestAssumption2}, we numerically illustrate Assumption~2 by investigating the performance of surrogate modeling (B and C) compared to uniform sampling (A).

In Section~\ref{sec:TestHighDim}, we numerically evaluate the impact of dimension on test problems. We included experiments with 
%the twelve variations of BASSO and 
three algorithms with available code (CMA-ES \cite{CMAES2001}, SNOBFIT \cite{snobfit2008}, and NOMAD\cite{NOMAD2011}). We chose these three algorithms because they ranked highly with good performance in \cite{BAM}. 

The Covariance Matrix Adaptation Evolution Strategy (CMA-ES) \cite{CMAES2001} is a population-based evolutionary algorithm that samples new points from a multivariate Gaussian distribution and evaluates their objective function values. The points are sorted by function values and the distribution parameters (i.e., the mean vector and the covariance matrix) are updated based on the ranking of function values.

The Stable Noisy Optimization by Branch and Fit algorithm (SNOBFIT) \cite{snobfit2008} combines global and local search by branching and fitting a surrogate model on subregions. 
%It is
% an algorithm for bound constrained (and soft constrained) noisy optimization of an expensive objective function. 
SNOBFIT builds local models of the function similar to trust region approaches.
%, and returns in each step a number of points whose evaluation is likely to improve these models or is expected to give better function values. 
%Since the total number of function evaluations is kept low, 
No guarantees can be given that a global minimum is located.

The Non-linear Optimization with Mesh Adaptive Direct Search algorithm (NOMAD) \cite{NOMAD2011} has been used for simulation-based optimization. 
The Mesh Adaptive Direct Search (MADS) algorithm \cite{MAD2006}
generates a series of grids with varying discretizations of the space of variables.
The adaptive mesh frame acts as a window that constrains the search to a specific region of the space,
and allows
an efficient allocation of computational resources.

Although none of these algorithms scale as well as the theoretical ideal, the numerical experiments show the effect of dimension on computational performance.

%We implemented twelve BASSO variations (Aa, Ab, Ac, Ad, Ba, Bb, Bc, Bd, Ca, Cb, Cc, and Cd) and examined their performance on test functions.  
%We performed experiments to evaluate whether 
%adaptive subregion probabilities $\tilde{p}_i$ (a, b, c, and d) satisfy Assumption~1 in Section~\ref{sec:testAssumption1}, and whether surrogate modeling (B and C) satisfy 
%Assumption~2  in Section~\ref{sec:TestAssumption2}.  We also performed experiments {\red on the twelve BASSO variations and on three available algorithms in the literature (CMA-ES, Snobfit, and NOMAD)} to evaluate the impact of dimension, see Section~\ref{sec:TestHighDim}. 

All twelve BASSO variations use the same branching strategy that is described in Section~\ref{sec:BASSOvars}.
The value for the maximum number of function evaluations without improving the incumbent value,  in Step~1 of BASSO, is set to 50. 
The stopping criterion is a fixed  number of function evaluations.

\subsection{Test Assumption 1: Adaptive Subregion Probability Assumption}
\label{sec:testAssumption1}

The four adaptive subregion probability variations (a, b, c, and d) form a proper probability distribution (satisfy Assumptions~1.1 and 1.2) by construction. 
To evaluate whether 
 the adaptive subregion probabilities $\tilde{p}_i$ satisfy \eqref{ratioassumption} in Assumption~1.3, we perform 250 function evaluations of the twelve variations of BASSO 
 %(combinations of 4 adaptive subregion probabilities (a, b, c, and d) and 3 methods to generate a point on the selected subregion (A, B, and C) 
 on the shifted sinusoidal problem and the Rosenbrock problem with dimensions $d=2 \textrm{\; and \;}5$.
During the execution of the algorithm at the $k$th function evaluation when $\tilde{p}$ is updated, we compute the conditional probability ratio on the left-hand side of \eqref{ratioassumption} with $t_z$ as the incumbent function value, and the right-hand side with $t_{z'}$ as the 20\% quantile of the observed objective function value, $t_z < t_{z'}$.

Figures \ref{fig:table_of_figures_a} - \ref{fig:table_of_figures_d} illustrate the 
 left-hand side (LHS) of the probability ratio in \eqref{ratioassumption}
and the right-hand side (RHS) for adaptive subregion probabilities a - d, respectively. 
%The graphs are constructed 
%with $t_z$ as the incumbent function value and
%with $t_{z'}$ as the 20\% quantile of the observed objective function values, $t_z < t_{z'}$.  
The numerical values associated with the graphs are provided in Tables \ref{table:table_assumption1.3_a} - \ref{table:table_assumption1.3_d}, where the incumbent value is labeled $y$, $y=t_z$, and the 20\% quantile is labeled $z$, $t_{z'}$.  

When the LHS is greater than the RHS, Assumption~1.3 is satisfied. The iterations where the estimated LHS is less than the estimated RHS are highlighted in
 Tables \ref{table:table_assumption1.3_a} -  \ref{table:table_assumption1.3_d},
  indicating a likely violation of    Assumption~1.3.

 We observe in Table~\ref{table:table_assumption1.3_c that the} adaptive subregion probability with a Gaussian process on the range distribution (c) performs the best, with fewer violations to Assumption~1.3 than the other three adaptive subregion probabilities. The experiment with adaptive subregion probability (c) has
only one event with a greater RHS on the shifted sinusoidal function and four events on the Rosenbrock function. The predictor function of the cumulative range distribution from the Gaussian process on the range appears to be effective in guiding the adaptive subregion probabilities towards promising subregions. 

%The values of $y,z,t_z \; \textrm{and} \; t_{z'}$  are presented in the Appendix.
As shown in  Table~\ref{table:table_assumption1.3_a},  the adaptive subregion probability with observed best function value (a) satisfies Assumption~1.3 more frequently on the Rosenbrock test function than on the shifted sinusoidal test function.  This may be due to the greedy nature of the adaptive subregion probability (a).  The shifted sinusoidal test function has more local minima than the Rosenbrock test function, so intuitively, a greedy approach would work better on a test function with few local minima.

On the other hand, adaptive subregion probabilities (b) and (d) satisfy Assumption~1.3 more frequently on the shifted sinusoidal test function than on the Rosenbrock test function (Tables \ref{table:table_assumption1.3_b}  and  \ref{table:table_assumption1.3_d}).  These two variations incorporate uncertainty, resulting in more exploration. Exploration is needed on the shifted sinusoidal test function, where the local minima are widely distributed.
 
From these numerical experiments, we conclude that the variation with Gaussian process on the range distribution (c) comes closest to roughly satisfying Assumption~1. Intuitively, variation (c) does the best job of balancing exploration with exploitation.

\subsection{Test Assumption 2: Surrogate Modeling Assumption}
\label{sec:TestAssumption2}
We perform another experiment to 
numerically test Assumption~2 using Gaussian processes (B) and quadratic regression (C) as  surrogate  models 
 and compared to uniform (A).  In order to  test  Assumption~2, we need to compute the probability of improvement at different objective function values and with different methods. 

To calculate the probability of improvements, we constructed a one-dimensional centered sinusoidal problem with two sets of sample points, where most of the points coincide, however the best function value in each subregion differs.  We call one set the base case sample set and the other the better-incumbent sample set.
See Figure~\ref{fig:sampleset} in the appendix. 
%Figure~\ref{fig:sampleset} in the Appendix illustrates the two sample sets on the test function, with three subregions. 

%We perform another experiment to 
%numerically test Assumption~2, using Gaussian processes (B) and quadratic regression (C) as  surrogate  models.
%In order to compute the probability of improvement used in  Assumption~2, we used a one-dimensional centered sinusoidal problem  with two sets of sample points, where most of the points coincide, however the best function value in each subregion differs.  We call one set the base case sample set and the other the better incumbent sample set.
%See Figure~\ref{fig:sampleset} in the appendix. 
%Figure~\ref{fig:sampleset} in the Appendix illustrates the two sample sets on the test function, with three subregions. 

As shown in Table~\ref{table:assumption2}, the probabilities of improvement using uniform sampling are smaller than those using the base case sample set for the Gaussian process and for quadratic regression. 
 This supports Assumption~2.1 in \eqref{eqn:surgassumption} that surrogate modeling with Gaussian processes and quadratic regression can outperform uniform sampling.
We also observe that the probabilities of improvement using the better-incumbent  sample set (with $z$) are higher than those using the base case sample set (with $z^\prime$). We calculate the probability of improvement  within a subregion to be below a value $y$, where  $y$ is the 5th best function value  of both sample sets. This supports Assumption 2.2 in \eqref{eqn:IncumbentBASSOsurg} that a better incumbent function value improves surrogate performance.

The probability of improvement using a uniform distribution is calculated by estimating $v_{\sigma_i}(y)/v_{\sigma_i}$.
The probability of improvement using Gaussian processes is 
%below a value $y$ within a subregion, as follows
\begin{equation}\label{eqn:PI}
PI_i(x) =\Phi\left(\frac{y - \hat{g}_{i,k}(x)}{\hat{s}_{i,k}(x)}\right)
\end{equation}
for $x\in \sigma_i$, where $\hat{g}_{i,k}(x)$ and $\hat{s}_{i,k}(x)$ are constructed from the Gaussian process for both  
the base case sample set, and  the  better-incumbent  sample set.
The probability of improvement in Table~\ref{table:assumption2} is assessed at the point that maximizes the expected improvement function $EI_i(x)$. 
Figure~\ref{fig:PI} in the appendix illustrates the expected improvement functions $EI_i(x)$, with maximal points, for both sample sets.

%as in \eqref{eq:EIx}, for each subregion, where Figure~\ref{fig:PIbase} is for the base case sample set, and Figure~\ref{fig:PIbetter} is for the better incumbent sample set.  Notice that there is a different point selected between the base case sample set and the better incumbent sample set, due to the different best function values.

%In order to demonstrate that Assumption~2 is satisfied, we calculate the probability of improvement, 
%at the new sample point, 
%comparing the probability of improvement using the Gaussian process and the probability of improvement using regularized quadratic regression to that of a uniformly distributed point (Assumption 2.1, \eqref{eqn:surgassumption}), and comparing the probability of improvement between the base case sample set and the better incumbent sample set (Assumption 2.2, \eqref{eqn:IncumbentBASSOsurg}). We calculate the probability of improvement  within a subregion to be below a value $y$, where  $y$ is the 5th best function value  of both sample sets.
%The probability of improvement using a uniform distribution is calculated by estimating $v_{\sigma_i}(y)/v_{\sigma_i}$, and using Gaussian processes as 
%below a value $y$ within a subregion, as follows
%\begin{equation}\label{eqn:PI}
%PI_i(x) =\Phi\left(\frac{y - \hat{g}_{i,k}(x)}{\hat{s}_{i,k}(x)}\right).
%\end{equation}

The probability of improvement using regularized quadratic regression is calculated using
\begin{equation}\label{eqn:PIlasso}
PI_i(x) = P\left(T(x) < \frac{y-\hat{q}_{i,k}(x)}{\sqrt{MSE_i + (SE_i(x))^2}}\right),
\end{equation}
for $x \in \sigma_i$, where  $y$ is the 5th best function value  of both sample sets.
The random variable $T(x)$ represents the true objective function value at $x$, $f(x)$, and follows a 
 $t$ distribution 
 %$T(x) \sim t$
 with degrees of freedom set to  $N^i_k-3$.
% number of model parameter, where
%$T(x) \sim t-distribution\textrm{(degree of freedom} = N^i_k - 2)$.
%Table \ref{table:assumption2} provides the estimates of the probability of improvement using uniform sampling,  the base case sample set, and  the better incumbent sample set. 
%{\red{Pete - please check this.}}
%
The
mean squared error for subregion $\sigma_i$ is
\begin{equation*}\label{eqn:MSE}
MSE_i = \frac{1}{N^i_k}  \sum_{j=1}^{N^i_k} (\hat{q}^i_k (x_{j})- f(x_j))^2,
\end{equation*}
and the standard error at $x\in \sigma_i$ is
\begin{equation*}
% SE_i(x) =   (\hat{q}^i_k (x)- f(x))^2 
SE_i(x) =  \sqrt{MSE_i(Q^{\textrm{T}}(\mathbf{X}^{\textrm{T}}\mathbf{X})^{-1}Q)},
\end{equation*}
where $Q$ is a vector of quadratic regression terms at point $x$.  The size of $Q$ is $3$ for $d =1$ and the size is  $1 + 2d + \Mycomb[d]{2}$ for $d > 1$.
%
%\noindent
For example, 
\begin{itemize}
    \item when $x$ is 1-dimensional, $Q = \{1,x_1,x_1^2\}$,
    \item when $x$ is 2-dimensional,
$Q = \{1,x_1,x_2,x_1^2,x_2^2,x_1x_2\}$,
    \item and when $x$ is 3-dimensional,
$Q = \{1,x_1,x_2,x_3,x_1^2,x_2^2,x_3^2,x_1x_2,x_2x_3,x_1x_3\}$,
\end{itemize}
\noindent
where $x_{\ell}$ is the $\ell$th component (decision variable) of point $x$.

%when $x_{\ell}$ is the $\ell$ component (decision variables) of point $x$\\
%for 1-dimensional vector $x, X = \{1,x_1,x_1^2\}$,\\
%for 2-dimensional vector
%$x, X = \{1,x_1,x_2,x_1^2,x_2^2,x_1x_2\}$,\\
%for 3-dimensional vector
%$x, X = \{1,x_1,x_2,x_3,x_1^2,x_2^2,x_3^2,x_1x_2,x_2x_3,x_1x_3\}$.

The design matrix $\mathbf{X}$ has a row for every point sampled in subregion $\sigma_i$, so the size of $\mathbf{X}$ is $N^i_k \times (1 + 2d + \Mycomb[d]{2})$.

%{\blue{Consider subregion $\sigma_i$ with 4 sample points 
%($N^i_k = 4$) for 1-dimensional problem. Those 4 sample points are $\{1,2,4,5\}$.
%We want to find $SE_i(x)$ where $x = 2$}.\\
%Then, vector $X$ of $x$ is $X = \{1,x_1,x_1^2\} = \{1,2,4\}$\\

%The size of Matrix $\mathbf{X}$ of subregion $\sigma_i$ will be  $N^i_k \times 3 = 4 \times 3 $

%$$ \mathbf{X} = \begin{bmatrix}
%1&1&1^2\\
%1&2&2^2\\
%1&4&4^2\\
%1&5&5^2
%\end{bmatrix} = %\begin{bmatrix}
%1&1&1\\
%1&2&4\\
%1&4&16\\
%1&5&25\\
%\end{bmatrix}$$
%When we do matrix multiplication the dimension will be 
%$X^{\textrm{T}}_{1 \times 3}(\mathbf{X}^{\textrm{T}}_{3 \times 4}\mathbf{X}_{4 \times 3})^{-1}X_{3 \times 1}$, the result will be scalar.
%}

%For example, for a problem with $d = 2$, in a subregion with $N^i_k$ sample points, let $x_{\ell,j}$ be the $\ell$th component  of the $j$th sample point $x$, $j=1, \ldots, N^i_k$. 
%
%The design matrix $\mathbf{X}$ for that subregion is,
%$$ \mathbf{X} = \begin{bmatrix}
%1&x_{1,1}&x_{2,1}&x_{1,1}^2,&x_{2,1}^2&x_{1,1}x_{2,1}\\
%1&x_{1,2}&x_{2,2}&x_{1,2}^2&x_{2,2}^2&x_{1,2}x_{2,2}\\
%.&.&.&.&.&.\\
%.&.&.&.&.&.\\
%.&.&.&.&.&.\\
%1&x_{1,{N^i_k}}&x_{2,{N^i_k}}&x_{1,{N^i_k}}^2&x_{2,{N^i_k}}^2&x_{1,{N^i_k}}x_{2,{N^i_k}}
%\end{bmatrix}.$$

The probability of improvement in Table~\ref{table:assumption2} is assessed at the point that minimizes the quadratic regression function.
Figure~\ref{fig:PI} in the appendix illustrates the quadratic regression  $\hat{q}(x)$, with minimal points, for both sample sets.

\begin{table}[ht!]
\centering
\small
\begin{tabular}{crrr}
  \hline
Estimate  & $\sigma_1$ & $\sigma_2$ & $\sigma_3$   \\ 
  \hline

{\scriptsize {$P(Y^{\rm \it BASSOunif}_{k+1}\leq y |Y^{\rm \it BASSOunif}_{k}=z, X^{\rm \it BASSOunif} \in \sigma_i)$}} & 0.252 & 0.833 & 0.252    \\ 
%          \qquad $X^{\rm \it BASSOunif} \in \sigma_i)$     \\ 
{\scriptsize {$P(Y^{\rm \it BASSOsurrGP}_{k+1}\leq y |Y^{\rm \it BASSOsurrGP}_{k}=z', X^{\rm \it BASSOsurrGP} \in \sigma_i)$}} & 0.276 & 1.000 &0.284    \\ 
%          \qquad $X^{\rm \it BASSOsurr} \in \sigma_i)$ \\
{\scriptsize {$P(Y^{\rm \it BASSOsurrGP}_{k+1}\leq y |Y^{\rm \it BASSOsurrGP}_{k}=z, X^{\rm \it BASSOsurrGP} \in \sigma_i)$}} & 0.411 & 1.000 &0.558    \\ 
 %         \qquad $X^{\rm \it BASSOsurr} \in \sigma_i)$ \\

{\scriptsize {$P(Y^{\rm \it BASSOsurrQ}_{k+1}\leq y |Y^{\rm \it BASSOsurrQ}_{k}=z', X^{\rm \it BASSOsurrQ} \in \sigma_i)$}} & 0.267 & 0.982 &0.998   \\ 
 %         \qquad $X^{\rm \it BASSOsurr} \in \sigma_i)$ \\

{\scriptsize {$P(Y^{\rm \it BASSOsurrQ}_{k+1}\leq y |Y^{\rm \it BASSOsurrQ}_{k}=z, X^{\rm \it BASSOsurrQ} \in \sigma_i)$}} &0.674 &0.986 &0.998   \\ 
 %         \qquad $X^{\rm \it BASSOsurr} \in \sigma_i)$ \\
      \hline
\end{tabular}
\caption{Result from numerical experiment to test Assumption~2.  Observe that the probability of improvement of BASSO with the surrogate modeling (Gaussian processes and quadratic regression) can outperform uniform sampling. The probabilities of improvement using the better-incumbent  sample set  (with $z$) are higher than those using the base case sample set (with $z^\prime$).}
\label{table:assumption2}
\end{table}

\subsection{Test Problems in High Dimensions}
\label{sec:TestHighDim}
Although the conditions leading to ideal scalability, as in the linearity result, are not completely satisfied, we explore the performance of BASSO implementations and  three algorithms  (CMA-ES \cite{CMAES2001}, SNOBFIT \cite{snobfit2008}, and NOMAD\cite{NOMAD2011})  on test problems to evaluate the impact of dimension.

Each algorithm is applied to six benchmark functions  (Ackley, Hartmann, Branin, centered sinusoidal, shifted sinusoidal, and Rosenbrock, details in Appendix~\ref{sec:testfunctions}).
Due to computational limitations, the Gaussian process surrogate (B) is restricted to problems with dimensions of 20 and 50. The rest of the algorithms are run on all problems across dimensions of 20, 50, 100, 500, and 1000.  Ten independent replications are performed using common random numbers and the average, minimum and maximum function values are recorded.

To compare the performance of the different algorithms, we use {\it performance profiles}, as in \cite{more2009}. 
We are interested in the fraction of problems $p \in P$ that can be solved within a given tolerance $\tau$ by solver $s \in S $ in $K$ function evaluations. 

For a given tolerance $0 \leq \tau \leq 1$ and common starting point $x_0$, a problem is
considered `solved’ by a solver  if its solution  improved the starting point by
at least a fraction of $(1 - \tau )$ of the largest possible reduction by any solver in the set of interest within a given number of function evaluations,  i.e., the convergence test in \eqref{eqn:convtest} is satisfied in the given number of function evaluations. The convergence test is 
\begin{equation}\label{eqn:convtest}
  f(x_0)- f_{s} \geq (1-\tau)(f(x_0)- f_L)
\end{equation}
where $f(x_0)$ is the average function value of the starting point of the 10 replications,  $f_s$ is the average function value over 10 replications
reported by solver $s$ on the $K$th function evaluation, and $f_L$ is the smallest objective function value obtained by any solver within $K$ function evaluations.
We use the binary metric $t_{p,s}^{\tau,K}$ to indicate if the convergence test in \eqref{eqn:convtest} is satisfied by solver $s \in S $  on  problem $p \in P$ with $K$ function evalutions at tolerance $\tau$. That is,
\begin{equation*}
   t_{p,s}^{\tau,K} = \left\{
   \begin{array}{l l}
   1 \quad \text{if the convergence test \eqref{eqn:convtest} is satisfied by solver $s$ on problem $p$} \\ \quad \; \; \text{within $K$ function evaluations at tolerance $\tau$} \\
   \ \  & \\
   0 \quad \mathrm{otherwise}
   \end{array}
   \right.  
\end{equation*}
The fraction of problems solved by solver $s$ is
\[
d_s = \frac{1}{|\mathcal{P}|}\sum_{p\in \mathcal{P}} t_{p,s}^{\tau,K}.
\]

\begin{figure}[htbp!]
     \centering
     \begin{subfigure}[b]{0.45\textwidth}
         \centering
         \includegraphics[width=\textwidth]{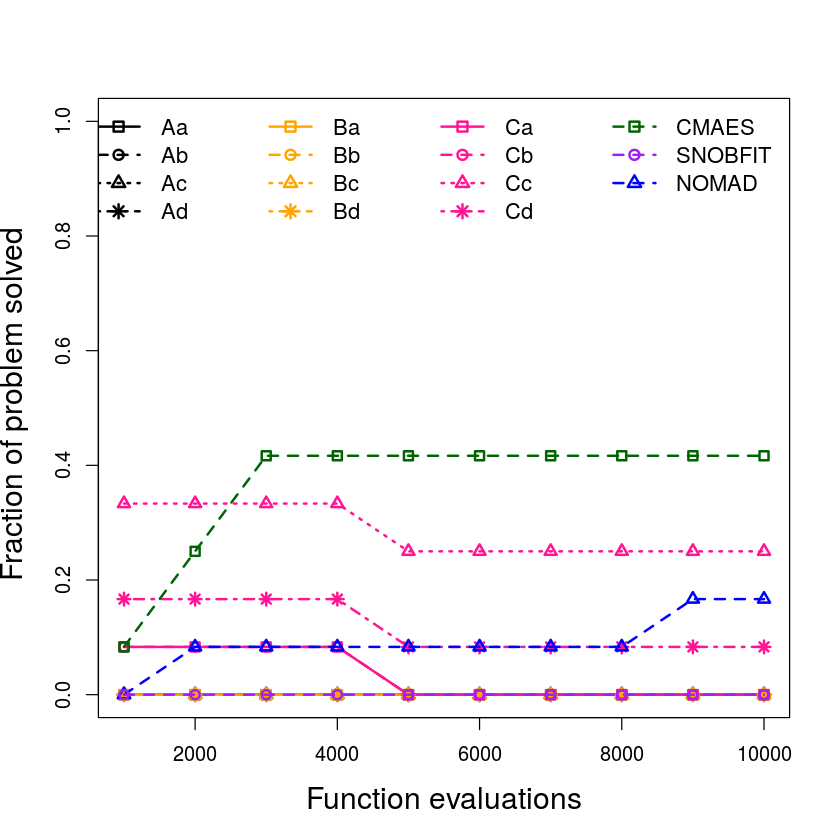}
         \caption{$\tau = 0.001$}
         \label{fig:meantau-3_12}
     \end{subfigure}
     \begin{subfigure}[b]{0.45\textwidth}
         \centering
         \includegraphics[width=\textwidth]{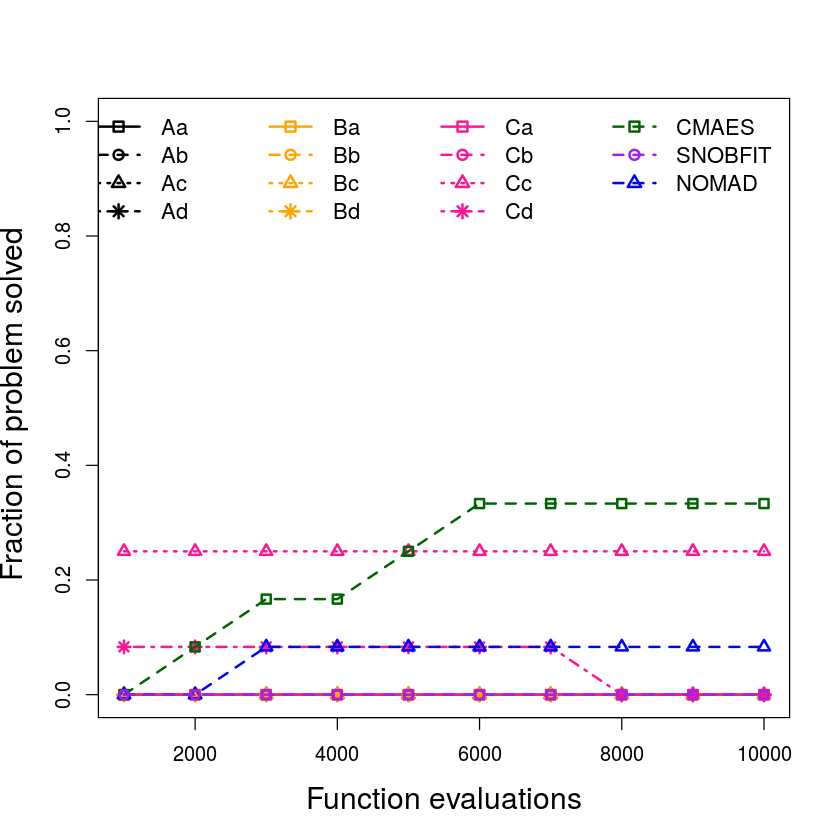}
         \caption{$\tau = 0.00001$}
         \label{fig:meantau-5_12}
     \end{subfigure}
        \caption{
        Performance profiles for 15 algorithms on 12 benchmark problems in dimensions 20 and 50, with $\tau = 0.001$ in (a) and $\tau = 0.00001$ in (b).
%        In (a), performance profile $d_s$ for $\tau = 0.001$.In (b), performance profile $d_s$ for $\tau = 0.0001$. We define $f(x_0)$ as the average function value at the beginning over 10 replications, $f_{s}$ as the average function value at chosen function evaluations over 10 replications and $f_{L}$ as the smallest objective function over 10 replications of all $s \in S$ at chosen function evaluation. We perform convergence test on on 12 benchmark problems on dimension 20 and 50.
        }

        \label{fig:20_50_mean}
\end{figure}

In Figure~\ref{fig:20_50_mean}, we compare the performance profiles of the 12 variations of BASSO (Aa, Ab, Ac, Ad, Ba, Bb, Bc, Bd, Ca, Cb, Cc, Cd) and  CMA-ES, SNOBFIT, and NOMAD, on the six test functions in dimensions 20 and 50. 
The performance profiles are presented for
$\tau = 1e-3$ and $\tau = 1e-5$.

%For completeness, the average, minimum and maximum function values are plotted across the number of function evaluations are provided in Figure~APPENDIX20-50 \ref{fig:table_of_figures_exper_2_1}, for the 15 algorithms, six test functions, and dimensions 20 and 50.

We observe in Figure~\ref{fig:20_50_mean}, that the BASSO variations with the regularized quadratic regression as the surrogate model (C) outperform the variations using no surrogate model (A) and variations using Gaussian processes (B). The performance profiles for (A) and (B) are barely visible at the bottom of the graph, as they are on top of each other. 
It was expected that both variations with surrogate models (B and C) would outperform the variation with no surrogate (A) because they both appeared to satisfy Assumption~2, 
but it was a surprise that the regularized quadratic regression performed better than Gaussian processes.  One intuitive explanation is that the quadratic regression captured the overall trend of the objective function whereas a Gaussian process is better at approximating the detailed behavior.

We also observe in Figure~\ref{fig:20_50_mean} that the adaptive subregion probability with a Gaussian process on the range distribution (c) performs the best, which is consistent with  our numerical experiment showing it had fewer violations to Assumption~1.3 than the other three adaptive subregion probabilities. 

\begin{figure}[htbp!]
     \centering
     \begin{subfigure}[b]{0.45\textwidth}
         \centering
         \includegraphics[width=\textwidth]{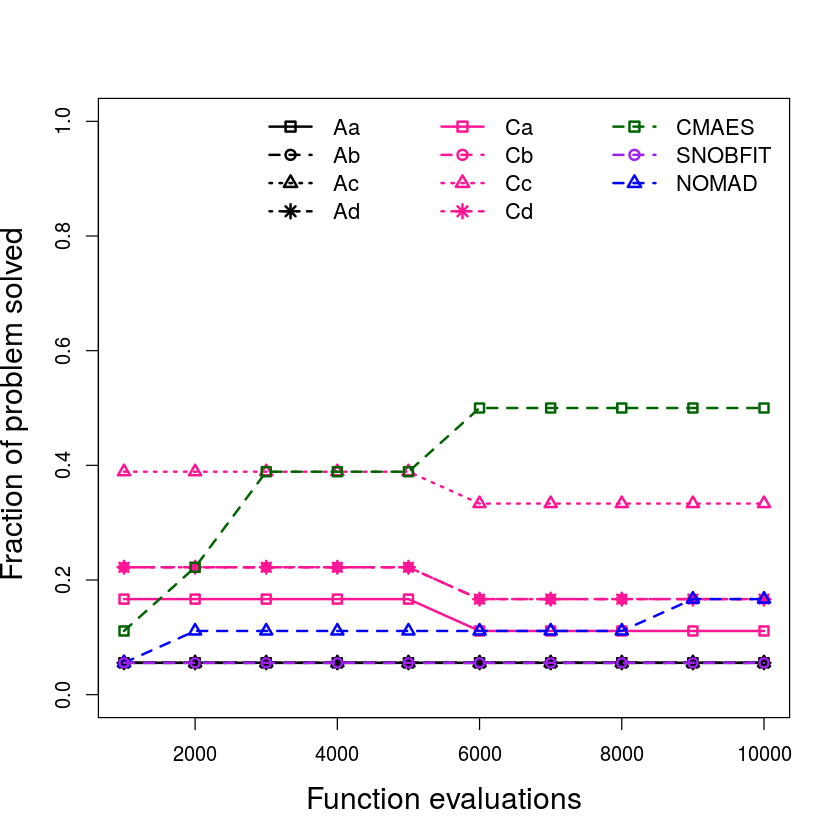}
         \caption{$\tau = 0.001$}
         \label{fig:meantau-3_18}
     \end{subfigure}
     \begin{subfigure}[b]{0.45\textwidth}
         \centering
         \includegraphics[width=\textwidth]{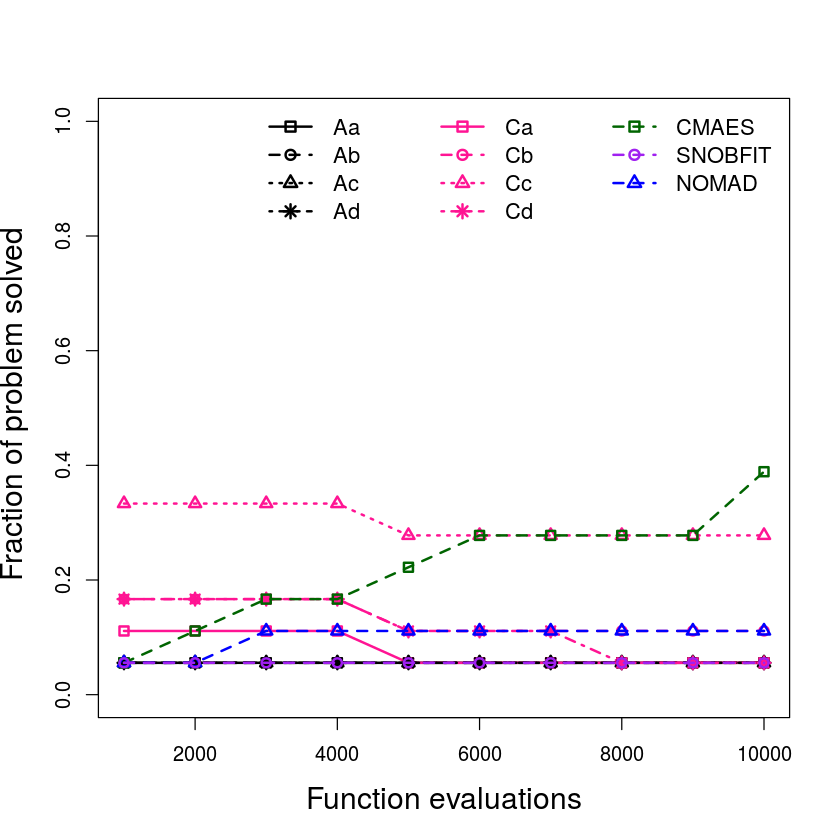}
         \caption{$\tau = 0.00001$}
         \label{fig:meantau-5_18}
     \end{subfigure}
        \caption{
        Performance profiles for 11 algorithms on 18 benchmark problems in dimensions 20, 50 and 100, with $\tau = 0.001$ in (a) and $\tau = 0.00001$ in (b).
        %In (a), performance profile $d_s$ for $\tau = 0.001$.In (b), performance profile $d_s$ for $\tau = 0.0001$. We define $f(x_0)$ as the average function value at the beginning over 10 replications, $f_{s}$ as the average function value at chosen function evaluations over 10 replications and $f_{L}$ as the smallest objective function over 10 replications of all $s \in S$ at chosen function evaluation. We perform convergence test on on 18 benchmark problems on dimension 20 ,50 and 100
        }
        \label{fig:mean_18}
\end{figure}

In Figure~\ref{fig:mean_18}, we expand the test problems to include dimensions 20, 50, and 100, but we drop the BASSO variations (B) with Gaussian processes. We see again that adaptive subregion probability (c) outperforms that other adaptive subregion probabilities.  Of the BASSO variations, (Cc) performs the best.  We observe that CMA-ES overtakes BASSO (Cc) after several thousand function evaluations.  
We speculate that the population-based aspect of CMA-ES takes awhile for the population to represent the promising regions, but once it does, then it does a better job of identifying the improving regions than the adaptive subregion probability (c).  This suggests that the adaptive subregion probability (c) is useful in the beginning of an execution, and perhaps can provide information to CMA-ES for a population-based approach.
%While we cannot apply our theory to CMA-ES, we speculate that the population-based aspect is performing better than the adaptive subregion probability (c) at identifying the improving regions.

For completeness, the average, minimum and maximum function values are plotted across the number of function evaluations  in Figure~\ref{fig:table_of_figures_exper_2_1} for the 15 algorithms, six test functions, and dimensions 20 and 50. Figure~\ref{fig:table_of_figures_exper_2_2} provides similar plots for 11 algorithms, six test functions, and dimensions 100, 500, and 1000.  The performance profiles allow a comparison for test functions up to 100 dimensions.  We decided to push the numerical experiments to 500 and 1,000 dimensions.  
As shown in Figure~\ref{fig:table_of_figures_exper_2_2}, the BASSO variation with regularized quadratic regression (C) and CMA-ES are effective on four of the six test functions in 500 dimensions.  The centered and shifted sinusoidal functions are challenging for all of the algorithms even at 100 dimensions.  The number of local minima grows exponentially with dimension for these two test functions, which is very challenging.
The insights from the BASSO variation with regularized quadratic regression (C) and CMA-ES are that they seem to balance exploration with exploitation very well.  They both include an element of randomness to allow for exploration, while updating the method of sampling to focus on improving regions.

\section{Conclusion}
\label{sec:conclusion}
We propose BASSO, a framework for adaptive random search that conceptualizes partitioning and surrogate modeling to solve black-box optimization problems. We derive an upper  bound on the expected number of function evaluations until a specified $\epsilon$ above the global minimum value is reached. We show that under certain conditions, the expected number of function evaluations is bounded by a  linear function of the domain dimension, achieving scalability to problems in high dimensions.
This conceptual framework balances exploration  through an adaptive subregion probability to focus the search on subregions with improving function values, and with exploitation through a surrogate model that optimizes within a subregion.  

We present several implementations of BASSO and numerically conclude that the best variation uses a regularized quadratic regression as a surrogate, and a one-dimensional Gaussian process on the range distribution for the adaptive subregion probability. 
The combination of exploration through the adaptive subregion probability and exploitation via the quadratic regression has potential for scalability.
Further research is needed to reduce the gap between implementation and the theoretical ideal performance.
The finite-time analysis of BASSO may inspire future implementations that can scale to high dimensional black-box optimization problems.\\

%\bmhead{Supplementary information}
% {\bf Supplementary information}

%\bmhead{Acknowledgments}
%\noindent
%{\bf Acknowledgments}
\section*{Acknowledgments}
This work has been funded in part by the U.S. National Science
Foundation grants CMMI-2204872, and CMMI-2046588.

\section*{Data Availability Declaration}
The computer code and accompanying test data will be made available upon request.
%\begin{itemize}
%\item Funding
%\item Conflict of interest/Competing interests (check journal-specific guidelines for which heading to use)
%\item Ethics approval and consent to participate
%%%\item Materials availability
%\item Code availability 
%\item Author contribution
%\end{itemize}

% \noindent
% If any of the sections are not relevant to your manuscript, please include the heading and write `Not applicable' for that section. 

\newpage
\begin{appendices}

\section{Overview of Hesitant Adaptive Search}
\label{sec:HAS}

HAS is defined by a sampling distribution $\delta$ with support on $S$ and a bettering probability $b(y)$, $0  < b(y)\leq 1$, defined for $y_* < y \leq y^*$.  
%The bettering probability $b(y)$
%depends on $y$, the best objective function value observed (incumbent).
On any iteration with incumbent  value $y$, HAS generates an improving point with probability $b(y)$ by drawing from the normalized restriction of $\delta$ on the improving level set.

\noindent\textbf{Hesitant Adaptive Search (HAS), cf.~\cite{Bulger1998} }

\textbf{Step 0}:  
 Initialize $X_0$ in $S$ according to the 
 sampling probability measure $\delta$ on $S$. 
 Set $k = 0$. Set $Y_0$ = $f(X_0)$.

\textbf{Step 1}:  Generate $X_{k+1}$ on the improving set
$S_{k+1} = \{ x \in S : f(x) < Y_k \}$ with probability $b(Y_k)$, and otherwise set $X_{k+1} = X_k.$ Set $Y_{k+1} = f(X_{k+1})$.

\textbf{Step 2}: If a stopping criterion is met, stop. Otherwise, increment $k$ and return to Step 1.

\section{Lemma for Stochastic Dominance}
\label{sec:Lemma30}
\textbf{Lemma  (cf. Lemma 30 in \cite{Shen2005})} 
 Let $\{Y^A_k,k = 0,1,2,\ldots\}$ and $\{Y^B_k,k = 0,1,2,\ldots\}$
be two sequences of objective
function values generated by algorithms A and B respectively for solving problem 
\eqref{eq:optproblem} 
where $Y^A_{k+1} \leq Y^A_{k}$ and $Y^B_{k+1} \leq Y^B_{k}$, for $k = 0,1,2,\ldots$. For $y_*\leq y \leq y^*$ and $k = 0,1,2,\ldots$, if

\begin{itemize}
\item[1.] $P(Y^A_{k+1}\leq y |Y^A_k\ =z) \geq  P(Y^B_{k+1}\leq y |Y^B_k\ =z),$
\item[2.] $P(Y^A_{k+1}\leq y |Y^A_k\ =z)$ is non-increasing in z, and
\item[3.] $P(Y^A_0\leq y) \geq  P(Y^B_0\leq y),$
\end{itemize}
then $P(Y^A_k\leq y) \geq  P(Y^B_k\leq y)$ for $k= 0,1,2,\ldots$ and  $y_*\leq y \leq y^*$.

\section{Testing Assumption 1.3 Numerically}
\label{sec:Test_1.3}

\begin{figure}
    \centering
    \setkeys{Gin}{width=\linewidth}
    \settowidth\rotheadsize{PARAMETERS 3} 
    % from makecell
\begin{tblr}{colspec = { Q[h] *{4}{Q[c,m, wd=30mm]}},
             colsep  = 1pt,
             cell{2-Z}{1} = {cmd=\rotcell, font=\footnotesize\bfseries},
             row{1} = {font=\bfseries},
             measure=vbox
            }
    \toprule
 & Shifted Sinusoidal (dimension=2)& Shifted Sinusoidal (dimension=5)& Rosenbrock (dimension=2)&  Rosenbrock (dimension=5)\\
    \midrule
\centering
A Uniform & 
{\includegraphics{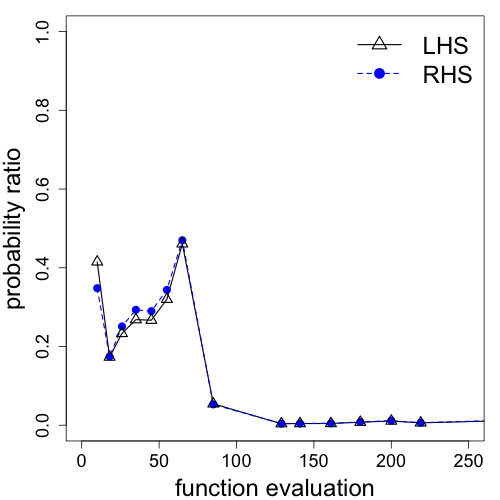}} 
        & {\includegraphics{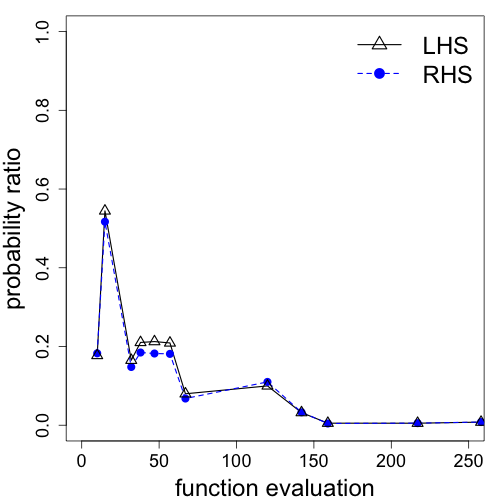}}
            & {\includegraphics{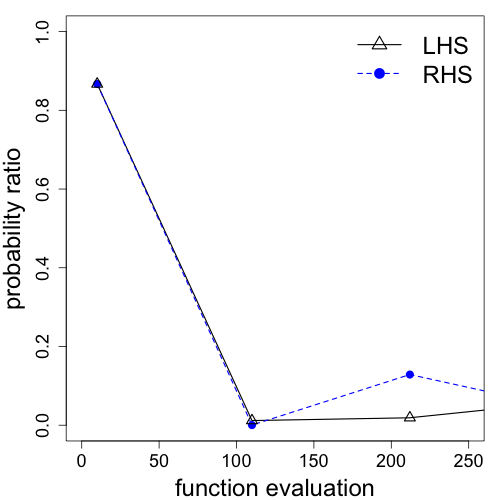}}
                & {\includegraphics{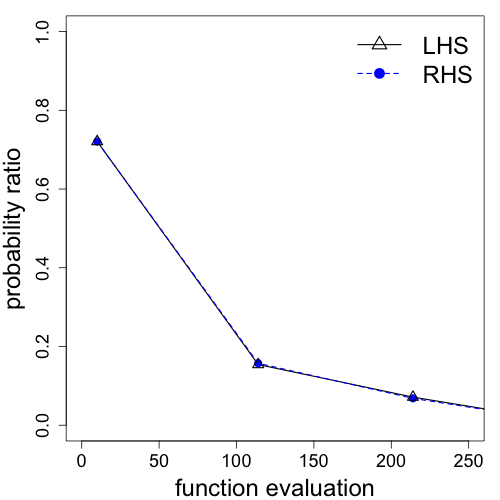}}     \\
\centering
B Gaussian Process& {\includegraphics{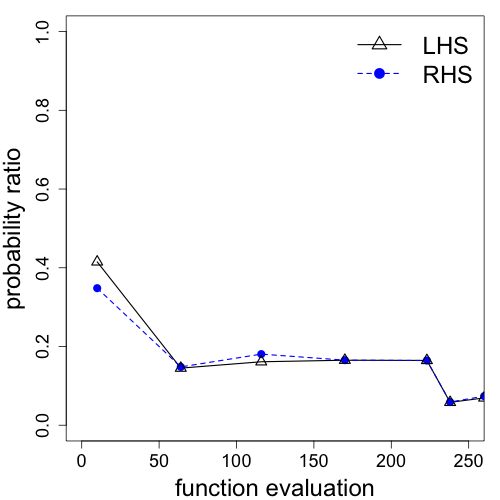}}
        & {\includegraphics{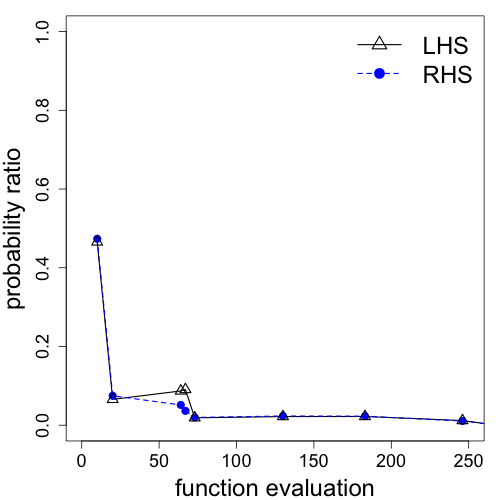}}
            & {\includegraphics{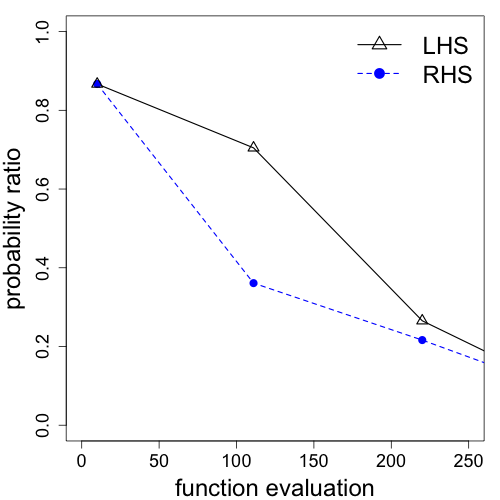}}
                & {\includegraphics{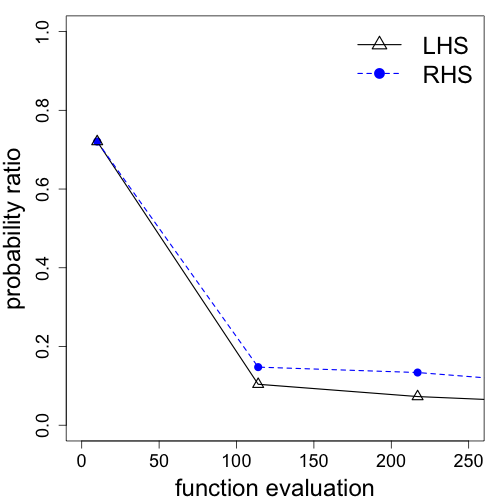}}    \\
\centering
C Quadratic Regression & {\includegraphics{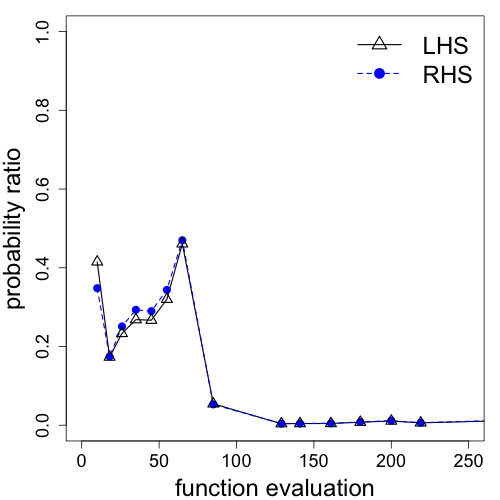}}
        & {\includegraphics{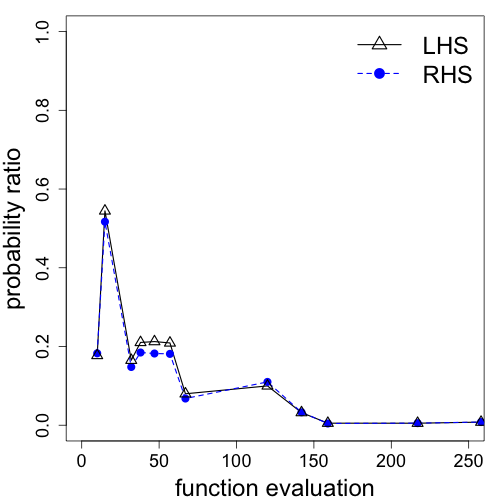}}
            & {\includegraphics{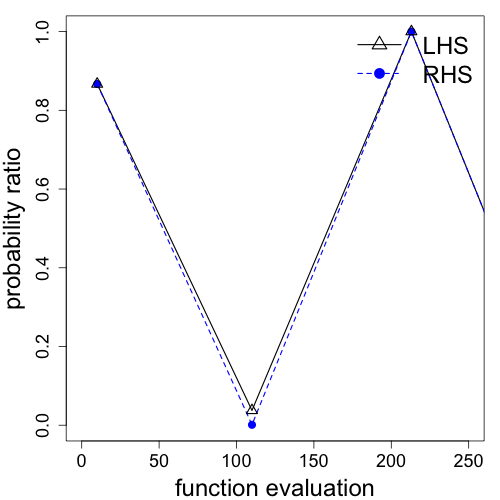}}
                & {\includegraphics{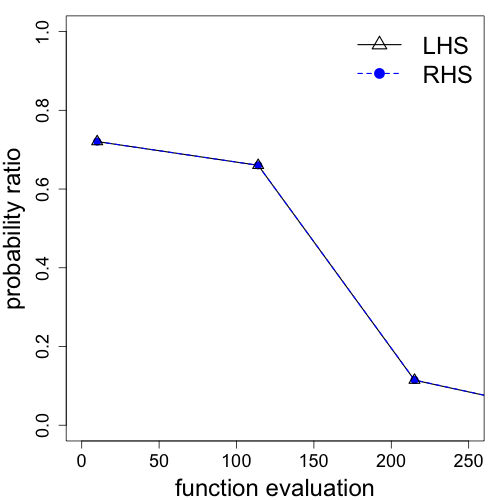}}\\
    \bottomrule
    \end{tblr}
\caption{BASSO variations (Aa, Ba, Ca) with adaptive subregion probabilities in \eqref{eq:ptildeincumbent}, observed best function value (a). When the left-hand side (LHS) of the probability ratio in \eqref{ratioassumption}
with $t_z$ as the incumbent function value is  larger than the right-hand side (RHS)
with $t_{z'}$ as the 20\% quantile of the observed objective function values, $t_z < t_{z'}$,   Assumption 1.3 is satisfied.
}
\label{fig:table_of_figures_a}
\end{figure}

\begin{figure}
    \centering
    \setkeys{Gin}{width=\linewidth}
    \settowidth\rotheadsize{PARAMETERS 3}    % from makecell
\begin{tblr}{colspec = { Q[h] *{4}{Q[c,m, wd=30mm]}},
             colsep  = 1pt,
             cell{2-Z}{1} = {cmd=\rotcell, font=\footnotesize\bfseries},
             row{1} = {font=\bfseries},
             measure=vbox
            }
    \toprule
 & Shifted Sinusoidal (dimension=2)& Shifted Sinusoidal (dimension=5)& Rosenbrock (dimension=2)&  Rosenbrock (dimension=5)\\
    \midrule
\centering
A Uniform & {\includegraphics{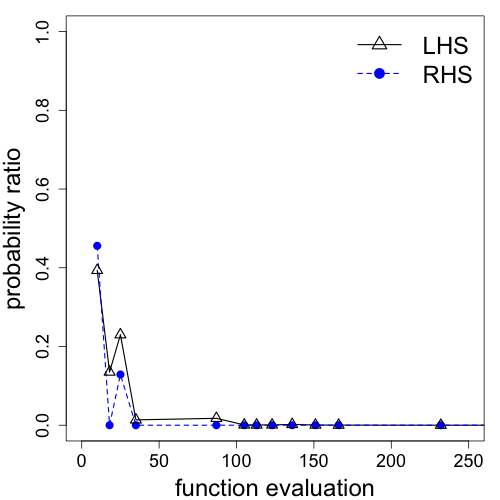}}
        & {\includegraphics{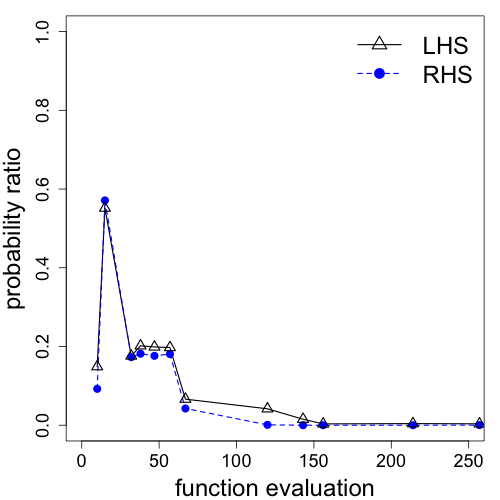}}
            & {\includegraphics{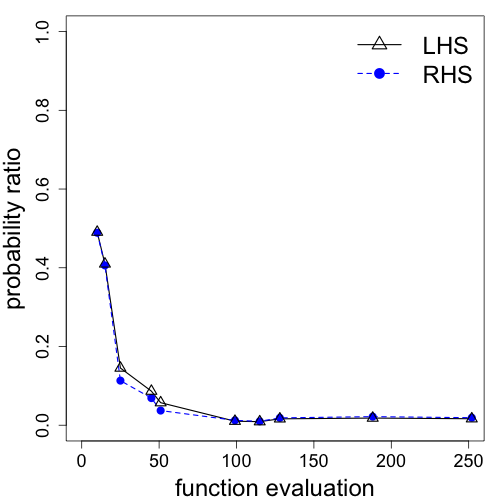}}
                & {\includegraphics{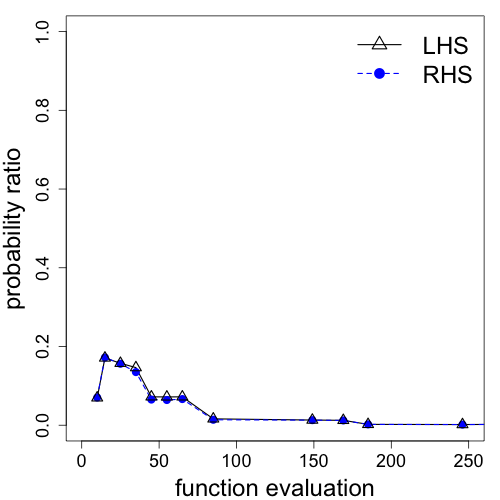}}     \\
\centering
B Gaussian Process& {\includegraphics{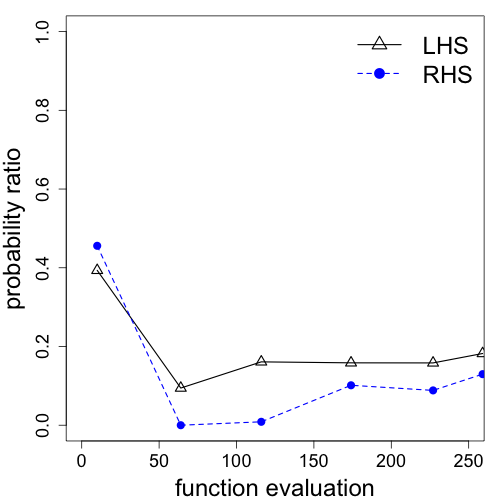}}
        & {\includegraphics{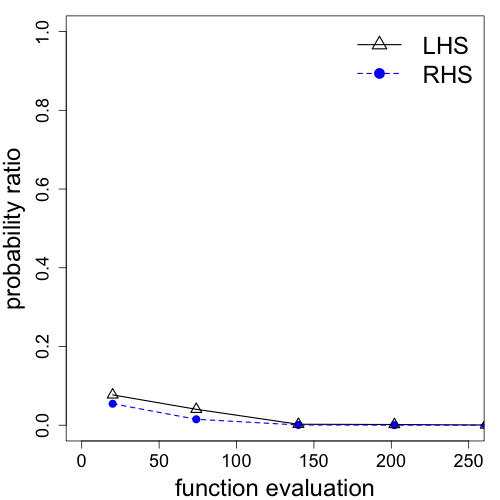}}
            & {\includegraphics{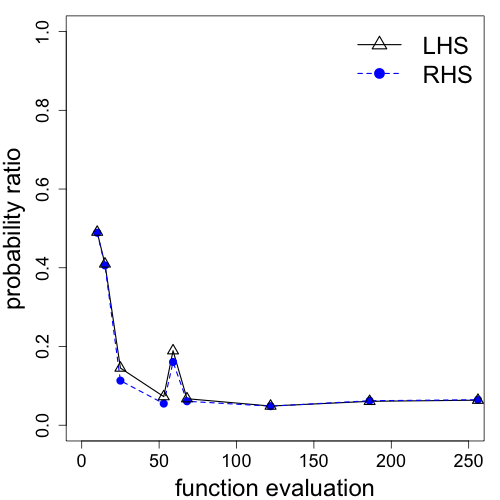}}
                & {\includegraphics{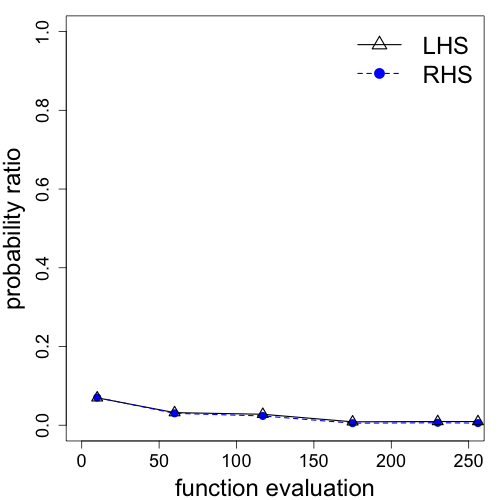}}    \\
\centering
C Quadratic Regression & {\includegraphics{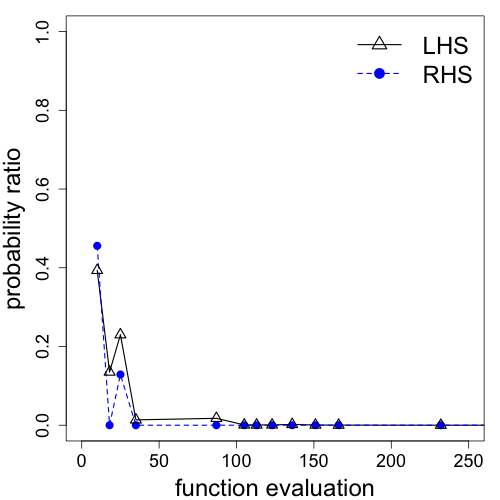}}
        & {\includegraphics{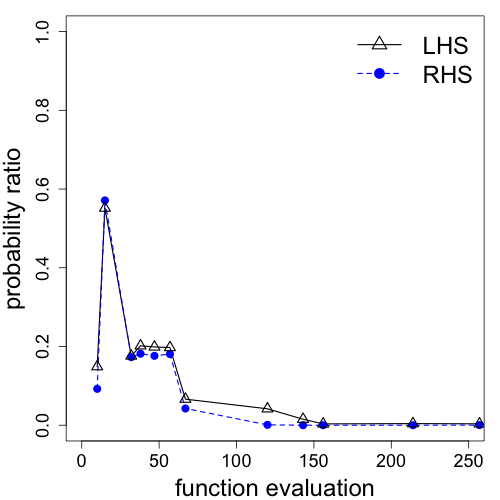}}
            & {\includegraphics{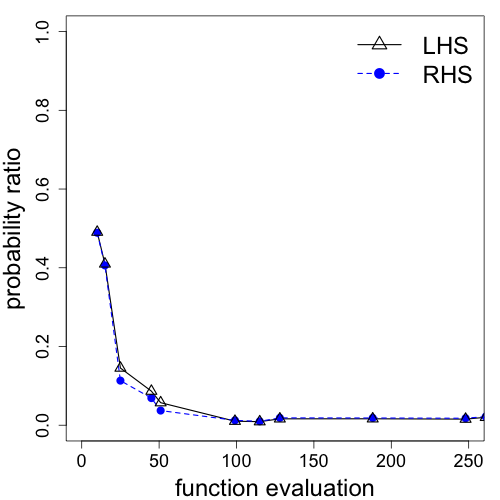}}
                & {\includegraphics{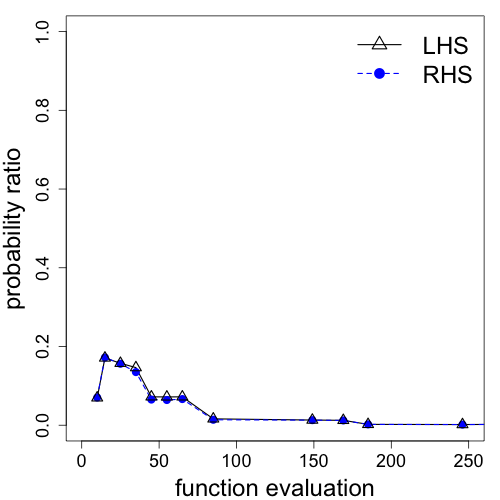}}\\
    \bottomrule
    \end{tblr}
\caption{
BASSO variations (Ab, Bb, Cb) with adaptive subregion probabilities in \eqref{eq:ptildesamplevariance}, sample variance (b). 
When the left-hand side (LHS) of the probability ratio in \eqref{ratioassumption}
with $t_z$ as the incumbent function value is  larger than the right-hand side (RHS)
with $t_{z'}$ as the 20\% quantile of the observed objective function values, $t_z < t_{z'}$,   Assumption 1.3 is satisfied.}
\label{fig:table_of_figures_b}
\end{figure}

\begin{figure}
    \centering
    \setkeys{Gin}{width=\linewidth}
    \settowidth\rotheadsize{PARAMETERS 3}    % from makecell
\begin{tblr}{colspec = { Q[h] *{4}{Q[c,m, wd=30mm]}},
             colsep  = 1pt,
             cell{2-Z}{1} = {cmd=\rotcell, font=\footnotesize\bfseries},
             row{1} = {font=\bfseries},
             measure=vbox
            }
    \toprule
 & Shifted Sinusoidal (dimension=2)& Shifted Sinusoidal (dimension=5)& Rosenbrock (dimension=2)&  Rosenbrock (dimension=5)\\
    \midrule
\centering
A Uniform & {\includegraphics{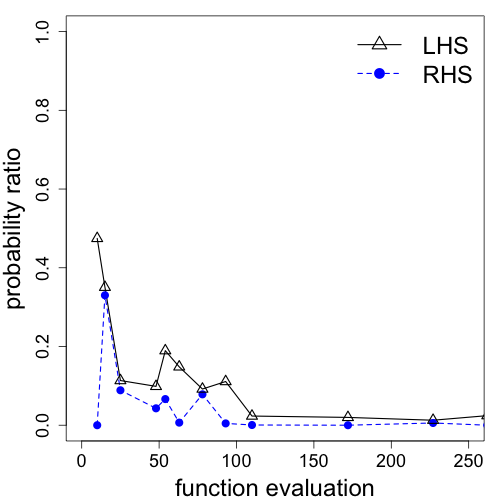}}
        & {\includegraphics{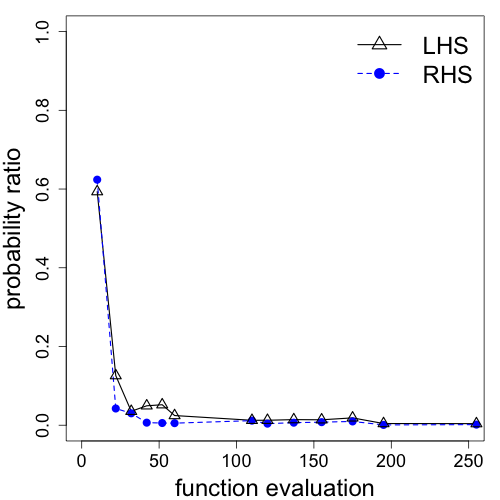}}
            & {\includegraphics{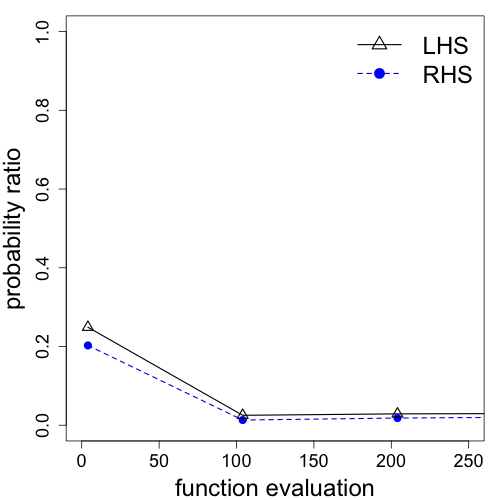}}
                & {\includegraphics{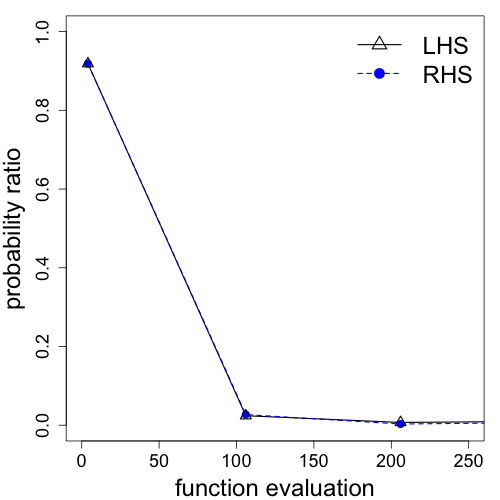}}     \\
\centering
B Gaussian Process & {\includegraphics{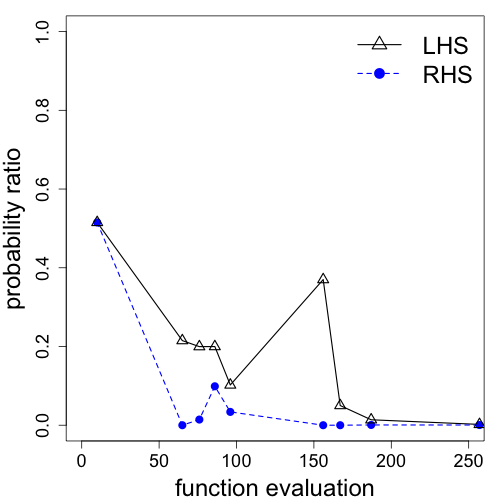}}
        & {\includegraphics{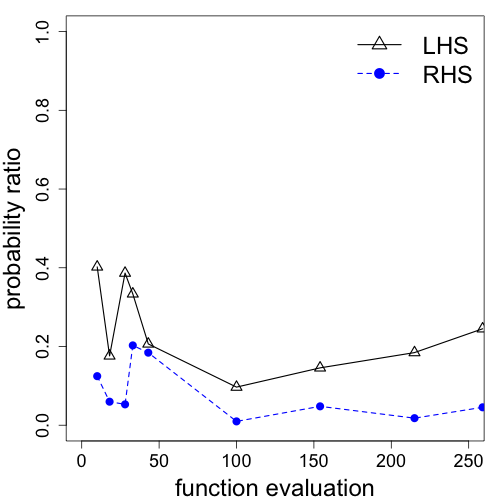}}
            & {\includegraphics{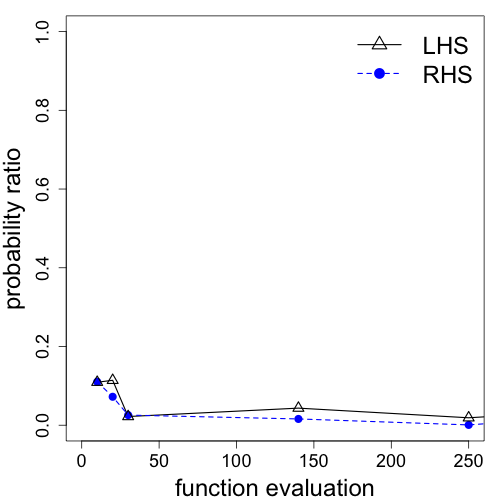}}
                & {\includegraphics{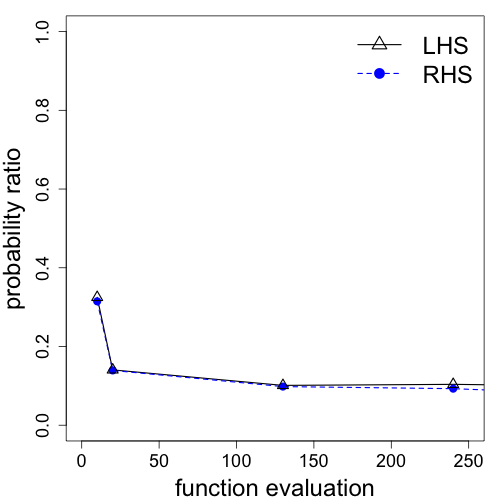}}    \\
\centering
C Quadratic Regression & {\includegraphics{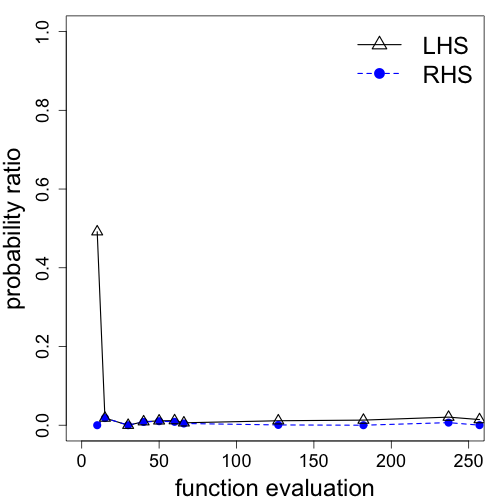}}
        & {\includegraphics{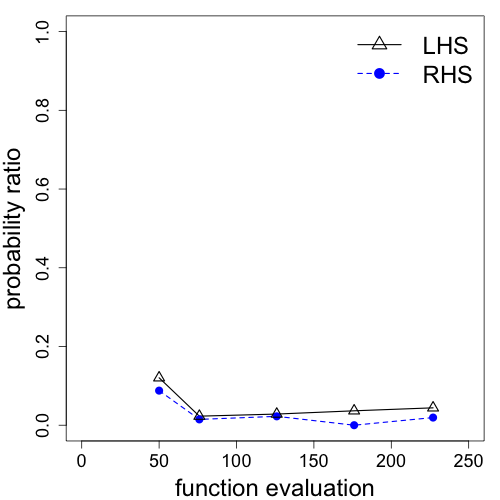}}
            & {\includegraphics{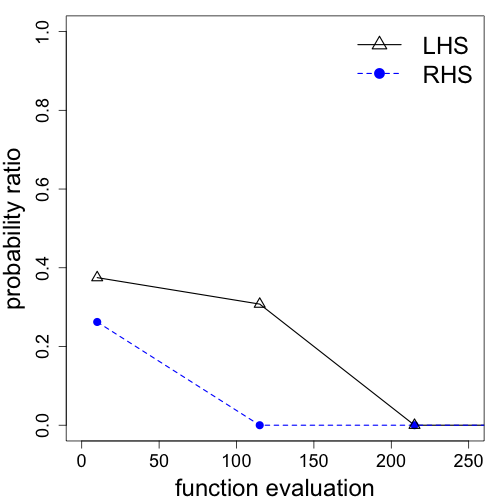}}
                & {\includegraphics{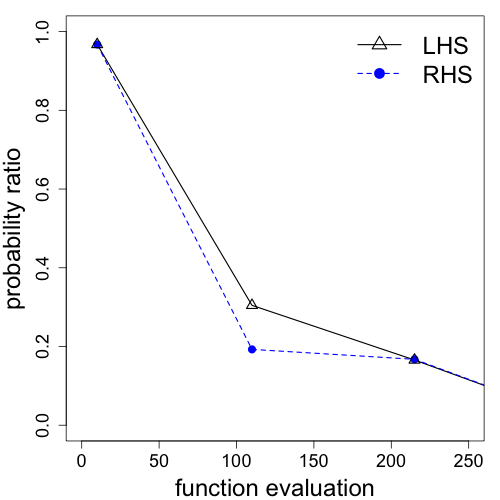}}\\
    \bottomrule
    \end{tblr}
\caption{
BASSO variations (Ac, Bc, Cc) with adaptive subregion probabilities in \eqref{eq:piC}, Gaussian process on the range distribution (c). 
When the left-hand side (LHS) of the probability ratio in \eqref{ratioassumption}
with $t_z$ as the incumbent function value is  larger than the right-hand side (RHS)
with $t_{z'}$ as the 20\% quantile of the observed objective function values, $t_z < t_{z'}$,   Assumption 1.3 is satisfied.
}
\label{fig:table_of_figures_c}
\end{figure}

\begin{figure}
    \centering
    \setkeys{Gin}{width=\linewidth}
    \settowidth\rotheadsize{PARAMETERS 3}    % from makecell
\begin{tblr}{colspec = { Q[h] *{4}{Q[c,m, wd=30mm]}},
             colsep  = 1pt,
             cell{2-Z}{1} = {cmd=\rotcell, font=\footnotesize\bfseries},
             row{1} = {font=\bfseries},
             measure=vbox
            }
    \toprule
& Shifted Sinusoidal (dimension=2)& Shifted Sinusoidal (dimension=5)& Rosenbrock (dimension=2)&  Rosenbrock (dimension=5)\\
    \midrule
\centering
A Uniform & {\includegraphics{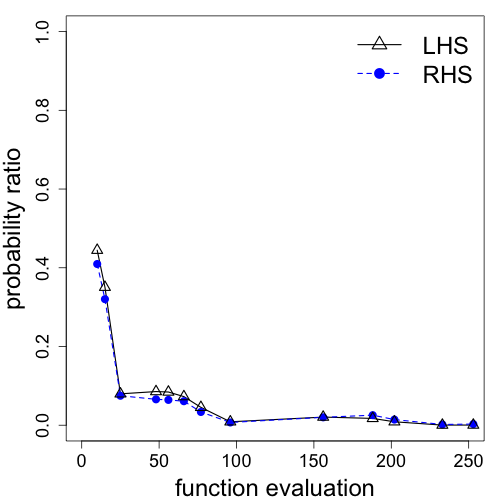}}
        & {\includegraphics{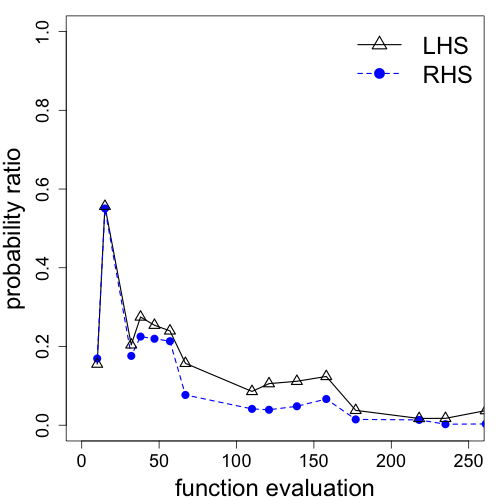}}
            & {\includegraphics{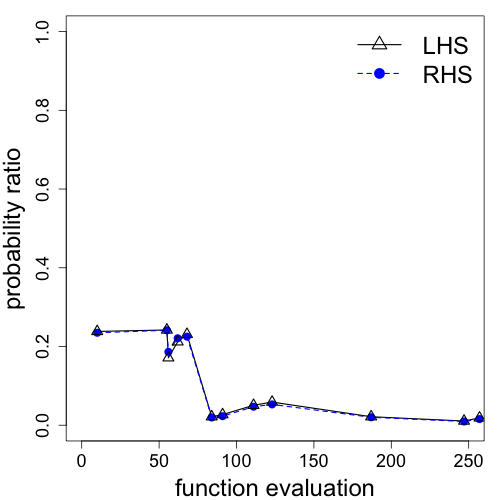}}
                & {\includegraphics{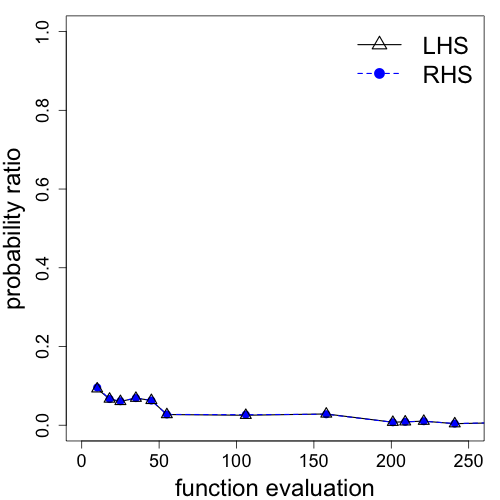}}     \\
\centering
B Gaussian Process & 
{\includegraphics{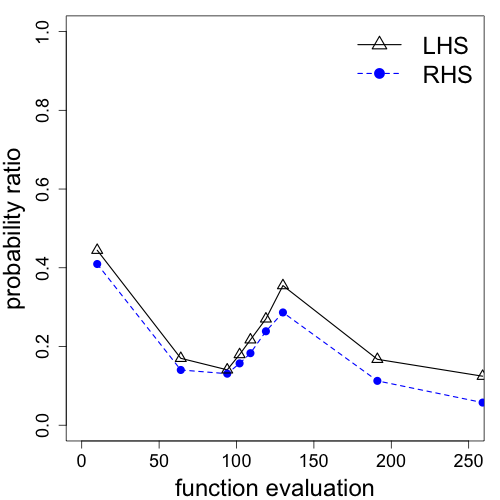}}
        & {\includegraphics{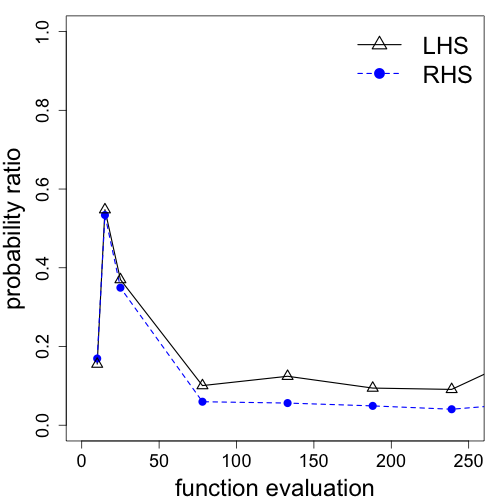}}
            & {\includegraphics{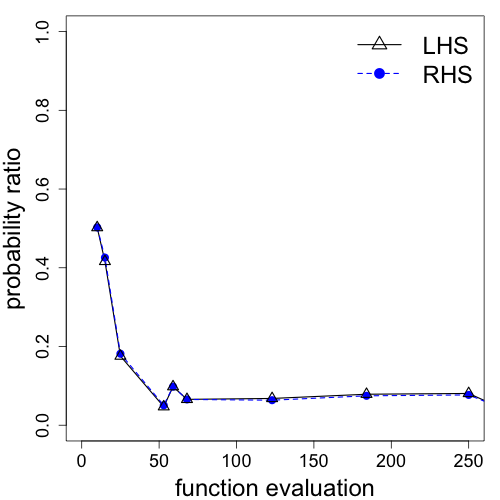}}
                & {\includegraphics{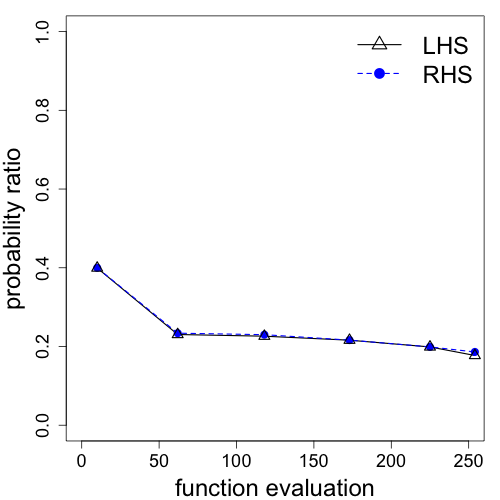}}    \\
\centering
C Quadratic Regression & {\includegraphics{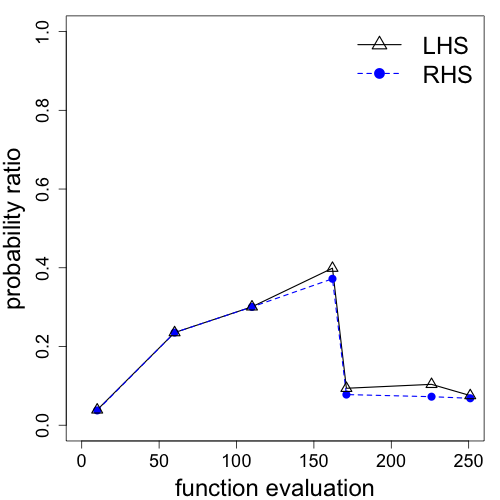}}
        & {\includegraphics{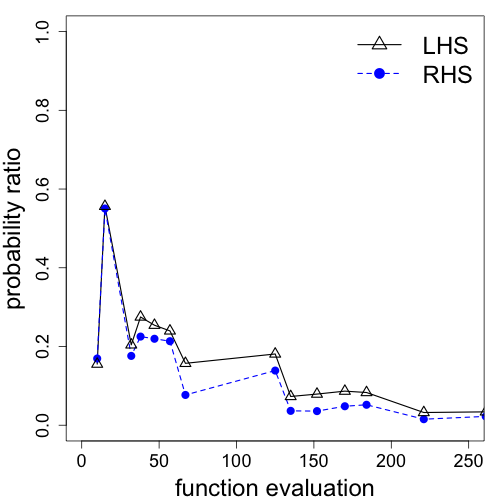}}
            & {\includegraphics{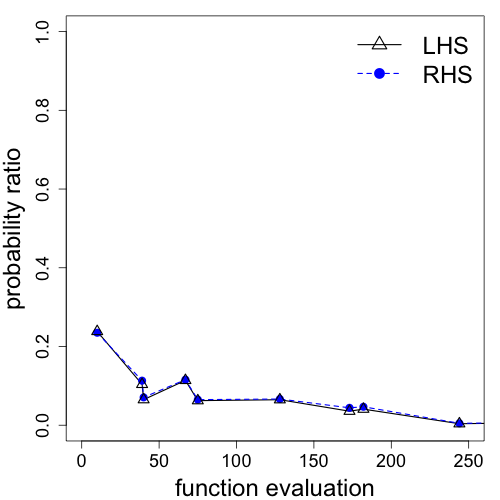}}
                & {\includegraphics{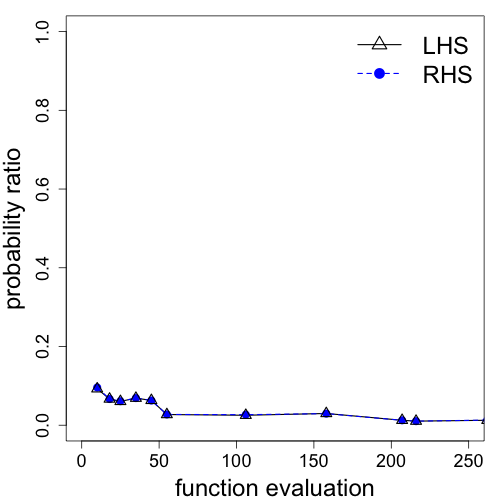}}\\
    \bottomrule
    \end{tblr}
\caption{
BASSO variations (Ad, Bd, Cd) with adaptive subregion probabilities in \eqref{eq:ptildeUBLB}, confidence bounds (d). 
When the left-hand side (LHS) of the probability ratio in \eqref{ratioassumption}
with $t_z$ as the incumbent function value is  larger than the right-hand side (RHS)
with $t_{z'}$ as the 20\% quantile of the observed objective function values, $t_z < t_{z'}$,   Assumption 1.3 is satisfied.
}
\label{fig:table_of_figures_d}
\end{figure}

Figures \ref{fig:table_of_figures_a} - \ref{fig:table_of_figures_d} illustrate the 
 left-hand side (LHS) of the probability ratio in \eqref{ratioassumption}
and the right-hand side (RHS) for adaptive subregion probabilities a-d. When the LHS is greater than the RHS, Assumption~1.3 is satisfied.

The graphs are constructed 
with $t_z$ as the incumbent function value and
with $t_{z'}$ as the 20\% quantile of the observed objective function values, $t_z < t_{z'}$.  The numerical values associated with the graphs are provided in Tables~\ref{table:table_assumption1.3_a} -  \ref{table:table_assumption1.3_d}. The iterations where the LHS is less than the RHS are highlighted, indicating a violation of    Assumption~1.3.

\begin{landscape}
\begin{table}[ht]
\fontsize{6pt}{6pt}\selectfont
\addtolength{\tabcolsep}{-2pt} % reduce space by -2pt
\centering
\begin{tabular}{lllll|lllll|lllll|lllll}
\hline
\multicolumn{20}{c}{\textbf{A Uniform}}\\
\hline
\multicolumn{5}{c|}{Shifted sinusoidal $d=2$}
&\multicolumn{5}{c|}{Shifted sinusoidal $d=5$}
& \multicolumn{5}{c|}{Rosenbrock $d=2$}
& \multicolumn{5}{c}{Rosenbrock $d=5$}
\\\hline
\multicolumn{1}{c}{$k$}& \multicolumn{1}{c}{LHS} & \multicolumn{1}{c}{RHS} & \multicolumn{1}{c}{$y$} & \multicolumn{1}{c}{$z$} 
& \multicolumn{1}{c}{$k$}& \multicolumn{1}{c}{LHS} & \multicolumn{1}{c}{RHS} & \multicolumn{1}{c}{$y$} & \multicolumn{1}{c}{$z$} 
& \multicolumn{1}{c}{$k$}& \multicolumn{1}{c}{LHS} & \multicolumn{1}{c}{RHS} & \multicolumn{1}{c}{$y$} & \multicolumn{1}{c}{$z$} 
& \multicolumn{1}{c}{$k$}& \multicolumn{1}{c}{LHS} & \multicolumn{1}{c}{RHS} & \multicolumn{1}{c}{$y$} & \multicolumn{1}{c}{$z$} \\ 
\hline
  10 & 0.415 & 0.348 & 1.63 & 2.769 & 
  10 & {\hl {0.177}} &{\hl { 0.183}} & 2.811 & 3.411 & 
  10 & 0.867 & 0.867 & 185.542 & 232.845 & 
  10 & 0.721 & 0.721 & 642.486 & 840.008 \\ 
  18 & 0.173 & 0.174 & 0.998 & 2.229 & 
  15 & 0.544 & 0.517 & 2.67 & 2.961 & 
  110 & 0.012 & 0 & 1.256 & 29.326 & 
  114 & {\hl {0.154}} & {\hl {0.157}} & 198.088 & 791.735 \\ 
  26 & {\hl {0.233}} & {\hl {0.251}} & 0.998 & 1.993 & 
  32 & 0.165 & 0.148 & 2.594 & 3.347 & 
  212 & 0.019 & 0.129 & 1.256 & 41.415 & 
  214 & 0.071 & 0.068 & 139.612 & 763.788 \\ 
  35 & {\hl {0.268}} & {\hl {0.293}} & 0.998 & 1.772 &
  38 & 0.21 & 0.185 & 2.594 & 3.287 & 
  312 & 0.059 & 0.042 & 1.256 & 37.179 & 
  319 & 0.004 & 0.004 & 33.084 & 718.587 \\ 
  45 & {\hl {0.267}} & {\hl {0.29}} & 0.998 & 1.772 & 
  47 & 0.212 & 0.182 & 2.594 & 3.284 &  &  &  &  &  &  &  &  &  &  \\ 
  55 & {\hl {0.32}} & {\hl {0.344}} & 0.998 & 1.63 & 
  57 & 0.209 & 0.181 & 2.594 & 3.287 &  &  &  &  &  &  &  &  &  &  \\ 
  65 & {\hl {0.461}} & {\hl {0.47}} & 0.998 & 1.332 & 
  67 & 0.08 & 0.068 & 1.92 & 3.284 &  &  &  &  &  &  &  &  &  &  \\ 
  85 & 0.054 & 0.053 & 0.095 & 1.233 &
  120 & {\hl {0.1}} & {\hl {0.11}} & 1.92 & 3.181 &  &  &  &  &  &  &  &  &  &  \\ 
  129 & 0.004 & 0.004 & 0.008 & 1.165 & 
  142 & {\hl {0.032}} & {\hl {0.033}} & 1.505 & 3.06 &  &  &  &  &  &  &  &  &  &  \\ 
  141 & 0.004 & 0.004 & 0.008 & 1.162 & 
  159 & 0.005 & 0.004 & 1.168 & 2.978 &  &  &  &  &  &  &  &  &  &  \\ 
  161 & 0.004 & 0.005 & 0.008 & 1.108 & 
  217 & 0.005 & 0.005 & 1.168 & 2.958 &  &  &  &  &  &  &  &  &  &  \\ 
  180 & 0.008 & 0.008 & 0.008 & 0.909 & 
  258 & 0.008 & 0.008 & 1.168 & 2.738 &  &  &  &  &  &  &  &  &  &  \\ 
  200 & 0.011 & 0.011 & 0.008 & 0.829 &  &  &  &  &  &  &  &  &  &  &  &  &  &  &  \\ 
  219 & 0.006 & 0.006 & 0.003 & 0.527 &  &  &  &  &  &  &  &  &  &  &  &  &  &  &  \\ 
  265 & {\hl {0.011}} & {\hl {0.012}} & 0.003 & 0.254 &  &  &  &  &  &  &  &  &  &  &  &  &  &  &  \\ 
   &  &  &  &  &  &  &  &  &  &  &  &  &  &  &  &  &  &  &  \\ 
   \hline
\multicolumn{20}{c}{\textbf{B Gaussian Process}}\\ 
   \hline
\multicolumn{5}{c|}{Shifted sinusoidal $d=2$}
&\multicolumn{5}{c|}{Shifted sinusoidal $d=5$}
& \multicolumn{5}{c|}{Rosenbrock $d=2$}
& \multicolumn{5}{c}{Rosenbrock $d=5$}
\\\hline
\multicolumn{1}{c}{$k$}& \multicolumn{1}{c}{LHS} & \multicolumn{1}{c}{RHS} & \multicolumn{1}{c}{$y$} & \multicolumn{1}{c}{$z$} 
& \multicolumn{1}{c}{$k$}& \multicolumn{1}{c}{LHS} & \multicolumn{1}{c}{RHS} & \multicolumn{1}{c}{$y$} & \multicolumn{1}{c}{$z$} 
& \multicolumn{1}{c}{$k$}& \multicolumn{1}{c}{LHS} & \multicolumn{1}{c}{RHS} & \multicolumn{1}{c}{$y$} & \multicolumn{1}{c}{$z$} 
& \multicolumn{1}{c}{$k$}& \multicolumn{1}{c}{LHS} & \multicolumn{1}{c}{RHS} & \multicolumn{1}{c}{$y$} & \multicolumn{1}{c}{$z$} \\ 
\hline
  10 & 0.415 & 0.348 & 1.63 & 2.769 & 
  10 &{\hl {0.466}} & {\hl {0.474}} & 3.286 & 3.5 & 
  10 & 0.867 & 0.867 & 185.542 & 232.845 & 
  10 & 0.721 & 0.721 & 642.486 & 840.008 \\ 
  64 & {\hl {0.145}} & {\hl {0.148}} & 1.002 & 2.769 & 
  20 & {\hl {0.066}} & {\hl {0.075}} & 2.531 & 3.489 & 
  111 & 0.705 & 0.361 & 158.721 & 645.963 & 
  114 & {\hl {0.104}} & {\hl {0.148}} & 358.04 & 2046.24 \\ 
  116 & {\hl {0.161}} & {\hl {0.181}} & 1.002 & 2.179 & 
  64 & 0.087 & 0.052 & 2.09 & 3.35 & 
  220 & 0.265 & 0.216 & 7.3 & 185.542 &
  217 & 0.073 & 0.134 & 358.04 & 2404.12 \\ 
  170 & 0.165 & 0.165 & 1.002 & 2.251 & 
  67 & 0.091 & 0.037 & 2.09 & 3.358 & 
  327 & {\hl {0.061}} & {\hl {0.063}} & 0.542 & 122.652 & 
  321 & {\hl {0.055}} & {\hl {0.102}} & 358.04 & 2587.64 \\ 
  223 & 0.165 & 0.164 & 1.002 & 2.316 & 
  73 & {\hl {0.019}} & {\hl {0.02}} & 1.443 & 3.325 &  &  &  &  &  &  &  &  &  &  \\ 
  238 & {\hl {0.058}} & {\hl {0.059}} & 0.278 & 2.179 & 
  130 & {\hl {0.022}} & {\hl {0.023}} & 1.415 & 3.136 &  &  &  &  &  &  &  &  &  &  \\ 
  260 & {\hl {0.069}} & {\hl {0.074}} & 0.278 & 1.932 & 
  183 & {\hl {0.022}} & {\hl {0.023}} & 1.415 & 3.131 &  &  &  &  &  &  &  &  &  &  \\ 
   &  &  &  &  & 
   246 & 0.012 & 0.011 & 1.321 & 3.047 &  &  &  &  &  &  &  &  &  &  \\ 
   &  &  &  &  & 
   265 & 0.002 & 0.002 & 1.016 & 2.995 &  &  &  &  &  &  &  &  &  &  \\ 
   &  &  &  &  &  &  &  &  &  &  &  &  &  &  &  &  &  &  &  \\ 
   \hline
\multicolumn{20}{c}{\textbf{C Quadratic Regression}}\\ 
   \hline
\multicolumn{5}{c|}{Shifted sinusoidal $d=2$}
&\multicolumn{5}{c|}{Shifted sinusoidal $d=5$}
& \multicolumn{5}{c|}{Rosenbrock $d=2$}
& \multicolumn{5}{c}{Rosenbrock $d=5$}
\\\hline
\multicolumn{1}{c}{$k$}& \multicolumn{1}{c}{LHS} & \multicolumn{1}{c}{RHS} & \multicolumn{1}{c}{$y$} & \multicolumn{1}{c}{$z$} 
& \multicolumn{1}{c}{$k$}& \multicolumn{1}{c}{LHS} & \multicolumn{1}{c}{RHS} & \multicolumn{1}{c}{$y$} & \multicolumn{1}{c}{$z$} 
& \multicolumn{1}{c}{$k$}& \multicolumn{1}{c}{LHS} & \multicolumn{1}{c}{RHS} & \multicolumn{1}{c}{$y$} & \multicolumn{1}{c}{$z$} 
& \multicolumn{1}{c}{$k$}& \multicolumn{1}{c}{LHS} & \multicolumn{1}{c}{RHS} & \multicolumn{1}{c}{$y$} & \multicolumn{1}{c}{$z$} \\ 
\hline
  10 & 0.415 & 0.348 & 1.63 & 2.769 & 
  10 & {\hl {0.177}} & {\hl {0.183}} & 2.811 & 3.411 & 
  10 & 0.867 & 0.867 & 185.542 & 232.845 & 
  10 & 0.721 & 0.721 & 642.486 & 840.008 \\ 
  18 & {\hl {0.173}} &{\hl { 0.174}} & 0.998 & 2.229 & 
  15 & 0.544 & 0.517 & 2.67 & 2.961 &
  110 & 0.038 & 0.001 & 1.067 & 5.379 & 
  114 & 0.66 & 0.66 & 16.327 & 20.32 \\ 
  26 & {\hl {0.233}} &{\hl { 0.251}} & 0.998 & 1.993 & 
  32 & 0.165 & 0.148 & 2.594 & 3.347 & 
  213 & 1 & 1 & 1.019 & 1.044 & 
  215 & {\hl {0.114}} & {\hl {0.115}} & 16.327 & 43.128 \\ 
  35 & {\hl {0.268}} & {\hl {0.293}} & 0.998 & 1.772 & 
  38 & 0.21 & 0.185 & 2.594 & 3.287 &
  316 & 0 & 0 & 0 & 1.052 & 
  320 & 0.025 & 0.025 & 3.918 & 20.309 \\ 
  45 & {\hl {0.267}} &{\hl { 0.29}} & 0.998 & 1.772 &
  47 & 0.212 & 0.182 & 2.594 & 3.284 &  &  &  &  &  &  &  &  &  &  \\ 
  55 & {\hl {0.32}} & {\hl {0.344}} & 0.998 & 1.63 &
  57 & 0.209 & 0.181 & 2.594 & 3.287 &  &  &  &  &  &  &  &  &  &  \\ 
  65 & {\hl {0.461}} & {\hl {0.47}} & 0.998 & 1.332 & 
  67 & 0.08 & 0.068 & 1.92 & 3.284 &  &  &  &  &  &  &  &  &  &  \\ 
  85 & 0.054 & 0.053 & 0.095 & 1.233 & 
  120 & {\hl {0.1}} & {\hl {0.11}} & 1.92 & 3.181 &  &  &  &  &  &  &  &  &  &  \\ 
  129 & 0.004 & 0.004 & 0.008 & 1.165 &
  142 & {\hl {0.032}} &{\hl { 0.033}} & 1.505 & 3.06 &  &  &  &  &  &  &  &  &  &  \\ 
  141 & 0.004 & 0.004 & 0.008 & 1.162 &
  159 & 0.005 & 0.004 & 1.168 & 2.978 &  &  &  &  &  &  &  &  &  &  \\ 
  161 & {\hl {0.004}} &{\hl { 0.005}} & 0.008 & 1.108 &
  217 & 0.005 & 0.005 & 1.168 & 2.958 &  &  &  &  &  &  &  &  &  &  \\ 
  180 & 0.008 & 0.008 & 0.008 & 0.909 & 
  258 & 0.008 & 0.008 & 1.168 & 2.738 &  &  &  &  &  &  &  &  &  &  \\ 
  200 & 0.011 & 0.011 & 0.008 & 0.829 &  &  &  &  &  &  &  &  &  &  &  &  &  &  &  \\ 
  219 & 0.006 & 0.006 & 0.003 & 0.527 &  &  &  &  &  &  &  &  &  &  &  &  &  &  &  \\ 
  265 & {\hl {0.011}} & {\hl {0.012}} & 0.003 & 0.254 &  &  &  &  &  &  &  &  &  &  &  &  &  &  &  \\ 
   \hline
\end{tabular}
\caption{Numerical values for testing Assumption 1.3  with adaptive subregion probabilities (a) observed best function value.
Highlighted cells violate Assumption 1.3.}
\label{table:table_assumption1.3_a}
\end{table}
\end{landscape}

\begin{landscape}
\begin{table}[ht]
\fontsize{6pt}{6pt}\selectfont
\addtolength{\tabcolsep}{-2pt} % reduce space by -2pt
\centering
\begin{tabular}{lllll|lllll|lllll|lllll}
\hline
\multicolumn{20}{c}{\textbf{A Uniform}}\\
\hline
\multicolumn{5}{c|}{Shifted sinusoidal $d=2$}
&\multicolumn{5}{c|}{Shifted sinusoidal $d=5$}
& \multicolumn{5}{c|}{Rosenbrock $d=2$}
& \multicolumn{5}{c}{Rosenbrock $d=5$}
\\\hline
\multicolumn{1}{c}{$k$}& \multicolumn{1}{c}{LHS} & \multicolumn{1}{c}{RHS} & \multicolumn{1}{c}{$y$} & \multicolumn{1}{c}{$z$} 
& \multicolumn{1}{c}{$k$}& \multicolumn{1}{c}{LHS} & \multicolumn{1}{c}{RHS} & \multicolumn{1}{c}{$y$} & \multicolumn{1}{c}{$z$} 
& \multicolumn{1}{c}{$k$}& \multicolumn{1}{c}{LHS} & \multicolumn{1}{c}{RHS} & \multicolumn{1}{c}{$y$} & \multicolumn{1}{c}{$z$} 
& \multicolumn{1}{c}{$k$}& \multicolumn{1}{c}{LHS} & \multicolumn{1}{c}{RHS} & \multicolumn{1}{c}{$y$} & \multicolumn{1}{c}{$z$} \\ 
\hline
  10 & {\hl {0.393}} & {\hl {0.456}} & 1.63 & 2.769 & 
  10 & 0.149 & 0.092 & 2.811 & 3.411 & 
  10 & 0.491 & 0.489 & 8.298 & 35.649 & 
  10 & 0.07 & 0.07 & 157.848 & 1012.868 \\ 
  18 & 0.135 & 0 & 0.998 & 2.229 & 
  15 & {\hl {0.552}} & {\hl {0.571}} & 2.67 & 2.961 & 
  15 & 0.41 & 0.406 & 8.298 & 35.649 & 
  15 & {\hl {0.171}} & {\hl {0.172}} & 105.59 & 287.408 \\ 
  25 & 0.23 & 0.129 & 0.998 & 1.912 & 
  32 & 0.176 & 0.173 & 2.594 & 3.347 & 
  25 & 0.145 & 0.113 & 3.834 & 35.649 & 
  25 & 0.157 & 0.156 & 105.59 & 305.584 \\ 
  35 & 0.014 & 0 & 0.108 & 1.847 &
  38 & 0.202 & 0.182 & 2.594 & 3.287 & 
  45 & 0.087 & 0.069 & 3.017 & 60.589 & 
  35 & 0.146 & 0.135 & 105.59 & 305.584 \\ 
  87 & 0.018 & 0 & 0.108 & 1.762 & 
  47 & 0.199 & 0.176 & 2.594 & 3.284 & 
  51 & 0.057 & 0.037 & 2.17 & 58.605 & 
  45 & 0.072 & 0.065 & 105.59 & 497.471 \\ 
  105 & 0.001 & 0 & 0.008 & 1.722 & 
  57 & 0.197 & 0.181 & 2.594 & 3.287 & 
  99 & {\hl {0.01}} & {\hl {0.012}} & 0.015 & 57.96 & 
  55 & 0.072 & 0.064 & 105.59 & 509.156 \\ 
  113 & 0.001 & 0 & 0.008 & 1.56 & 
  67 & 0.066 & 0.043 & 1.92 & 3.284 & 
  115 & {\hl {0.009}} & {\hl {0.01}} & 0.015 & 29.507 & 
  65 & 0.072 & 0.066 & 105.59 & 497.471 \\ 
  123 & 0.001 & 0 & 0.008 & 1.108 & 
  120 & 0.042 & 0.001 & 1.92 & 3.272 & 
  128 & {\hl {0.016}} &{\hl { 0.019}} & 0.015 & 14.095 & 
  85 & 0.016 & 0.014 & 43.581 & 384.793 \\ 
  136 & 0.002 & 0 & 0.008 & 0.838 & 
  143 & 0.015 & 0 & 1.505 & 3.082 & 
  188 & {\hl {0.018}} &{\hl { 0.022}} & 0.015 & 12.251 & 
  149 & 0.013 & 0.013 & 40.763 & 385.923 \\ 
  151 & 0 & 0 & 0.008 & 0.521 & 
  156 & 0.003 & 0 & 1.168 & 2.978 & 
  252 & {\hl {0.016}} & {\hl {0.019}} & 0.015 & 14.095 & 
  169 & 0.012 & 0.012 & 39.685 & 384.793 \\ 
  166 & 0 & 0 & 0.002 & 0.361 & 
  214 & 0.004 & 0 & 1.168 & 2.871 &  &  &  &  &  & 
  185 & 0.002 & 0.002 & 18.616 & 340.773 \\ 
  232 & 0 & 0 & 0.001 & 0.398 & 
  257 & 0.004 & 0 & 1.168 & 2.742 &  &  &  &  &  & 
  246 & 0.002 & 0.002 & 18.616 & 379.67 \\ 
  264 & 0 & 0 & 0.001 & 0.348 &  &  &  &  &  &  &  &  &  &  & 
  268 & 0.002 & 0.002 & 18.616 & 332.86 \\ 
   &  &  &  &  &  &  &  &  &  &  &  &  &  &  &  &  &  &  &  \\ 
   \hline
\multicolumn{20}{c}{\textbf{B Gaussian Process}}\\ 
   \hline
\multicolumn{5}{c|}{Shifted sinusoidal $d=2$}
&\multicolumn{5}{c|}{Shifted sinusoidal $d=5$}
& \multicolumn{5}{c|}{Rosenbrock $d=2$}
& \multicolumn{5}{c}{Rosenbrock $d=5$}
\\\hline
\multicolumn{1}{c}{$k$}& \multicolumn{1}{c}{LHS} & \multicolumn{1}{c}{RHS} & \multicolumn{1}{c}{$y$} & \multicolumn{1}{c}{$z$} 
& \multicolumn{1}{c}{$k$}& \multicolumn{1}{c}{LHS} & \multicolumn{1}{c}{RHS} & \multicolumn{1}{c}{$y$} & \multicolumn{1}{c}{$z$} 
& \multicolumn{1}{c}{$k$}& \multicolumn{1}{c}{LHS} & \multicolumn{1}{c}{RHS} & \multicolumn{1}{c}{$y$} & \multicolumn{1}{c}{$z$} 
& \multicolumn{1}{c}{$k$}& \multicolumn{1}{c}{LHS} & \multicolumn{1}{c}{RHS} & \multicolumn{1}{c}{$y$} & \multicolumn{1}{c}{$z$} \\ 
\hline
  10 & {\hl {0.393}} & {\hl {0.456}} & 1.63 & 2.769 & 
  20 & 0.077 & 0.055 & 1.883 & 3.209 & 
  10 & 0.491 & 0.489 & 8.298 & 35.649 & 
  10 & 0.07 & 0.07 & 157.848 & 1012.868 \\ 
  64 & 0.094 & 0 & 1.002 & 2.636 & 
  74 & 0.04 & 0.015 & 1.883 & 3.517 & 
  15 & 0.41 & 0.406 & 8.298 & 35.649 & 
  60 & 0.032 & 0.03 & 157.848 & 2414.04 \\ 
  116 & 0.161 & 0.009 & 1.002 & 2.116 & 
  140 & 0.002 & 0.001 & 1.207 & 3.533 &
  25 & 0.145 & 0.113 & 3.834 & 35.649 & 
  117 & 0.028 & 0.024 & 157.848 & 2207.16 \\ 
  174 & 0.158 & 0.101 & 1.002 & 2.116 & 
  202 & 0.002 & 0 & 1.207 & 3.637 & 
  53 & 0.073 & 0.055 & 3.62 & 72.723 &
  175 & 0.009 & 0.005 & 82.823 & 2085.52 \\ 
  227 & 0.158 & 0.088 & 1.002 & 2.145 & 
  261 & 0.001 & 0 & 1.207 & 3.686 & 
  59 & 0.189 & 0.161 & 3.62 & 35.649 & 
  230 & 0.009 & 0.006 & 82.823 & 1988.828 \\ 
  259 & 0.182 & 0.13 & 1.002 & 2.116 &  &  &  &  &  & 
  68 & 0.067 & 0.061 & 0.416 & 15.787 & 
  256 & 0.009 & 0.005 & 82.823 & 2045.24 \\ 
   &  &  &  &  &  &  &  &  &  & 
   122 & 0.049 & 0.048 & 0.416 & 35.649 &  &  &  &  &  \\ 
   &  &  &  &  &  &  &  &  &  & 
   186 & 0.061 & 0.062 & 0.416 & 33.315 &  &  &  &  &  \\ 
   &  &  &  &  &  &  &  &  &  & 
   256 & {\hl {0.064}} & {\hl {0.065}} & 0.416 & 31.607 &  &  &  &  &  \\ 
   &  &  &  &  &  &  &  &  &  &  &  &  &  &  &  &  &  &  &  \\ 
   \hline
\multicolumn{20}{c}{\textbf{C Quadratic Regression}}\\ 
   \hline
\multicolumn{5}{c|}{Shifted sinusoidal $d=2$}
&\multicolumn{5}{c|}{Shifted sinusoidal $d=5$}
& \multicolumn{5}{c|}{Rosenbrock $d=2$}
& \multicolumn{5}{c}{Rosenbrock $d=5$}
\\\hline
\multicolumn{1}{c}{$k$}& \multicolumn{1}{c}{LHS} & \multicolumn{1}{c}{RHS} & \multicolumn{1}{c}{$y$} & \multicolumn{1}{c}{$z$} 
& \multicolumn{1}{c}{$k$}& \multicolumn{1}{c}{LHS} & \multicolumn{1}{c}{RHS} & \multicolumn{1}{c}{$y$} & \multicolumn{1}{c}{$z$} 
& \multicolumn{1}{c}{$k$}& \multicolumn{1}{c}{LHS} & \multicolumn{1}{c}{RHS} & \multicolumn{1}{c}{$y$} & \multicolumn{1}{c}{$z$} 
& \multicolumn{1}{c}{$k$}& \multicolumn{1}{c}{LHS} & \multicolumn{1}{c}{RHS} & \multicolumn{1}{c}{$y$} & \multicolumn{1}{c}{$z$} \\ 
\hline
  10 & {\hl {0.393}} & {\hl {0.456}} & 1.63 & 2.769 & 
  10 & 0.149 & 0.092 & 2.811 & 3.411 & 
  10 & 0.491 & 0.489 & 8.298 & 35.649 & 
  10 & 0.07 & 0.07 & 157.848 & 1012.868 \\ 
  18 & 0.135 & 0 & 0.998 & 2.229 & 
  15 & {\hl {0.552}}& {\hl {0.571}} & 2.67 & 2.961 & 
  15 & 0.41 & 0.406 & 8.298 & 35.649 & 
  15 & {\hl {0.171}} & {\hl {0.172}} & 105.59 & 287.408 \\ 
  25 & 0.23 & 0.129 & 0.998 & 1.912 & 
  32 & 0.176 & 0.173 & 2.594 & 3.347 & 
  25 & 0.145 & 0.113 & 3.834 & 35.649 & 
  25 & 0.157 & 0.156 & 105.59 & 305.584 \\ 
  35 & 0.014 & 0 & 0.108 & 1.847 & 
  38 & 0.202 & 0.182 & 2.594 & 3.287 & 
  45 & 0.087 & 0.069 & 3.017 & 60.589 & 
  35 & 0.146 & 0.135 & 105.59 & 305.584 \\ 
  87 & 0.018 & 0 & 0.108 & 1.762 &
  47 & 0.199 & 0.176 & 2.594 & 3.284 & 
  51 & 0.057 & 0.037 & 2.17 & 58.605 & 
  45 & 0.072 & 0.065 & 105.59 & 497.471 \\ 
  105 & 0.001 & 0 & 0.008 & 1.722 & 
  57 & 0.197 & 0.181 & 2.594 & 3.287 & 
  99 &{\hl { 0.01}} & {\hl {0.012}} & 0.015 & 57.96 & 
  55 & 0.072 & 0.064 & 105.59 & 509.156 \\ 
  113 & 0.001 & 0 & 0.008 & 1.56 & 
  67 & 0.066 & 0.043 & 1.92 & 3.284 & 
  115 & {\hl {0.009}} & {\hl {0.01}} & 0.015 & 29.507 &
  65 & 0.072 & 0.066 & 105.59 & 497.471 \\ 
  123 & 0.001 & 0 & 0.008 & 1.108 &
  120 & 0.042 & 0.001 & 1.92 & 3.272 & 
  128 & {\hl {0.016}} & {\hl {0.019}} & 0.015 & 14.095 & 
  85 & 0.016 & 0.014 & 43.581 & 384.793 \\ 
  136 & 0.002 & 0 & 0.008 & 0.838 & 
  143 & 0.015 & 0 & 1.505 & 3.082 & 
  188 & {\hl {0.016}} & {\hl {0.018}} & 0.015 & 14.4 & 
  149 & 0.013 & 0.013 & 40.763 & 385.923 \\ 
  151 & 0 & 0 & 0.008 & 0.521 & 
  156 & 0.003 & 0 & 1.168 & 2.978 & 
  248 & {\hl {0.016}} & {\hl {0.018}} & 0.015 & 14.825 & 
  169 & 0.012 & 0.012 & 39.685 & 384.793 \\ 
  166 & 0 & 0 & 0.002 & 0.361 & 
  214 & 0.004 & 0 & 1.168 & 2.871 & 
  261 & 0.02 & 0.022 & 0.015 & 9.111 & 
  185 & 0.002 & 0.002 & 18.616 & 340.773 \\ 
  232 & 0 & 0 & 0.001 & 0.398 & 
  257 & 0.004 & 0 & 1.168 & 2.742 &  &  &  &  &  & 
  246 & 0.002 & 0.002 & 18.616 & 379.67 \\ 
  264 & 0 & 0 & 0.001 & 0.348 &  &  &  &  &  &  &  &  &  &  & 
  268 & 0.002 & 0.002 & 18.616 & 332.86 \\ 
   \hline
\end{tabular}
\caption{Numerical values for testing Assumption 1.3  with with adaptive subregion probabilities (b) sample variance.
Highlighted cells violate Assumption 1.3.
}
\label{table:table_assumption1.3_b}
\end{table}
\end{landscape}

\begin{landscape}
\begin{table}[ht]
\fontsize{6pt}{6pt}\selectfont
\addtolength{\tabcolsep}{-2pt} % reduce space by -2pt
\centering
\begin{tabular}{lllll|lllll|lllll|lllll}
\hline
\multicolumn{20}{c}{\textbf{A Uniform}}\\
\hline
\multicolumn{5}{c|}{Shifted sinusoidal $d=2$}
&\multicolumn{5}{c|}{Shifted sinusoidal $d=5$}
& \multicolumn{5}{c|}{Rosenbrock $d=2$}
& \multicolumn{5}{c}{Rosenbrock $d=5$}
\\\hline
\multicolumn{1}{c}{$k$}& \multicolumn{1}{c}{LHS} & \multicolumn{1}{c}{RHS} & \multicolumn{1}{c}{$y$} & \multicolumn{1}{c}{$z$} 
& \multicolumn{1}{c}{$k$}& \multicolumn{1}{c}{LHS} & \multicolumn{1}{c}{RHS} & \multicolumn{1}{c}{$y$} & \multicolumn{1}{c}{$z$} 
& \multicolumn{1}{c}{$k$}& \multicolumn{1}{c}{LHS} & \multicolumn{1}{c}{RHS} & \multicolumn{1}{c}{$y$} & \multicolumn{1}{c}{$z$} 
& \multicolumn{1}{c}{$k$}& \multicolumn{1}{c}{LHS} & \multicolumn{1}{c}{RHS} & \multicolumn{1}{c}{$y$} & \multicolumn{1}{c}{$z$} \\ 
\hline
  10 & 0.475 & 0 & 1.63 & 2.769 & 
  10 & {\hl {0.593}} & {\hl {0.624}} & 3.429 & 3.583 & 
  4 & 0.249 & 0.203 & 0.971 & 4.073 & 
  4 & 0.918 & 0.918 & 1278.517 & 1398.628 \\ 
  15 & 0.351 & 0.33 & 1.259 & 2.442 & 
  22 & 0.126 & 0.042 & 2.068 & 3.443 & 
  104 & 0.025 & 0.013 & 0.167 & 27.401 & 
  106 & {\hl {0.024}} & {\hl {0.027}} & 68.333 & 648.264 \\ 
  25 & 0.113 & 0.089 & 0.63 & 2.15 & 
  32 & 0.035 & 0.03 & 1.513 & 3.331 & 
  204 & 0.029 & 0.018 & 0.167 & 41.594 & 
  206 & 0.007 & 0.003 & 30.144 & 636.586 \\ 
  48 & 0.099 & 0.043 & 0.229 & 1.685 & 
  42 & 0.049 & 0.006 & 1.513 & 3.246 & 
  304 & 0.030 & 0.021 & 0.167 & 46.290 & 
  308 & 0.010 & 0.007 & 30.144 & 623.913 \\ 
  54 & 0.189 & 0.066 & 0.229 & 1.685 & 
  52 & 0.052 & 0.006 & 1.513 & 3.246 &  &  &  &  &  &  &  &  &  &  \\ 
  63 & 0.148 & 0.007 & 0.096 & 1.4 & 
  60 & 0.025 & 0.005 & 1.306 & 2.972 &  &  &  &  &  &  &  &  &  &  \\ 
  78 & 0.092 & 0.078 & 0.053 & 0.788 & 
  110 & 0.013 & 0.012 & 1.115 & 2.33 &  &  &  &  &  &  &  &  &  &  \\ 
  93 & 0.111 & 0.005 & 0.053 & 0.653 & 
  120 & 0.013 & 0.004 & 1.115 & 2.395 &  &  &  &  &  &  &  &  &  &  \\ 
  110 & 0.023 & 0.001 & 0.006 & 0.512 & 
  137 & 0.014 & 0.006 & 1.115 & 2.33 &  &  &  &  &  &  &  &  &  &  \\ 
  172 & 0.02 & 0 & 0.006 & 0.535 & 
  155 & 0.014 & 0.008 & 1.115 & 2.249 &  &  &  &  &  &  &  &  &  &  \\ 
  227 & 0.013 & 0.006 & 0.006 & 0.512 & 
  175 & 0.019 & 0.01 & 1.115 & 2.1 &  &  &  &  &  &  &  &  &  &  \\ 
  262 & 0.024 & 0 & 0.006 & 0.377 & 
  195 & 0.004 & 0.001 & 0.708 & 1.967 &  &  &  &  &  &  &  &  &  &  \\ 
   &  &  &  &  & 
   255 & 0.004 & 0.002 & 0.708 & 1.967 &  &  &  &  &  &  &  &  &  &  \\ 
   \hline
\multicolumn{20}{c}{\textbf{B Gaussian Process}}\\ 
   \hline
\multicolumn{5}{c|}{Shifted sinusoidal $d=2$}
&\multicolumn{5}{c|}{Shifted sinusoidal $d=5$}
& \multicolumn{5}{c|}{Rosenbrock $d=2$}
& \multicolumn{5}{c}{Rosenbrock $d=5$}
\\\hline
\multicolumn{1}{c}{$k$}& \multicolumn{1}{c}{LHS} & \multicolumn{1}{c}{RHS} & \multicolumn{1}{c}{$y$} & \multicolumn{1}{c}{$z$} 
& \multicolumn{1}{c}{$k$}& \multicolumn{1}{c}{LHS} & \multicolumn{1}{c}{RHS} & \multicolumn{1}{c}{$y$} & \multicolumn{1}{c}{$z$} 
& \multicolumn{1}{c}{$k$}& \multicolumn{1}{c}{LHS} & \multicolumn{1}{c}{RHS} & \multicolumn{1}{c}{$y$} & \multicolumn{1}{c}{$z$} 
& \multicolumn{1}{c}{$k$}& \multicolumn{1}{c}{LHS} & \multicolumn{1}{c}{RHS} & \multicolumn{1}{c}{$y$} & \multicolumn{1}{c}{$z$} \\ 
\hline
  10 & 0.515 & 0.515 & 1.209 & 1.778 & 
  10 & 0.402 & 0.125 & 2.727 & 3.222 & 
  10 & {\hl {0.109}} & {\hl {0.110}} & 1.382 & 23.492 & 
  10 & 0.325 & 0.315 & 196.646 & 416.047 \\ 
  65 & 0.215 & 0 & 1.209 & 4.479 &
  18 & 0.176 & 0.06 & 2.627 & 3.29 &
  20 & 0.114 & 0.072 & 1.382 & 23.492 & 
  20 & 0.141 & 0.139 & 120.619 & 416.047 \\ 
  76 & 0.2 & 0.014 & 1.1 & 3.229 & 
  28 & 0.386 & 0.053 & 2.627 & 3.322 & 
  30 & {\hl {0.022}} & {\hl {0.026}} & 0.102 & 34.613 & 
  130 & 0.101 & 0.098 & 91.088 & 404.833 \\ 
  86 & 0.2 & 0.099 & 1.1 & 2.524 &
  33 & 0.334 & 0.203 & 2.627 & 3.322 & 
  140 & 0.043 & 0.016 & 0.102 & 47.682 & 
  240 & 0.104 & 0.093 & 91.088 & 400.131 \\ 
  96 & 0.102 & 0.034 & 0.242 & 2.163 & 
  43 & 0.206 & 0.184 & 2.441 & 3.29 & 
  250 & 0.019 & 0.001 & 0.102 & 31.743 &
  350 & 0.098 & 0.076 & 91.088 & 400.131 \\ 
  156 & 0.37 & 0 & 0.242 & 2.516 & 
  100 & 0.097 & 0.01 & 1.549 & 3.369 & 
  360 & 0.045 & 0.032 & 0.102 & 17.118 &  &  &  &  &  \\ 
  167 & 0.05 & 0 & 0.114 & 1.864 &
  154 & 0.146 & 0.048 & 1.549 & 3.238 &  &  &  &  &  &  &  &  &  &  \\ 
  187 & 0.014 & 0.001 & 0.009 & 1.745 & 
  215 & 0.185 & 0.018 & 1.549 & 3.14 &  &  &  &  &  &  &  &  &  &  \\ 
  257 & 0.002 & 0.001 & 0.009 & 1.825 &
  259 & 0.245 & 0.046 & 1.549 & 3.039 &  &  &  &  &  &  &  &  &  &  \\ 
\hline
\multicolumn{20}{c}{\textbf{C Quadratic Regression}}\\ 
\hline
\multicolumn{5}{c|}{Shifted sinusoidal $d=2$}
&\multicolumn{5}{c|}{Shifted sinusoidal $d=5$}
& \multicolumn{5}{c|}{Rosenbrock $d=2$}
& \multicolumn{5}{c}{Rosenbrock $d=5$}
\\\hline
\multicolumn{1}{c}{$k$}& \multicolumn{1}{c}{LHS} & \multicolumn{1}{c}{RHS} & \multicolumn{1}{c}{$y$} & \multicolumn{1}{c}{$z$} 
& \multicolumn{1}{c}{$k$}& \multicolumn{1}{c}{LHS} & \multicolumn{1}{c}{RHS} & \multicolumn{1}{c}{$y$} & \multicolumn{1}{c}{$z$} 
& \multicolumn{1}{c}{$k$}& \multicolumn{1}{c}{LHS} & \multicolumn{1}{c}{RHS} & \multicolumn{1}{c}{$y$} & \multicolumn{1}{c}{$z$} 
& \multicolumn{1}{c}{$k$}& \multicolumn{1}{c}{LHS} & \multicolumn{1}{c}{RHS} & \multicolumn{1}{c}{$y$} & \multicolumn{1}{c}{$z$} \\ 
\hline
  10 & 0.49 & 0 & 1.65 & 2.77 & 
  50 & 0.121 & 0.088 & 2.077 & 3.246 & 
  10 & 0.375 & 0.262 & 1.006 & 5.158 & 
  10 & 0.967 & 0.967 & 1143.565 & 1180.265 \\ 
  15 & 0.02 & 0.02 & 0.15 & 1.65 & 
  76 & 0.023 & 0.015 & 1.574 & 3.246 & 
  115 & 0.308 & 0.000 & 0.453 & 3.366 &
  110 & 0.305 & 0.193 & 14.028 & 25.120 \\ 
  30 & 0 & 0 & 0.03 & 1.19 & 
  126 & 0.028 & 0.023 & 1.574 & 3.26 & 
  215 & 0.000 & 0.000 & 0.000 & 1.046 & 
  215 & {\hl {0.166}} & {\hl {0.167}} & 12.824 & 27.317 \\ 
  40 & 0.01 & 0.01 & 0.03 & 1.26 & 
  176 & 0.037 & 0 & 1.437 & 3.265 & 
  315 & 0.000 & 0.000 & 0.000 & 1.009 & 
  315 & 0.022 & 0.022 & 3.992 & 20.138 \\ 
  50 & 0.01 & 0.01 & 0.03 & 1.23 & 
  227 & 0.044 & 0.019 & 1.423 & 3.158 &  &  &  &  &  &  &  &  &  &  \\ 
  60 & 0.01 & 0.01 & 0.03 & 1.19 &  &  &  &  &  &  &  &  &  &  &  &  &  &  &  \\ 
  66 & 0.01 & 0 & 0.01 & 0.9 &  &  &  &  &  &  &  &  &  &  &  &  &  &  &  \\ 
  127 & 0.01 & 0 & 0.01 & 0.81 &  &  &  &  &  &  &  &  &  &  &  &  &  &  &  \\ 
  182 & 0.01 & 0 & 0.01 & 0.77 &  &  &  &  &  &  &  &  &  &  &  &  &  &  &  \\ 
  237 & 0.02 & 0.01 & 0.01 & 0.74 &  &  &  &  &  &  &  &  &  &  &  &  &  &  &  \\ 
  257 & 0.01 & 0 & 0.01 & 0.62 &  &  &  &  &  &  &  &  &  &  &  &  &  &  &  \\ 
   \hline
\end{tabular}
\caption{Numerical values for testing Assumption 1.3  with adaptive subregion probabilities (c) Gaussian process on the range.
Highlighted cells violate Assumption 1.3.
}
\label{table:table_assumption1.3_c}
\end{table}
\end{landscape}

\begin{landscape}
\begin{table}[ht]
\fontsize{6pt}{6pt}\selectfont
\addtolength{\tabcolsep}{-2pt} % reduce space by -2pt
\centering
\begin{tabular}{lllll|lllll|lllll|lllll}
\hline
\multicolumn{20}{c}{\textbf{A Uniform}}\\
\hline
\multicolumn{5}{c|}{Shifted sinusoidal $d=2$}
&\multicolumn{5}{c|}{Shifted sinusoidal $d=5$}
& \multicolumn{5}{c|}{Rosenbrock $d=2$}
& \multicolumn{5}{c}{Rosenbrock $d=5$}
\\\hline
\multicolumn{1}{c}{$k$}& \multicolumn{1}{c}{LHS} & \multicolumn{1}{c}{RHS} & \multicolumn{1}{c}{$y$} & \multicolumn{1}{c}{$z$} 
& \multicolumn{1}{c}{$k$}& \multicolumn{1}{c}{LHS} & \multicolumn{1}{c}{RHS} & \multicolumn{1}{c}{$y$} & \multicolumn{1}{c}{$z$} 
& \multicolumn{1}{c}{$k$}& \multicolumn{1}{c}{LHS} & \multicolumn{1}{c}{RHS} & \multicolumn{1}{c}{$y$} & \multicolumn{1}{c}{$z$} 
& \multicolumn{1}{c}{$k$}& \multicolumn{1}{c}{LHS} & \multicolumn{1}{c}{RHS} & \multicolumn{1}{c}{$y$} & \multicolumn{1}{c}{$z$} \\ 
\hline
  10 & 0.445 & 0.41 & 1.63 & 2.769 & 
  10 & {\hl {0.155}} & {\hl {0.169}} & 2.811 & 3.411 & 
  10 & 0.238 & 0.235 & 7.63 & 123.944 & 
  10 & {\hl {0.093}} & {\hl {0.096}} & 319.605 & 2467.76 \\ 
  15 & 0.351 & 0.32 & 1.259 & 2.442 & 
  15 & 0.556 & 0.551 & 2.67 & 2.961 & 
  55 & 0.242 & 0.241 & 6.874 & 104.011 & 18 & 0.067 & 0.067 & 143.092 & 854.155 \\ 
  25 & 0.08 & 0.075 & 0.63 & 2.15 & 
  32 & 0.204 & 0.176 & 2.594 & 3.347 & 
  56 & {\hl {0.171}} & {\hl {0.186}} & 6.874 & 104.011 & 
  25 & 0.061 & 0.061 & 143.092 & 943.599 \\ 
  48 & 0.086 & 0.066 & 0.229 & 1.319 & 
  38 & 0.275 & 0.225 & 2.594 & 3.287 & 
  62 & {\hl {0.212}} & {\hl {0.221}} & 6.44 & 100.838 & 
  35 & 0.069 & 0.069 & 143.092 & 854.155 \\ 
  56 & 0.084 & 0.064 & 0.229 & 1.319 & 
  47 & 0.254 & 0.22 & 2.594 & 3.284 & 
  68 & 0.231 & 0.225 & 4.275 & 57.553 & 
  45 & 0.063 & 0.063 & 143.092 & 937.713 \\ 
  66 & 0.073 & 0.061 & 0.172 & 1.259 &
  57 & 0.24 & 0.213 & 2.594 & 3.287 & 
  84 & 0.022 & 0.02 & 0.196 & 62.788 & 
  55 & 0.027 & 0.027 & 88.036 & 854.155 \\ 
  77 & 0.045 & 0.034 & 0.053 & 1.028 & 
  67 & 0.157 & 0.077 & 1.92 & 3.284 & 
  91 & 0.027 & 0.023 & 0.196 & 57.553 &
  106 & {\hl {0.025}} &{\hl {0.027}} & 88.036 & 830.225 \\ 
  96 & 0.009 & 0.007 & 0.008 & 0.829 & 
  110 & 0.086 & 0.041 & 1.51 & 3.038 &
  111 & 0.051 & 0.047 & 0.196 & 18.663 & 
  158 & 0.029 & 0.028 & 88.036 & 873.822 \\ 
  156 & 0.021 & 0.02 & 0.008 & 0.353 & 
  121 & 0.106 & 0.039 & 1.51 & 3.078 &
  123 & 0.059 & 0.053 & 0.07 & 16.368 & 
  201 & 0.008 & 0.008 & 45.281 & 827.083 \\ 
  188 & {\hl{0.017}} & {\hl{0.026}} & 0.008 & 0.267 & 
  139 & 0.112 & 0.048 & 1.51 & 2.98 & 
  187 & 0.022 & 0.02 & 0.07 & 15.322 & 
  209 & 0.009 & 0.009 & 45.281 & 740.015 \\ 
  202 & 0.009 & 0.014 & 0.004 & 0.248 &
  158 & 0.124 & 0.067 & 1.51 & 2.811 & 
  247 & 0.011 & 0.01 & 0.003 & 15.322 &
  221 & 0.01 & 0.01 & 45.281 & 645.793 \\ 
  233 & 0 & 0.002 & 0.001 & 0.196 & 
  177 & 0.038 & 0.015 & 1.221 & 2.61 & 
  257 & 0.019 & 0.016 & 0.003 & 11.401 & 
  241 & 0.004 & 0.004 & 28.19 & 554.874 \\ 
  253 & 0 & 0.003 & 0.001 & 0.172 &
  218 & 0.017 & 0.013 & 1.105 & 2.171 &  &  &  &  &  & 
  265 & 0.006 & 0.006 & 28.19 & 411.649 \\ 
   &  &  &  &  & 
   235 & 0.017 & 0.003 & 0.673 & 2.053 &  &  &  &  &  &  &  &  &  &  \\ 
   &  &  &  &  & 
   261 & 0.037 & 0.003 & 0.673 & 1.958 &  &  &  &  &  &  &  &  &  &  \\ 
   &  &  &  &  &  &  &  &  &  &  &  &  &  &  &  &  &  &  &  \\ 
\hline
\multicolumn{20}{c}{\textbf{B Gaussian Process}}\\ 
\hline
\multicolumn{5}{c|}{Shifted sinusoidal $d=2$}
&\multicolumn{5}{c|}{Shifted sinusoidal $d=5$}
& \multicolumn{5}{c|}{Rosenbrock $d=2$}
& \multicolumn{5}{c}{Rosenbrock $d=5$}
\\\hline
\multicolumn{1}{c}{$k$}& \multicolumn{1}{c}{LHS} & \multicolumn{1}{c}{RHS} & \multicolumn{1}{c}{$y$} & \multicolumn{1}{c}{$z$} 
& \multicolumn{1}{c}{$k$}& \multicolumn{1}{c}{LHS} & \multicolumn{1}{c}{RHS} & \multicolumn{1}{c}{$y$} & \multicolumn{1}{c}{$z$} 
& \multicolumn{1}{c}{$k$}& \multicolumn{1}{c}{LHS} & \multicolumn{1}{c}{RHS} & \multicolumn{1}{c}{$y$} & \multicolumn{1}{c}{$z$} 
& \multicolumn{1}{c}{$k$}& \multicolumn{1}{c}{LHS} & \multicolumn{1}{c}{RHS} & \multicolumn{1}{c}{$y$} & \multicolumn{1}{c}{$z$} \\ 
\hline
  10 & 0.445 & 0.41 & 1.63 & 2.769 & 
  10 & {\hl {0.155}} & {\hl {0.169}} & 2.811 & 3.411 & 
  10 & {\hl {0.502 }}& {\hl {0.503}} & 8.298 & 35.649 & 
  10 & {\hl {0.399 }}& {\hl {0.4}} & 427.597 & 863.028 \\ 
  64 & 0.17 & 0.14 & 1.002 & 2.636 & 
  15 & 0.548 & 0.534 & 2.489 & 2.811 &
  15 & {\hl {0.416}} & {\hl {0.426}} & 8.298 & 35.649 & 
  62 & {\hl {0.23}} & {\hl {0.234}} & 427.597 & 1407.291 \\ 
  94 & 0.141 & 0.131 & 0.848 & 2.097 & 
  25 & 0.37 & 0.349 & 2.44 & 2.961 & 
  25 & {\hl {0.176 }}& {\hl {0.182}} & 3.834 & 35.649 &
  118 & {\hl {0.226}} &{\hl { 0.23}} & 427.597 & 1407.291 \\ 
  102 & 0.18 & 0.157 & 0.848 & 1.933 & 
  78 & 0.101 & 0.06 & 2.068 & 3.411 & 
  53 & {\hl {0.047}} & {\hl {0.051}} & 1.22 & 72.152 & 
  173 & 0.216 & 0.216 & 427.597 & 1533.6 \\ 
  109 & 0.217 & 0.183 & 0.848 & 1.9 & 
  133 & 0.124 & 0.056 & 2.068 & 3.442 & 
  59 & 0.098 & 0.098 & 1.22 & 33.315 & 
  225 & 0.199 & 0.199 & 427.597 & 1695.92 \\ 
  119 & 0.27 & 0.239 & 0.848 & 1.615 & 
  188 & 0.094 & 0.049 & 2.068 & 3.463 & 
  68 & 0.066 & 0.065 & 0.416 & 15.14 & 
  254 &{\hl {0.177 }}& {\hl {0.186}} & 427.597 & 1761.241 \\ 
  130 & 0.355 & 0.286 & 0.572 & 1.331 & 
  239 & 0.091 & 0.041 & 2.068 & 3.506 & 
  123 & 0.068 & 0.064 & 0.416 & 29.507 &  &  &  &  &  \\ 
  191 & 0.167 & 0.113 & 0.262 & 1.456 & 
  273 & 0.153 & 0.051 & 2.068 & 3.463 & 
  184 & 0.079 & 0.075 & 0.416 & 21.804 &  &  &  &  &  \\ 
  259 & 0.125 & 0.058 & 0.105 & 1.361 &  &  &  &  &  & 
  250 & 0.081 & 0.077 & 0.416 & 20 &  &  &  &  &  \\ 
   &  &  &  &  &  &  &  &  &  & 
   266 & 0.053 & 0.052 & 0.243 & 11.429 &  &  &  &  &  \\ 
   &  &  &  &  &  &  &  &  &  &  &  &  &  &  &  &  &  &  &  \\ 
\hline
\multicolumn{20}{c}{\textbf{C Quadratic Regression}}\\ 
\hline
\multicolumn{5}{c|}{Shifted sinusoidal $d=2$}
&\multicolumn{5}{c|}{Shifted sinusoidal $d=5$}
& \multicolumn{5}{c|}{Rosenbrock $d=2$}
& \multicolumn{5}{c}{Rosenbrock $d=5$}
\\\hline
\multicolumn{1}{c}{$k$}& \multicolumn{1}{c}{LHS} & \multicolumn{1}{c}{RHS} & \multicolumn{1}{c}{$y$} & \multicolumn{1}{c}{$z$} 
& \multicolumn{1}{c}{$k$}& \multicolumn{1}{c}{LHS} & \multicolumn{1}{c}{RHS} & \multicolumn{1}{c}{$y$} & \multicolumn{1}{c}{$z$} 
& \multicolumn{1}{c}{$k$}& \multicolumn{1}{c}{LHS} & \multicolumn{1}{c}{RHS} & \multicolumn{1}{c}{$y$} & \multicolumn{1}{c}{$z$} 
& \multicolumn{1}{c}{$k$}& \multicolumn{1}{c}{LHS} & \multicolumn{1}{c}{RHS} & \multicolumn{1}{c}{$y$} & \multicolumn{1}{c}{$z$} \\ 
\hline
  10 & 0.039 & 0.038 & 0.493 & 2.228 & 
  10 & {\hl {0.155}} & {\hl {0.169}} & 2.811 & 3.411 & 
  10 & 0.238 & 0.235 & 7.63 & 123.944 & 
  10 & {\hl {0.093 }}&{\hl { 0.096}} & 319.605 & 2467.76 \\ 
  60 & 0.235 & 0.235 & 0.493 & 0.913 & 
  15 & 0.556 & 0.551 & 2.67 & 2.961 & 
  39 & {\hl {0.104}} & {\hl {0.113}} & 4.086 & 100 & 
  18 & 0.067 & 0.067 & 143.092 & 854.155 \\ 
  110 & 0.301 & 0.3 & 0.493 & 0.904 & 
  32 & 0.204 & 0.176 & 2.594 & 3.347 & 
  40 & {\hl {0.065}} &{\hl { 0.071}} & 2.563 & 72.526 & 
  25 & 0.061 & 0.061 & 143.092 & 943.599 \\ 
  162 & 0.399 & 0.372 & 0.493 & 0.913 &
  38 & 0.275 & 0.225 & 2.594 & 3.287 & 
  67 & {\hl {0.114}} & {\hl {0.116}} & 1.653 & 29.617 & 
  35 & 0.069 & 0.069 & 143.092 & 854.155 \\ 
  171 & 0.094 & 0.078 & 0.077 & 0.912 & 
  47 & 0.254 & 0.22 & 2.594 & 3.284 & 
  75 & {\hl {0.062 }}&{\hl { 0.065}} & 0.83 & 28.374 & 
  45 & 0.063 & 0.063 & 143.092 & 937.713 \\ 
  226 & 0.104 & 0.072 & 0.077 & 0.898 &
  57 & 0.24 & 0.213 & 2.594 & 3.287 & 
  128 & {\hl {0.065}} &{\hl { 0.067}} & 0.83 & 27.782 & 
  55 & 0.027 & 0.027 & 88.036 & 854.155 \\ 
  251 & 0.075 & 0.068 & 0.077 & 0.895 & 
  67 & 0.157 & 0.077 & 1.92 & 3.284 & 
  173 & {\hl {0.036}} & {\hl {0.044}} & 0.643 & 16.368 & 
  106 & {\hl {0.025}} & {\hl {0.027 }}& 88.036 & 830.225 \\ 
   &  &  &  &  & 
   125 & 0.181 & 0.139 & 1.92 & 2.982 & 
   182 & {\hl {0.041}} &{\hl { 0.047}} & 0.48 & 10.173 &
   158 & 0.03 & 0.029 & 88.036 & 830.225 \\ 
   &  &  &  &  & 
   135 & 0.073 & 0.036 & 1.452 & 2.981 & 
   244 & {\hl {0.004}} &{\hl { 0.005}} & 0.028 & 13.268 & 
   207 & 0.013 & 0.012 & 49.609 & 618.208 \\ 
   &  &  &  &  & 
   152 & 0.079 & 0.036 & 1.452 & 2.961 & 
   266 & {\hl {0.005}} & {\hl {0.006}} & 0.028 & 9.159 & 
   216 & 0.01 & 0.01 & 41.59 & 567.647 \\ 
   &  &  &  &  & 
   170 & 0.087 & 0.048 & 1.452 & 2.792 &  &  &  &  &  & 
   262 & 0.013 & 0.013 & 41.59 & 449.238 \\ 
   &  &  &  &  & 
   184 & 0.083 & 0.052 & 1.423 & 2.67 &  &  &  &  &  &  &  &  &  &  \\ 
   &  &  &  &  & 
   221 & 0.032 & 0.015 & 1.192 & 2.496 &  &  &  &  &  &  &  &  &  &  \\ 
   &  &  &  &  & 
   261 & 0.034 & 0.022 & 1.192 & 2.273 &  &  &  &  &  &  &  &  &  &  \\ 
   \hline
\end{tabular}
\caption{Numerical values for testing Assumption 1.3  with adaptive subregion probabilities (d) confidence bound.
Highlighted cells violate Assumption 1.3.
}
\label{table:table_assumption1.3_d}
\end{table}
\end{landscape}

\section{Testing Assumption~2 Numerically}
\label{sec:Test}

\begin{figure}[htbp!]
    \centering
    \includegraphics[width=0.8\textwidth]{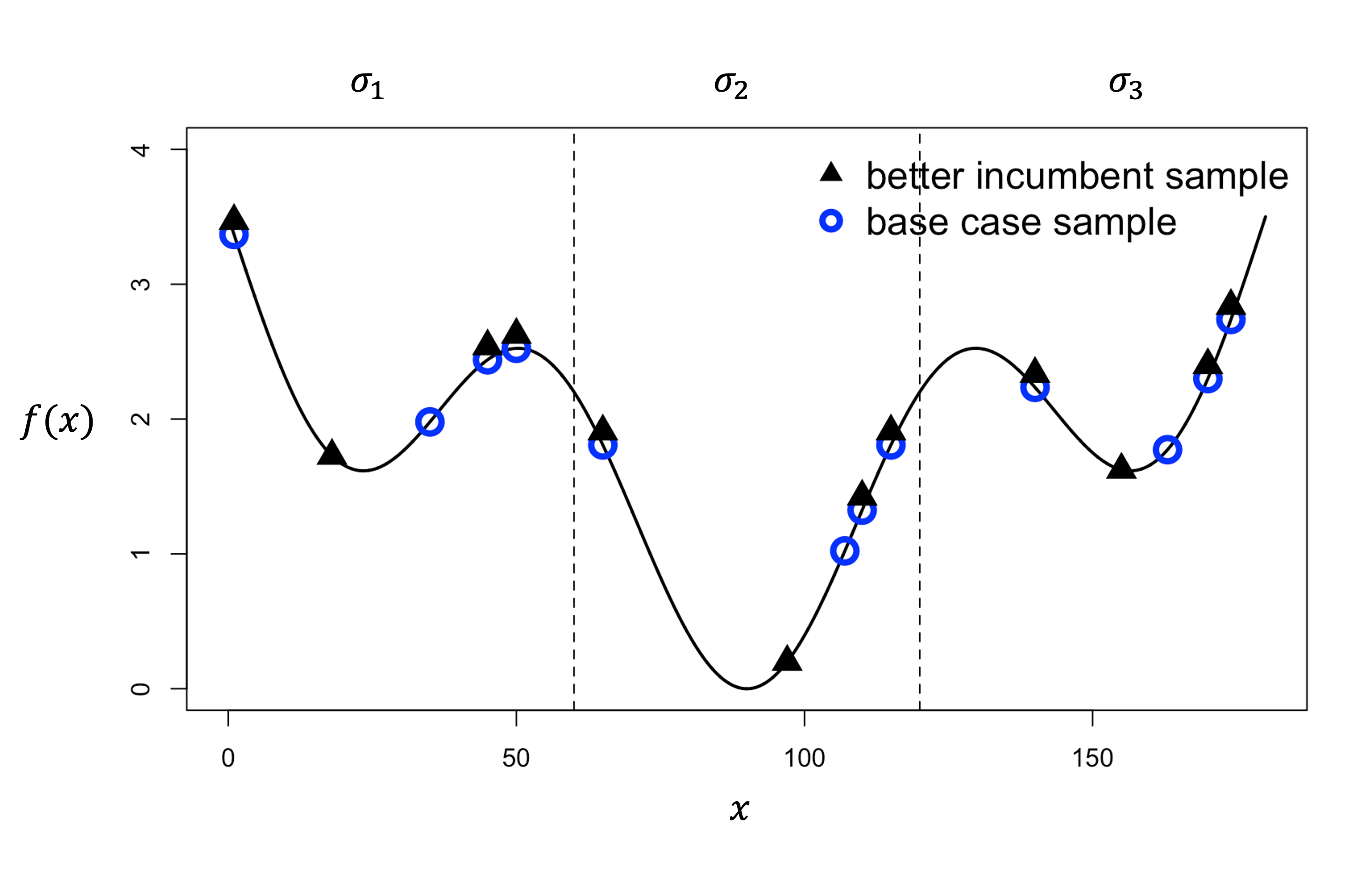}
    \caption{This figure shows two different sets of sample points (blue circle base case sample set and black triangle better incumbent sample set) used to build the Gaussian processes with 3 of the 4 sample points in each subregion coinciding.  The black triangle sample set has better incumbent points than the blue circle sample set.}
    \label{fig:sampleset}
\end{figure}

\begin{figure}[htbp!]
     \centering
     \begin{subfigure}[b]{0.8\textwidth}
         \centering
         \includegraphics[width=\textwidth]{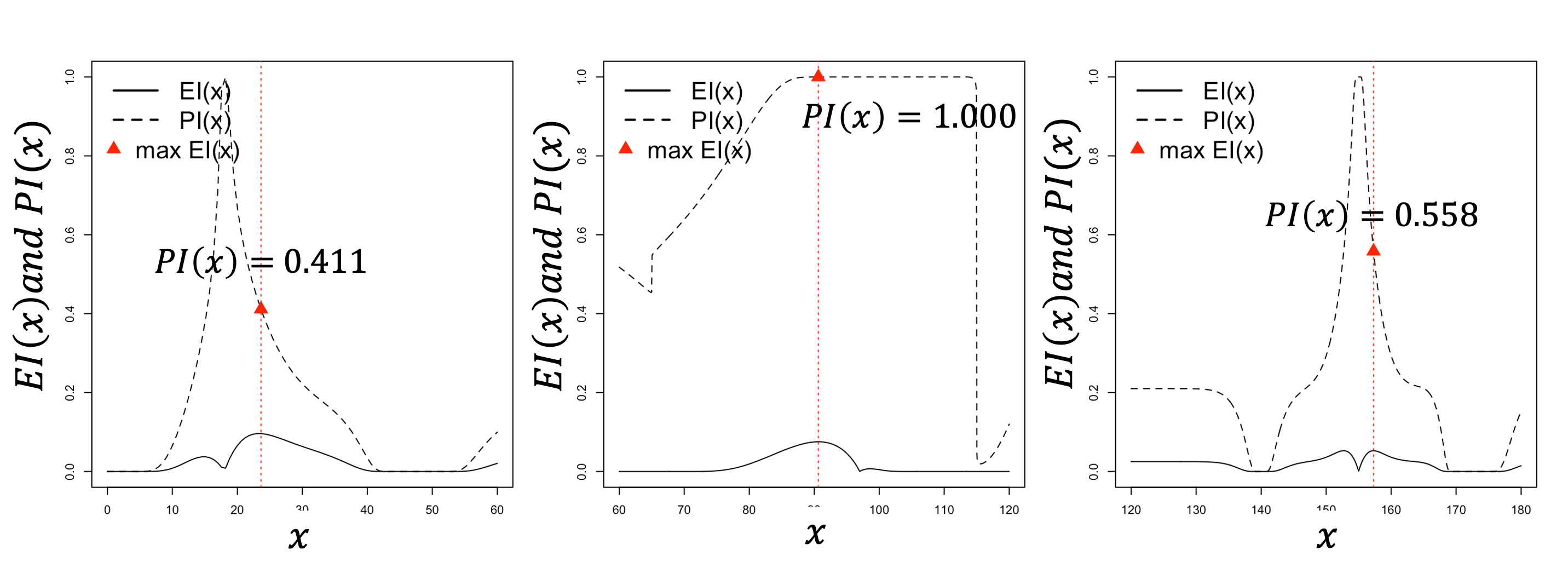}
         \caption{Gaussian process using  the black triangle better incumbent sample set.
        % sample points ( sample set with better incumbent) in Figure ~\ref{fig:sampleset}
         }
         \label{fig:PIbetter}
     \end{subfigure}
     \begin{subfigure}[b]{0.8\textwidth}
         \centering
         \includegraphics[width=\textwidth]{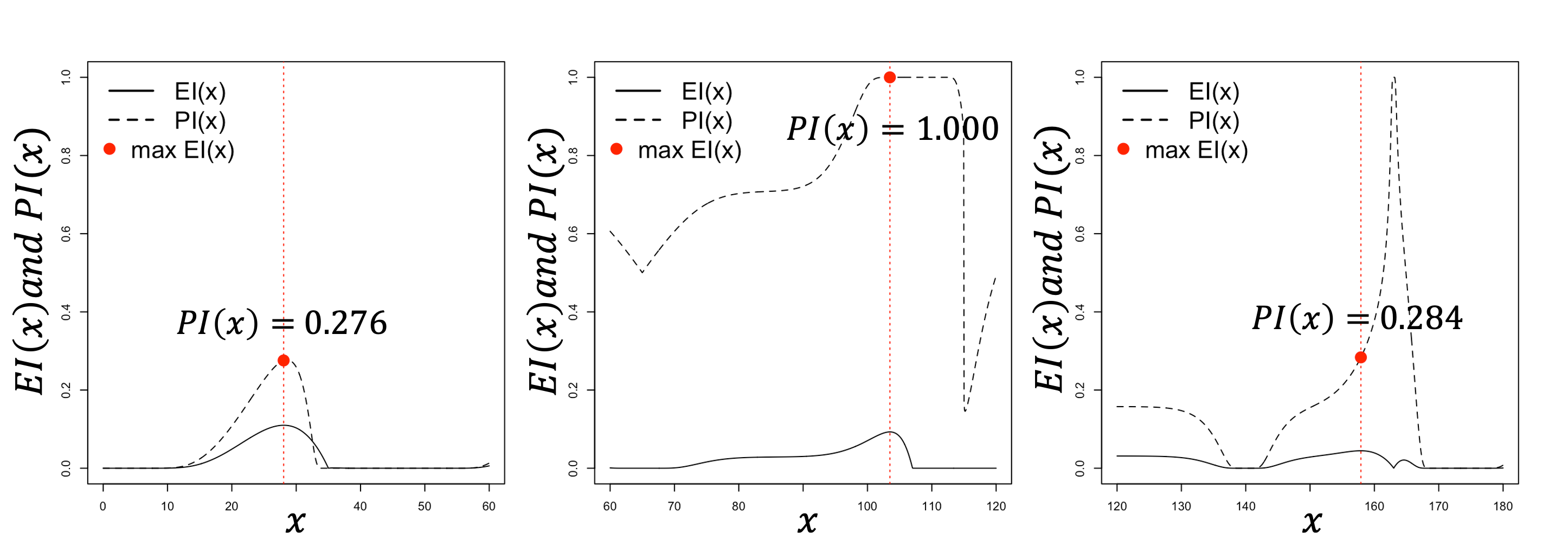}
         \caption{Gaussian process using the blue circle base case sample set.}
         \label{fig:PIbase}
         \end{subfigure}
     \begin{subfigure}[b]{0.8\textwidth}
         \centering
         \includegraphics[width=\textwidth]{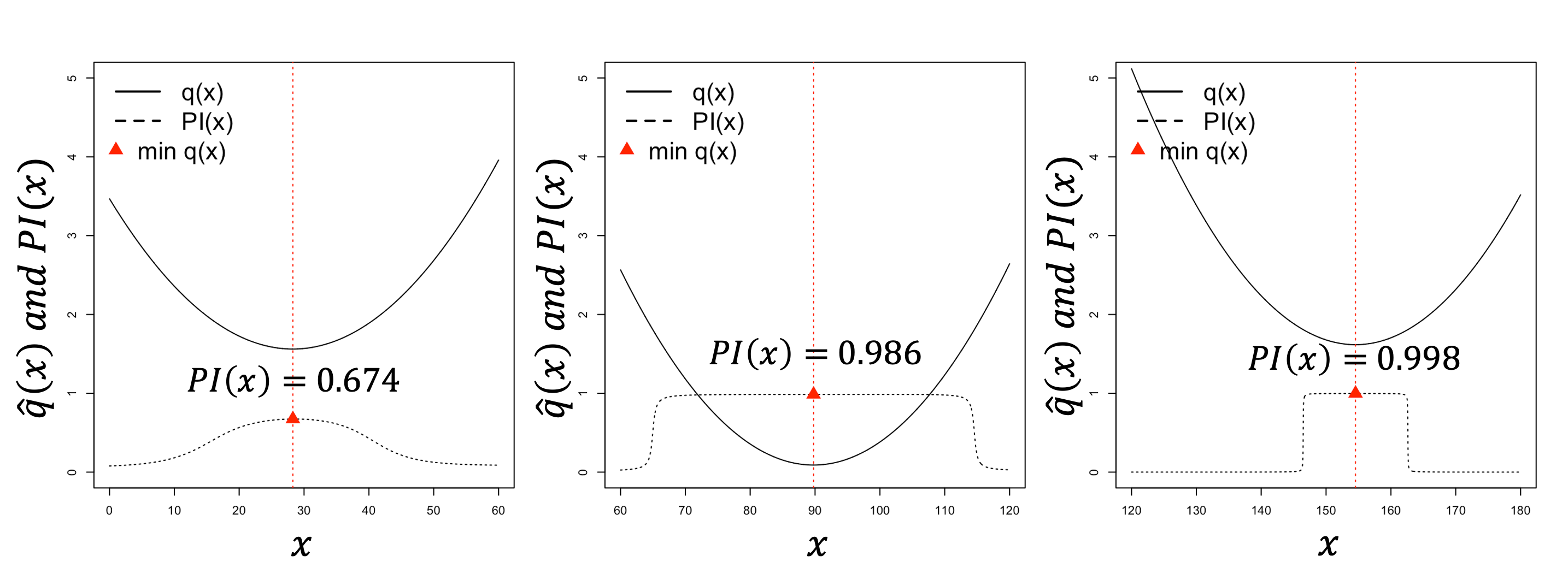}
         \caption{ Quadratic Regression using the black triangle better incumbent sample set. 
        % sample points ( sample set with better incumbent) in Figure ~\ref{fig:sampleset}
         }
         \label{fig:PIbetter_Q}
     \end{subfigure}
    \begin{subfigure}[b]{0.8\textwidth}
         \centering
         \includegraphics[width=\textwidth]{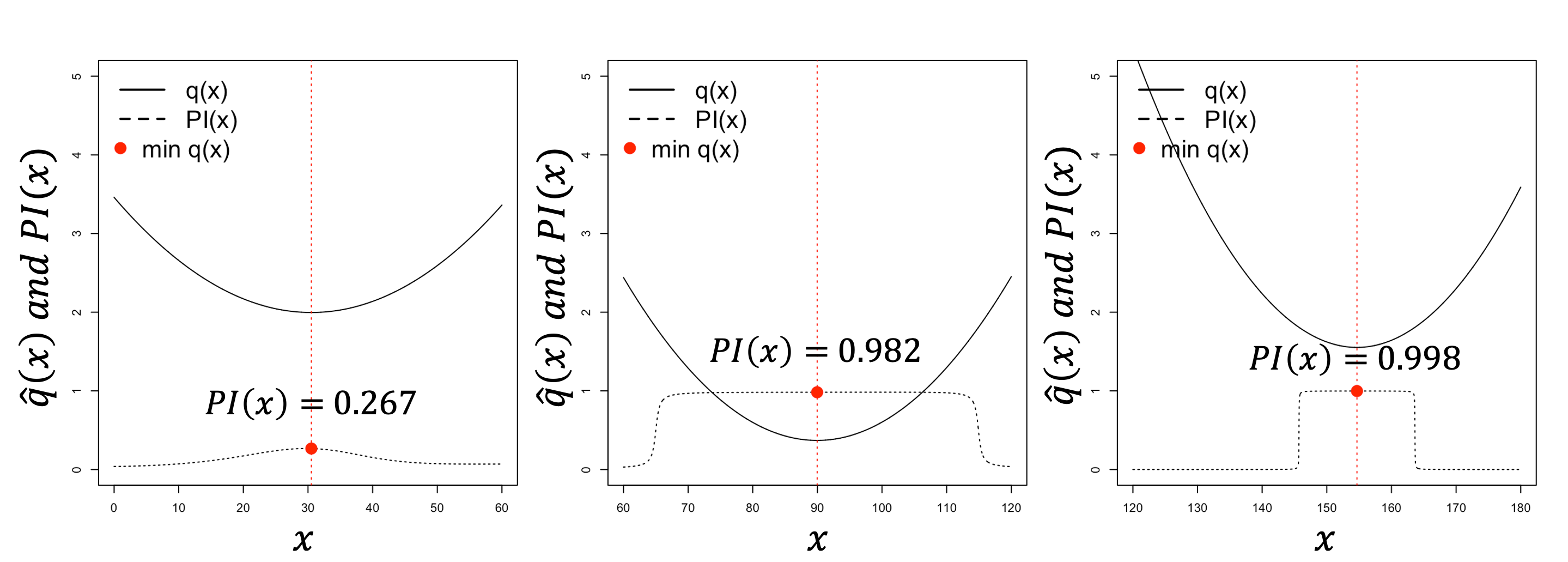}
         \caption{Quadratic regression using the blue circle base case sample set.}
         \label{fig:PIbase_Q}
         \end{subfigure}
        \caption{ Details of Testing Assumption 2.}
        %present the probability of the improvement function ($PI_i(x)$) and the expected improvement function ($EI_i(x)$) of the subregion $\sigma_1$, $\sigma_2$ and $\sigma_3$, from left to right. Observe that, comparing in the same subregion, the $PI_i(x)$  at the point that maximizes $EI_i(x)$ of the sample set with the better incumbent sample set (top row) is higher than the one with the base case sample set (bottom row).Figures (c) and (d) present the probability of the improvement function ($PI_i(x)$) and the quadratic regression  ($\hat{q}_i(x)$) of the subregion $\sigma_1$, $\sigma_2$ and $\sigma_3$, from left to right.Comparing in the same subregion, the $PI_i(x)$  at the point that minimizes $\hat{q}_i(x)$ of the sample set with the better incumbent sample set (top row) is higher than the one with the base case sample set (bottom row) in subregion $\sigma_1$ and $\sigma_2$.}
        \label{fig:PI}
\end{figure}

We generate two sets of sample points,  where each subregion has four sample points. Three of the sample points coincide, and the base case sample set has a higher best objective function value than the better incumbent sample set, as
illustrated in Figure~\ref{fig:sampleset}.
%,  blue solid circles as the base case sample set, and three of the black triangle points coincide with three of the blue circles, but the best black triangle point is better (i.e., lower) than the best blue circle point.
For each subregion, we build two Gaussian processes where one is constructed with the base case sample set and the other with the better incumbent sample set.  Figure~\ref{fig:PI} plots the expected improvement function $EI_i(x)$, as in \eqref{eq:EIx}, for each subregion, where Figure \ref{fig:PIbase} is for the base case sample set, and Figure \ref{fig:PIbetter} is for the  better incumbent sample set. The new sample point for each subregion is the point that maximizes the expected improvement function, illustrated by the vertical dotted line. Notice that there is a different point selected between the base case sample set and the better incumbent sample set, due to the different best function values.

%In order to demonstrate that Assumption~2 is satisfied, we need to calculate the probability of improvement at the new sample point, comparing the probability of improvement using the Gaussian process to that of a uniformly distributed point (Assumption 2.1, \eqref{eqn:surgassumption}), and comparing the probability of improvement between the base case sample set and the better incumbent sample set (Assumption 2.2, \eqref{eqn:IncumbentBASSOsurg}). We calculate the probability of improvement  below a value $y$ within a subregion, as in \eqref{eqn:PI}, where we take $y$ as the 6th best function value  of both sample sets. 

\section{Test Functions and Detailed Plots}
\label{sec:testfunctions}
\textbf{Rosenbrock Problem:} $(-2 \leq x_i \leq 2,\  i=1, \ldots, d.)$ The global minimum is at $x_*=(1,\ldots,1)$ and $f(x_*) = 0$.
$$
%\min_{x\in[-2,2]} 
f(\mathbf{x}) = \sum_{i=1}^{d-1}\left((1-x_i)^2 + 100\left(x_{i+1}-x_i^2\right)^2\right)
%\text{~and~}x\in\mathbb{R}^{dim}
$$
%The global minimum is at $(1,\ldots,1)$ and $f(x^*) = 0$.
% {\red{I'm changing optimum to minimum, and $x^*$ to $x_*$}}

\vspace{12pt}
\noindent\textbf{Centered Sinusoidal Problem:} $(0 \leq x_i \leq 180,\  i=1, \ldots, d.)$ The global minimum is at $x_*= (90,\ldots,90)$ and $f(x_*)=0$. \cite{ali2005}
$$
f(\mathbf{x}) = -\left[2.5\prod_{i=1}^{d} \sin \left(\frac{\pi x_i}{180}\right) + \prod_{i=1}^{d} \sin\left(\frac{\pi x_i}{180}\right)\right] +3.5 $$

\vspace{12pt}
\noindent\textbf{Shifted Sinusoidal Problem:} $(0 \leq x_i \leq 180, \ i =1 \ldots, d.)$ The global minimum is at $x_*=(30,\ldots,30)$ and $f(x_*)=0$. \cite{ali2005}
$$f(\mathbf{x}) = -\left[2.5\prod_{i=1}^{d} \sin\left(\frac{\pi(x_i+60)}{180}\right) + \prod_{i=1}^{d} \sin\left(\frac{\pi(x_i+60)}{36}\right)\right] +3.5 $$
% {\red{I added a subscript $i$ in the sin terms.}}

\vspace{12pt}
\noindent\textbf{Ackley Problem:} $(-32 \leq x_i \leq 32, \ i =1 \ldots, d.)$ The global minimum is at $x_*=(0,\ldots,0)$ and $f(x_*)=0$. \cite{Jamil2013}
$$f(\mathbf{x}) = -20\;\exp\left(-0.2 \sqrt{\frac{1}{d}
                \sum_{i=1}^{d} x_i^2}\right) - \exp\left(\frac{1}{d}
                \sum_{i=1}^{d} \cos(2\pi\;x_i)\right) + 20 + \exp(1)$$

 % {\red{Pete, what is $+ a + \exp(1)$ in the above test function?  Is it a typo?  what is $a$?}}               
\vspace{12pt}
\noindent\textbf{Repeated Branin Problem:} 
% {\red{Pete, this is not clear to me.  I'm going to rewrite it and please check me that I am correct.}}
$(-1 \leq x_i \leq 1, \ i =1,\dots,d).$ \cite{oh2018}
$$f(\mathbf{x}) = \frac{1}{\floor{d/2}}\sum^{\floor{d/2}}_{i=1} f_{branin}(x_{2i-1},x_{2i}),$$
where the two-dimensional Branin function  is,\\ 
$$f_{branin}(x_1,x_2) = \left(x_2 - \left(\frac{5}{4\pi^2}\right)x_1^2 + \left(\frac{5}{\pi}\right)x_1 -6\right)^2 + 10\left(1-\frac{1}{8\pi}\right)\textrm{cos}(x_1)+10.$$

% ***********************
% $(-1 \leq x_i \leq 1, \ i =1).$ 
% $$f(\mathbf{x}) = \frac{1}{\floor{d/2}}\sum^{d/2}_{i=1} f_{branin}(x_{2i-1},x_{2i}),$$
% where Branin function formula is,\\ 
% $$f(x) = (x_2 - (\frac{5}{4\pi^2})x_1^2 + (\frac{5}{\pi})x_1 -6)^2 + 10(1-\frac{1}{8\pi})\textrm{cos}(x_1)+10.$$

\vspace{12pt}
\noindent\textbf{Repeated Hartmann Problem:}
$(-1 \leq x_i \leq 1, \ i =6,\ldots,d).$ 
$$f(\mathbf{x}) = \frac{1}{\floor{d/6}}\sum^{d/6}_{i=1} f_{hartmann6}(x_{6i-5},x_{6i-4},x_{6i-3},x_{6i-2},x_{6i-1},x_{6i}),$$
where the 6-dimensional Hartmann function  is,
$$f_{hartmann6}(x_1,\ldots,x_6) = -\sum_{i=1}^{4} c_i \; \textrm{exp}(-\sum_{j=1}^{6}a_{ij}(x_j - p_{ij})^2).$$
The $a_{ij},c_i$ and $p_{ij}$ for $i=1,\ldots,4,\  j=1,
\ldots, 6$ can be found in \cite{laguna2005}.

\begin{figure}
    \centering 
    \setkeys{Gin}{width=\linewidth}
    \settowidth\rotheadsize{PARAMETERS 3} 
    % from makecell
\begin{tblr}{colspec = { Q[h] *{2}{Q[c,m, wd=30mm]}},
             colsep  = 1pt,
             rowsep  = 0.1pt,
             cell{2-Z}{1} = {cmd=\rotcell, font=\footnotesize\bfseries},
             row{1} = {font=\scriptsize\bfseries},
             measure=vbox
            }
    \toprule
  & dim = 20 & dim = 50 \\
    \midrule
\centering
Ackley & 
{\includegraphics{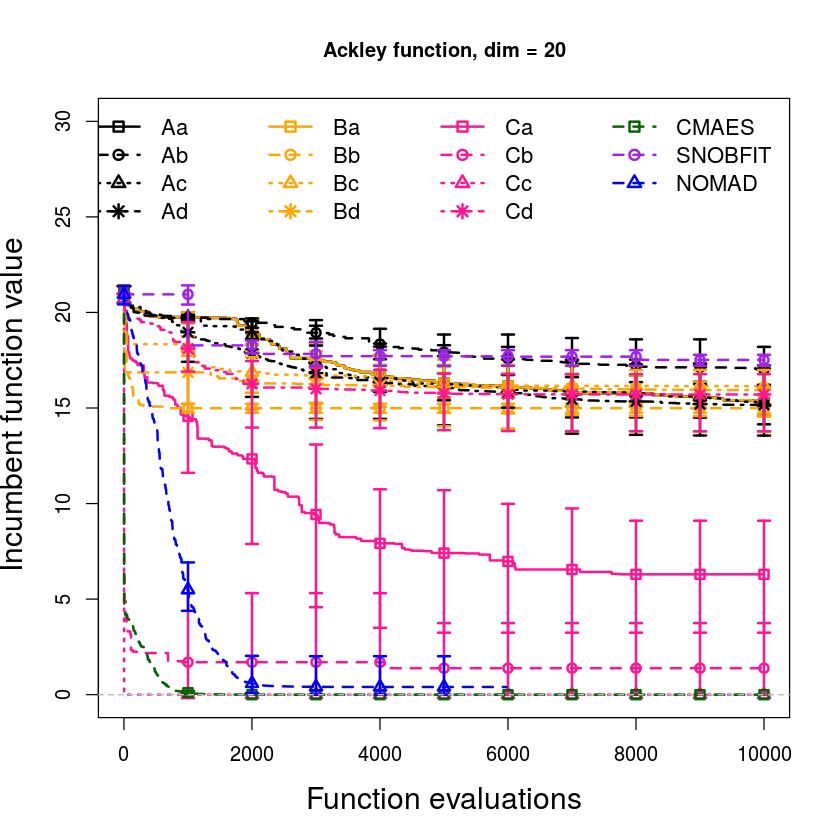}} 
        & {\includegraphics{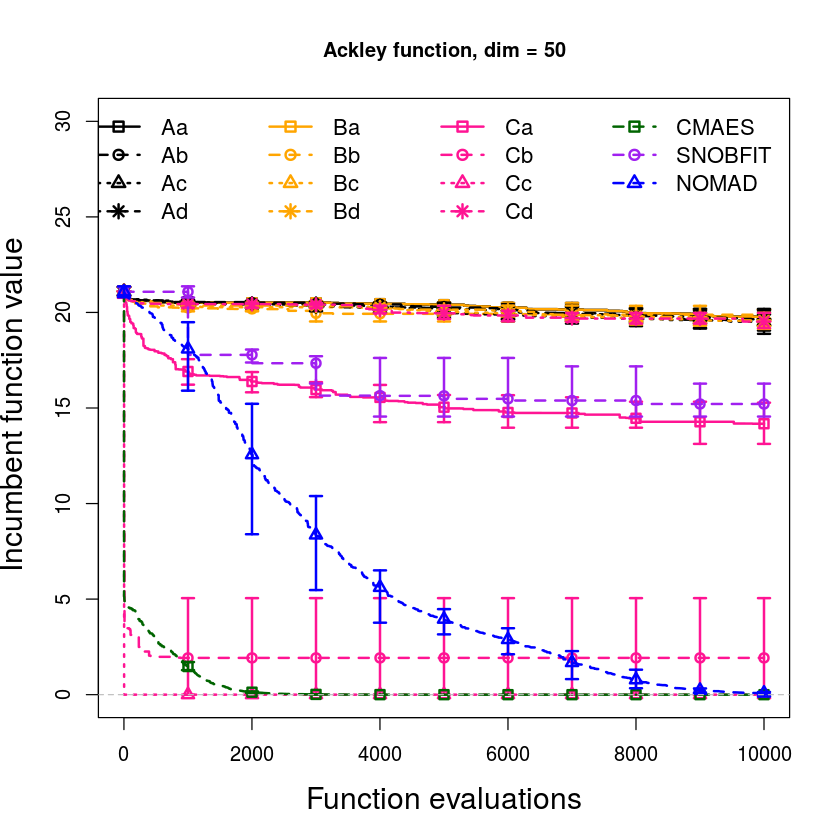}}
\\
\centering
Branin & {\includegraphics{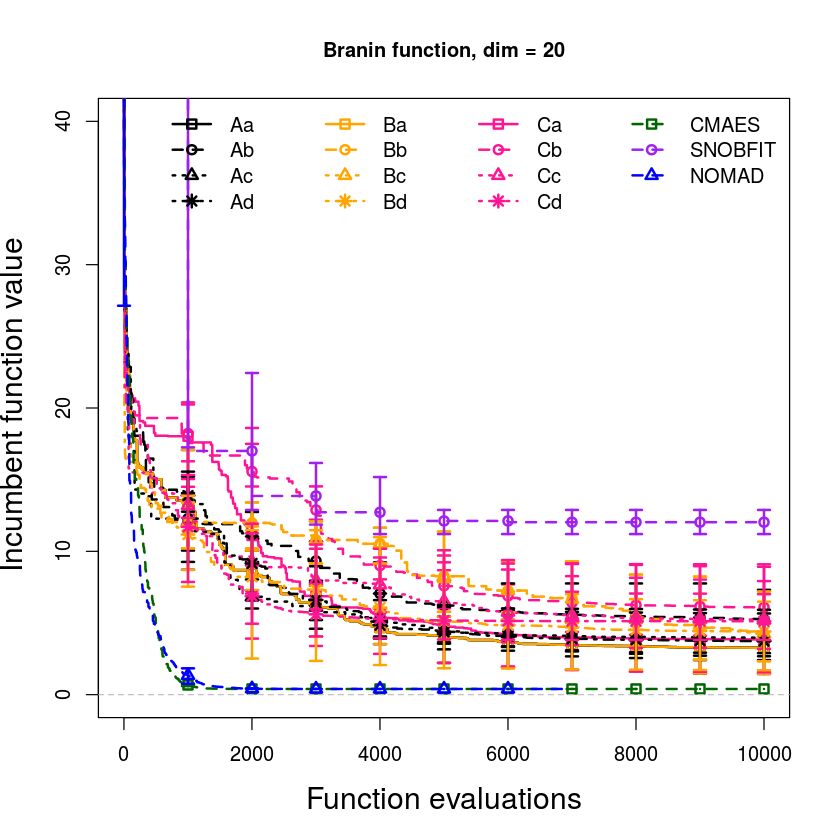}} 
        & {\includegraphics{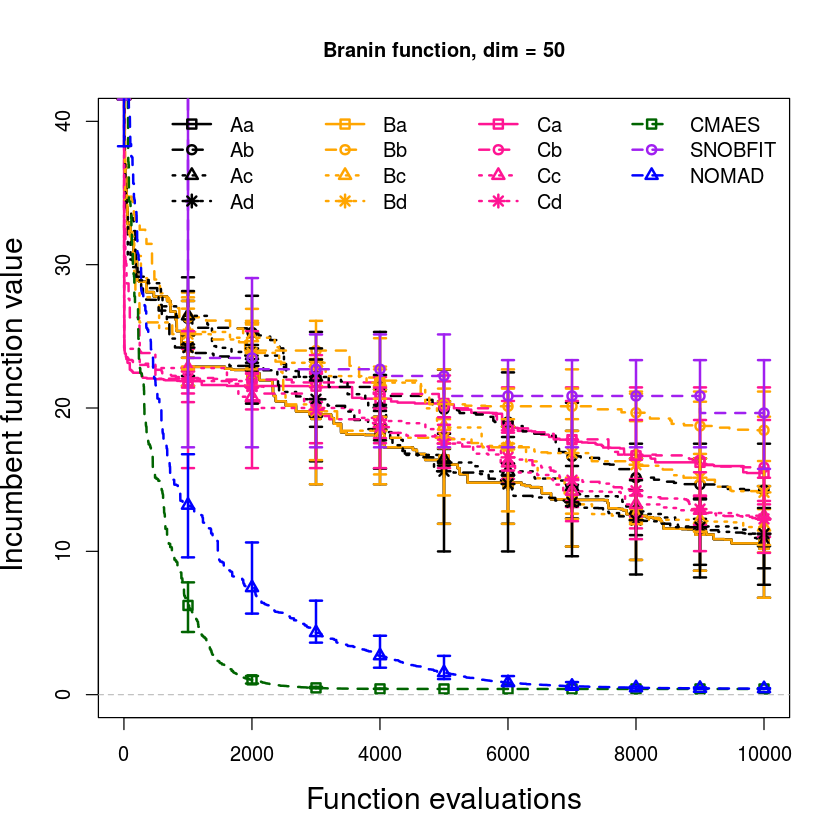}}
\\
\centering
Hartmann & {\includegraphics{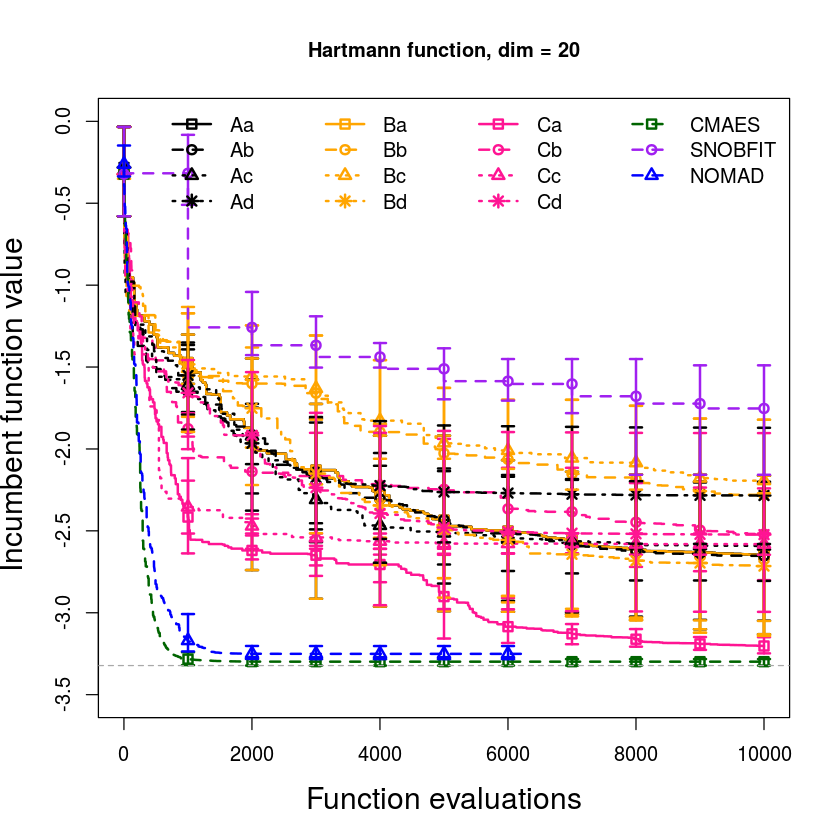}} 
        & {\includegraphics{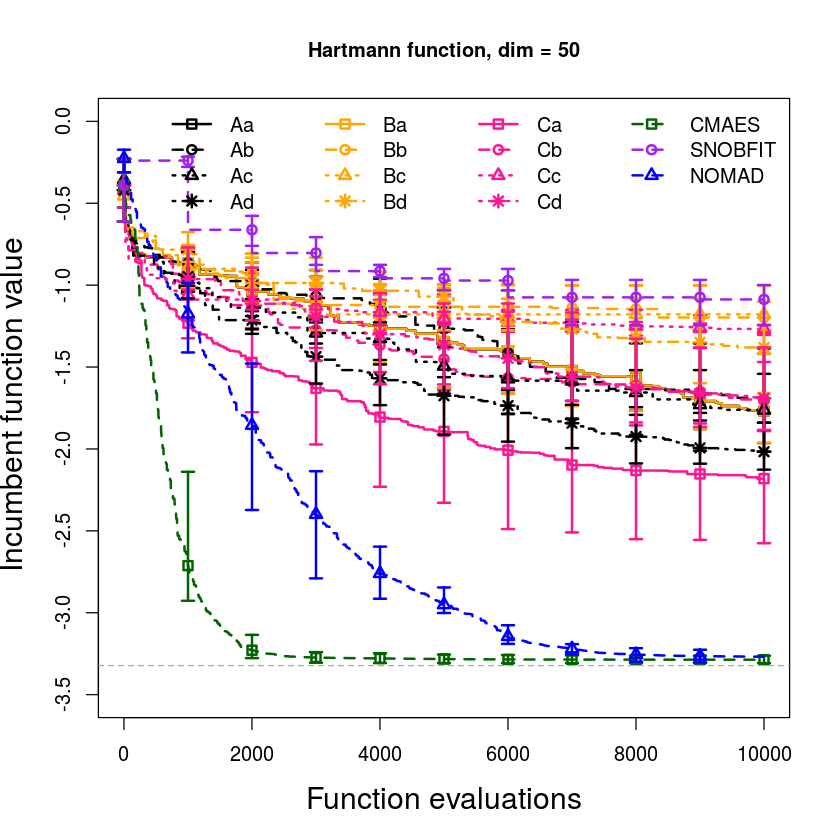}}
\\
\centering
Rosenbrock & {\includegraphics{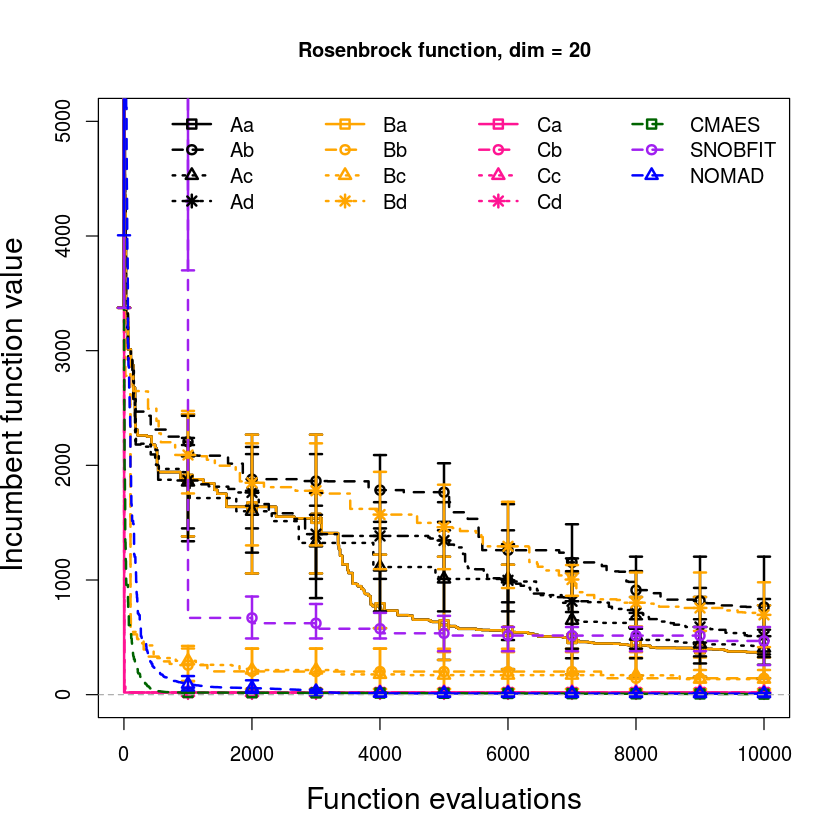}} 
        & {\includegraphics{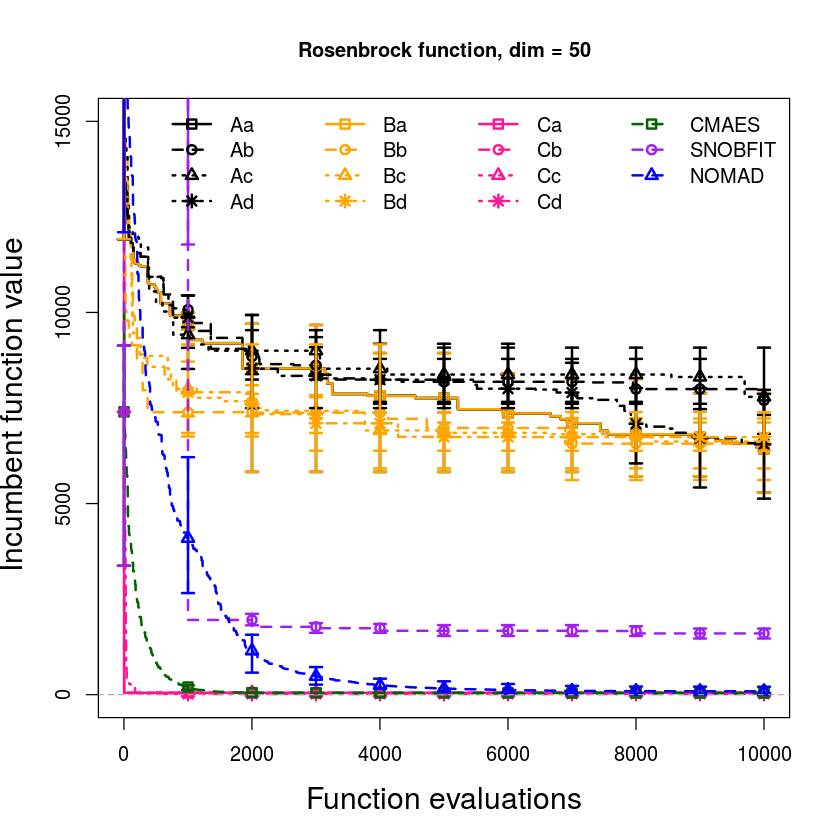}}
\\
\centering
Centered Sinusoidal &{\includegraphics{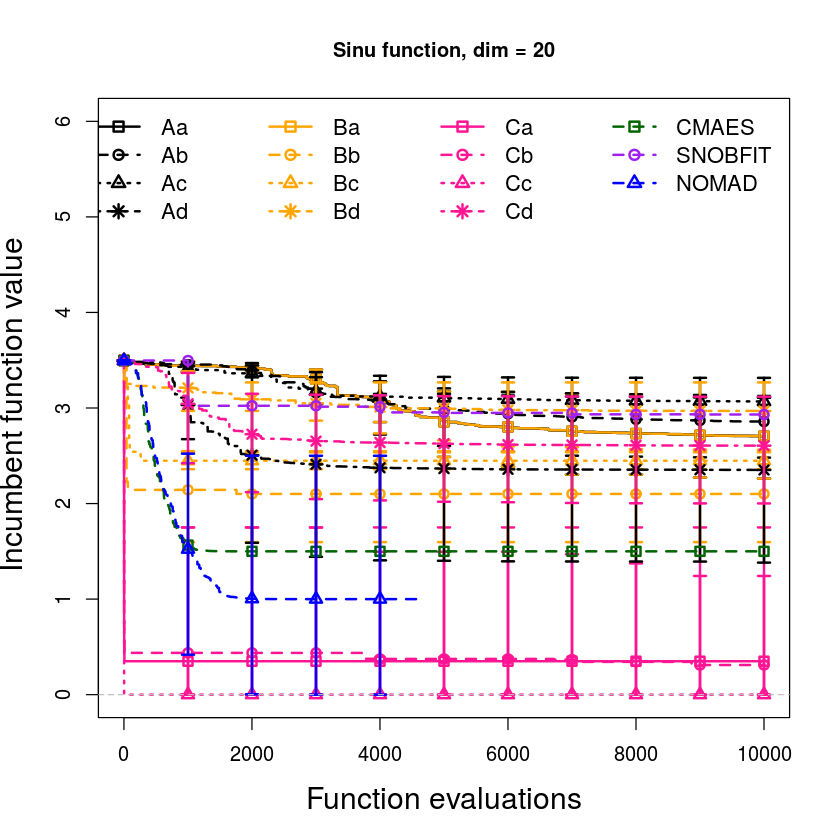}} 
        & {\includegraphics{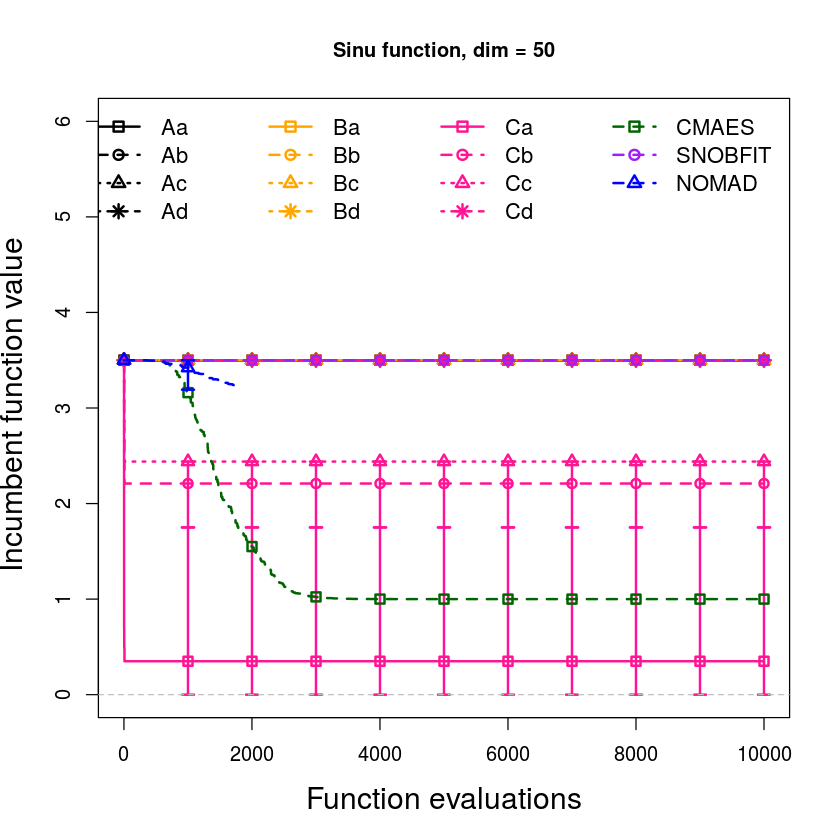}}
\\
\centering
Shifted Sinusoidal &{\includegraphics{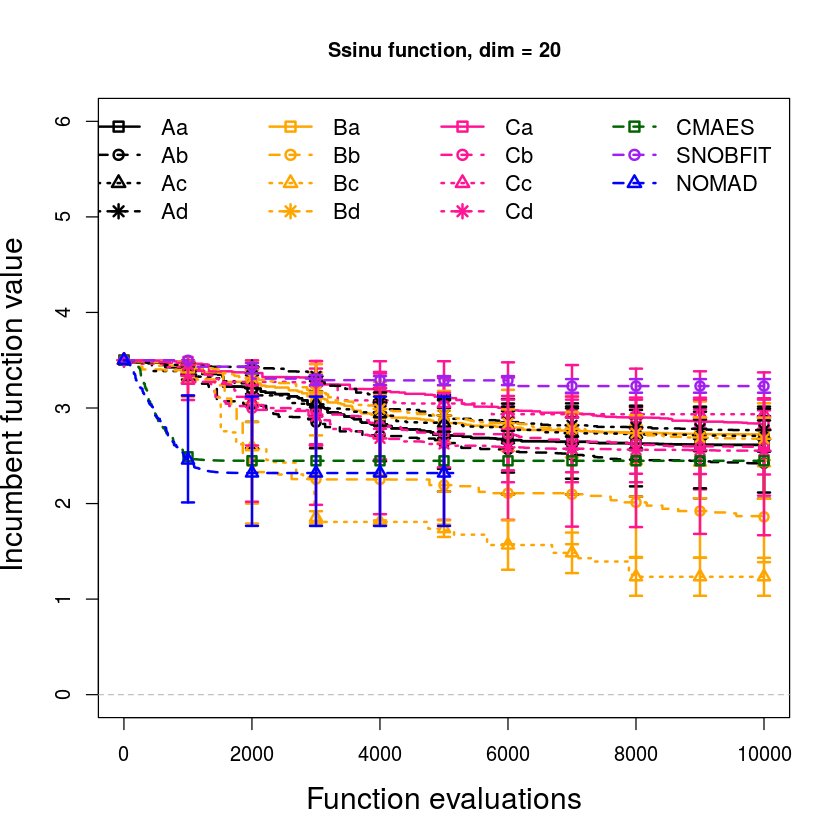}} 
        & {\includegraphics{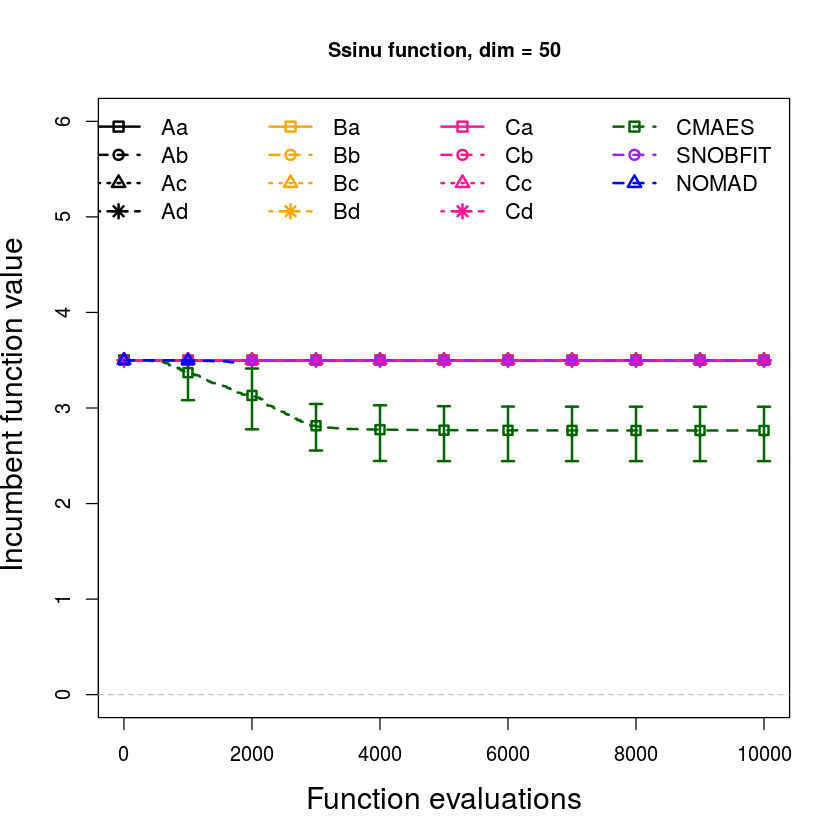}}
\\
    \bottomrule
    \end{tblr}
\caption{
%The incumbent function value vs. number of function evaluations.
 Numerical results at dimensions 20 and 50.}
\label{fig:table_of_figures_exper_2_1}
\end{figure}

\begin{figure}
    \centering 
    \setkeys{Gin}{width=\linewidth}
    \settowidth\rotheadsize{PARAMETERS 3} 
    % from makecell
\begin{tblr}{colspec = { Q[h] *{3}{Q[c,m, wd=30mm]}},
             colsep  = 1pt,
             rowsep = 0.1pt,
             cell{2-Z}{1} = {cmd=\rotcell, font=\footnotesize\bfseries},
             row{1} = {font=\scriptsize\bfseries},
             measure=vbox
            }
    \toprule
  & dim = 100 & dim = 500 & dim = 1000 \\
    \midrule
\centering
Ackley & 
{\includegraphics{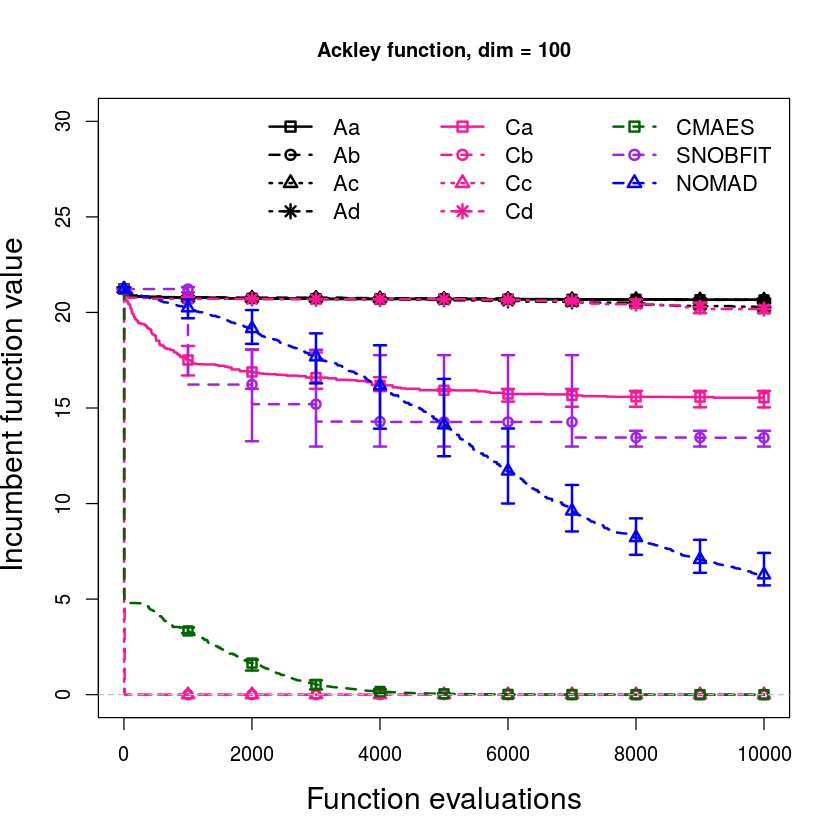}} 
        & {\includegraphics{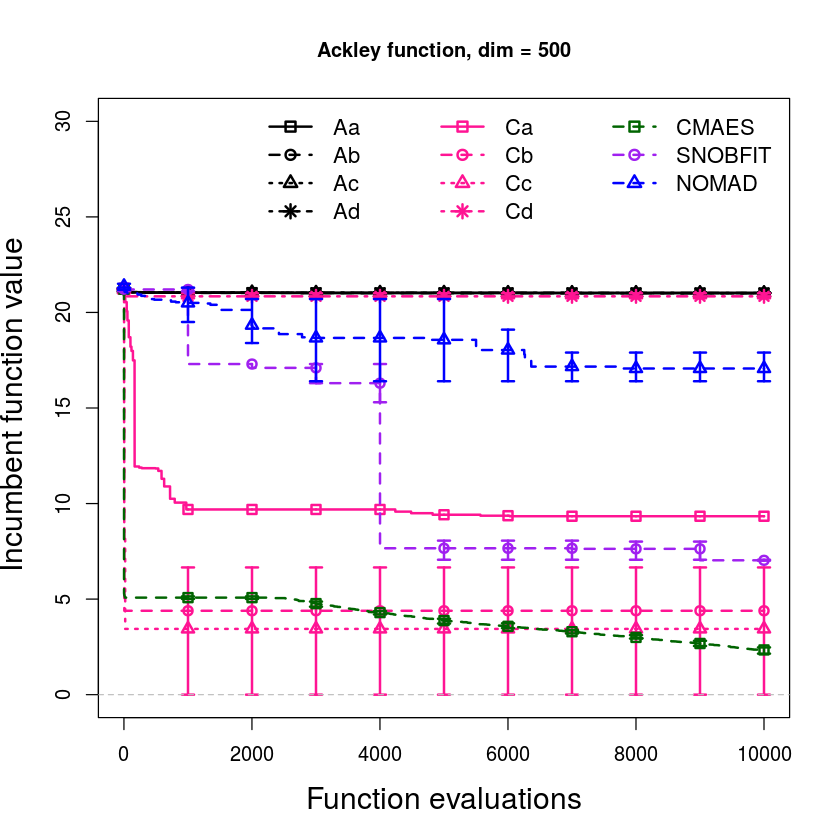}}
                & {\includegraphics{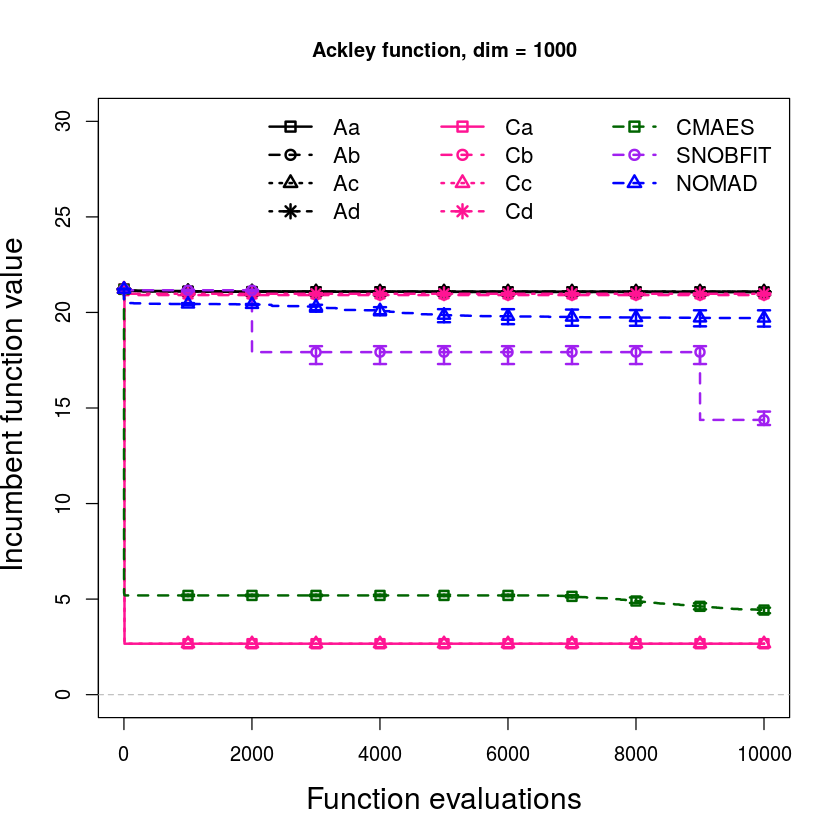}}
\\
\centering
Branin & {\includegraphics{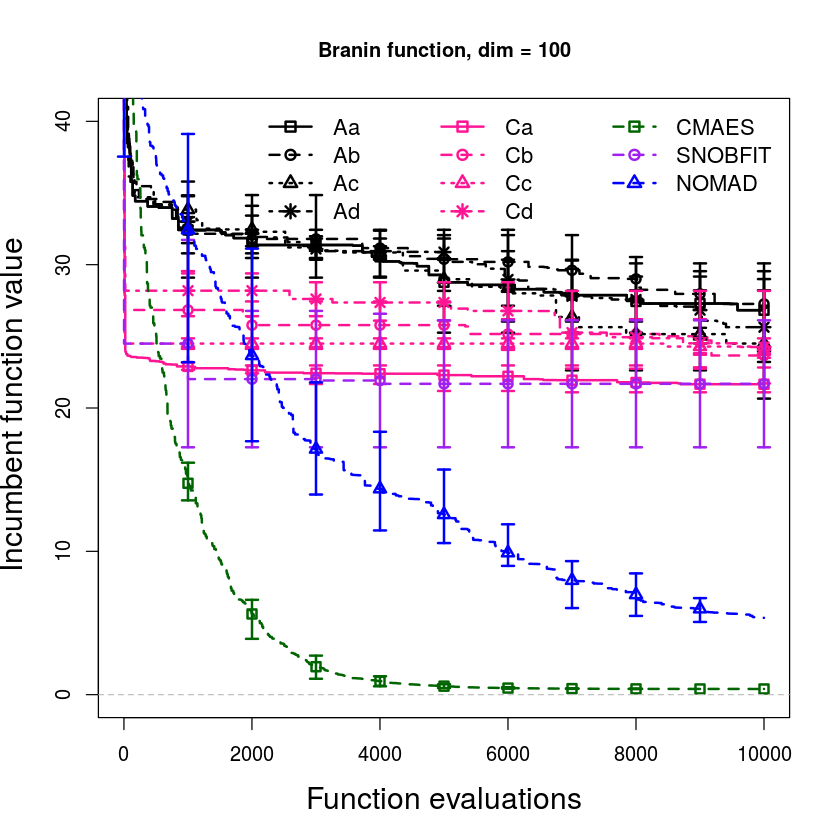}} 
        & {\includegraphics{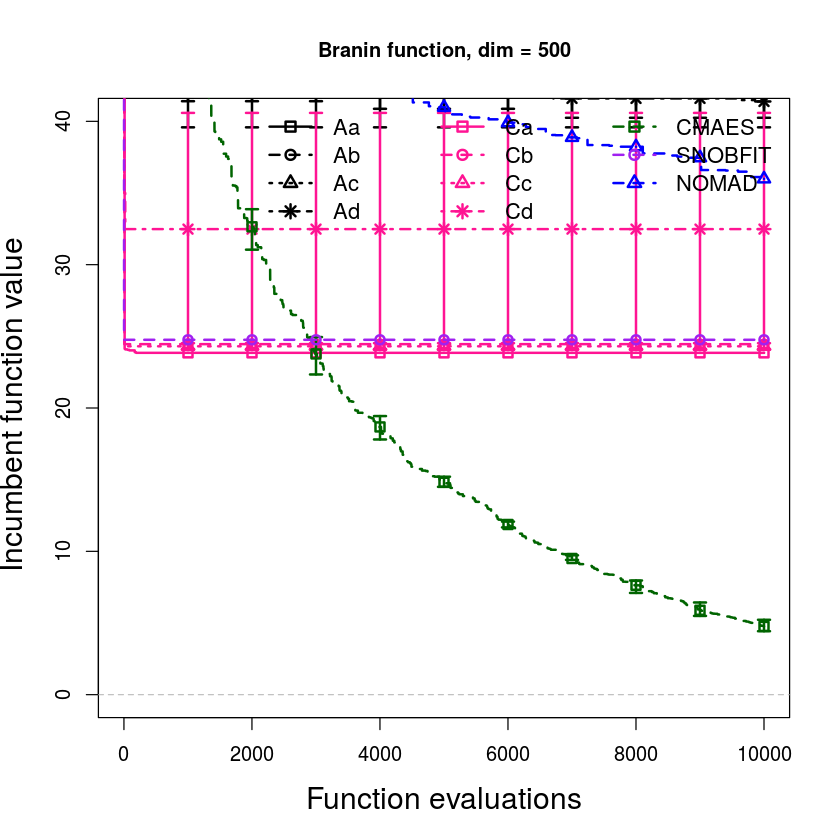}}
                & {\includegraphics{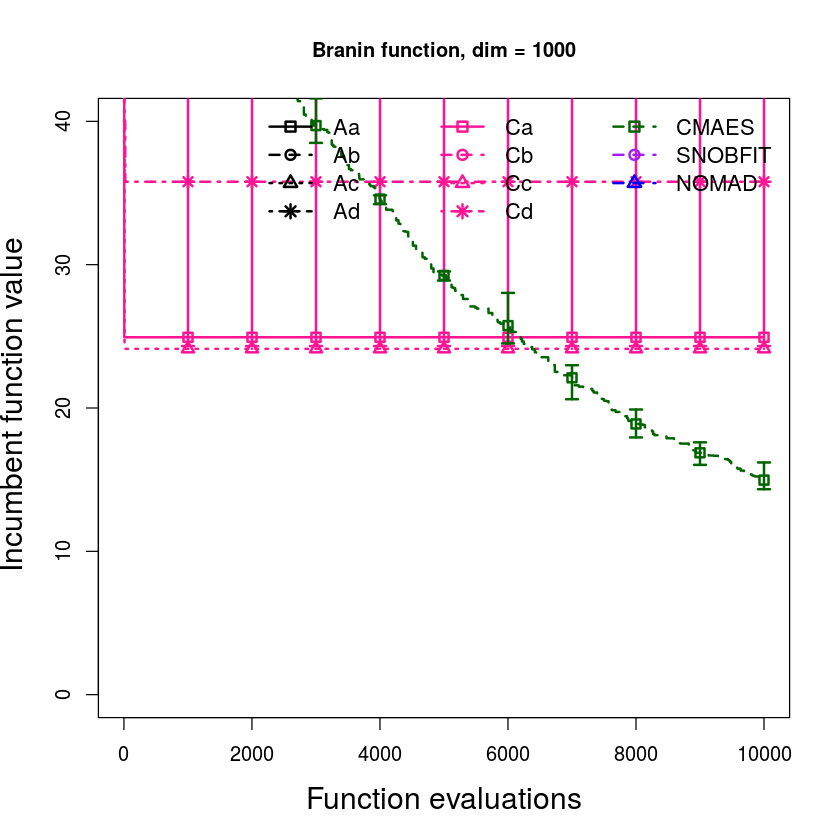}}
\\
\centering
Hartmann & {\includegraphics{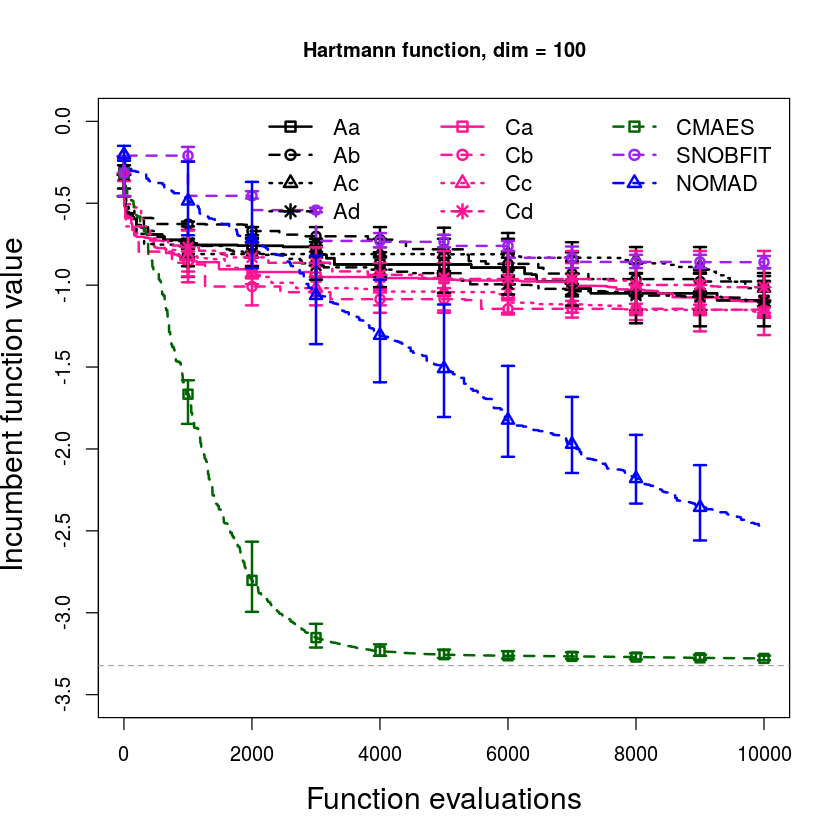}} 
        & {\includegraphics{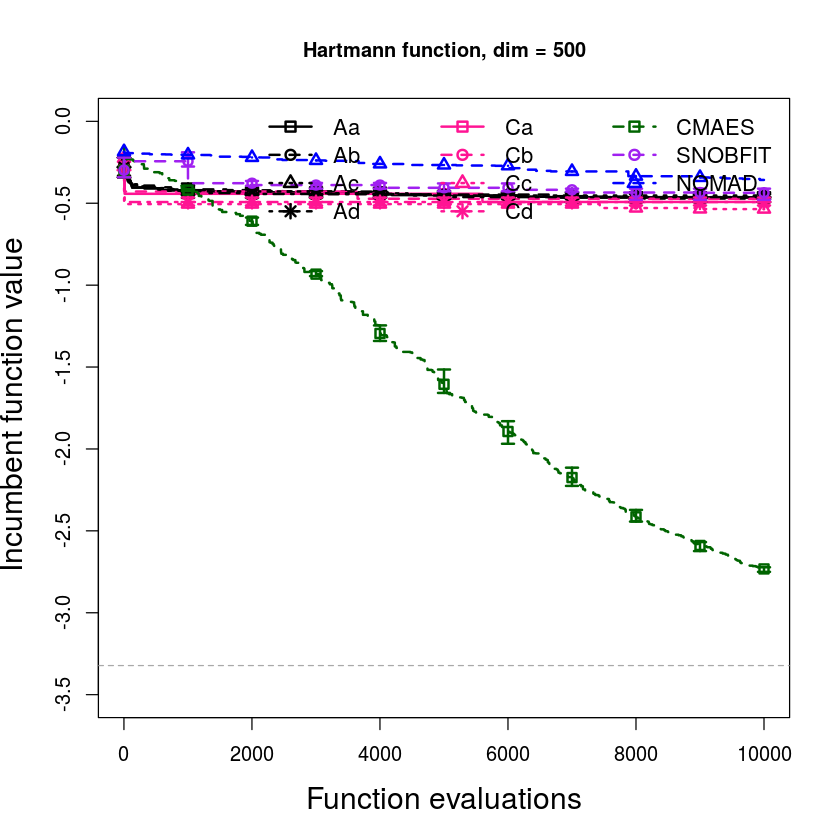}}
                & {\includegraphics{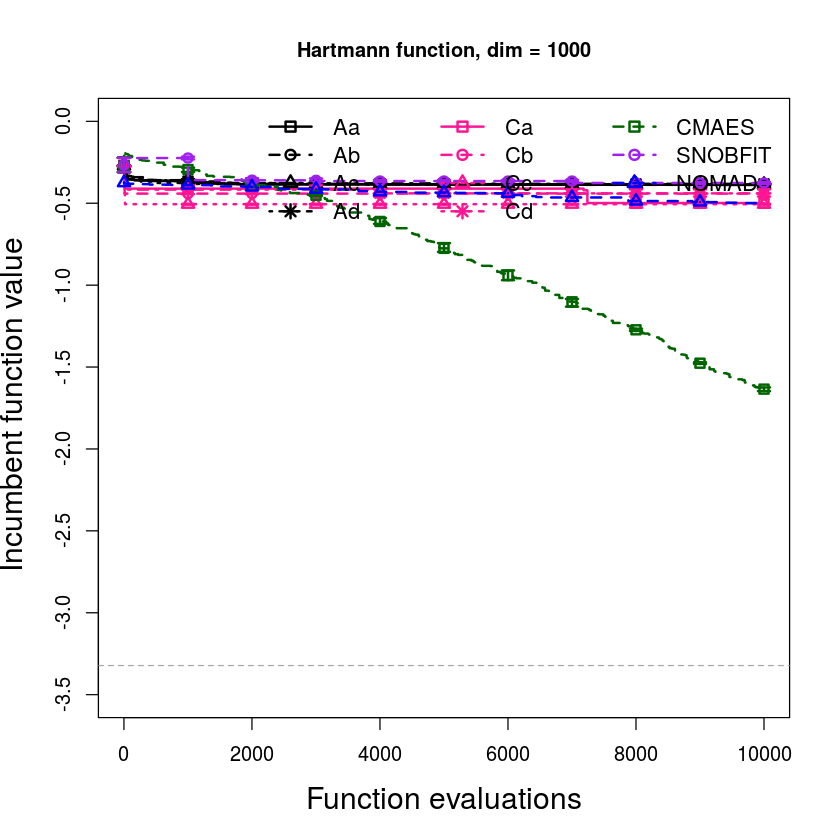}}
\\
\centering
Rosenbrock & {\includegraphics{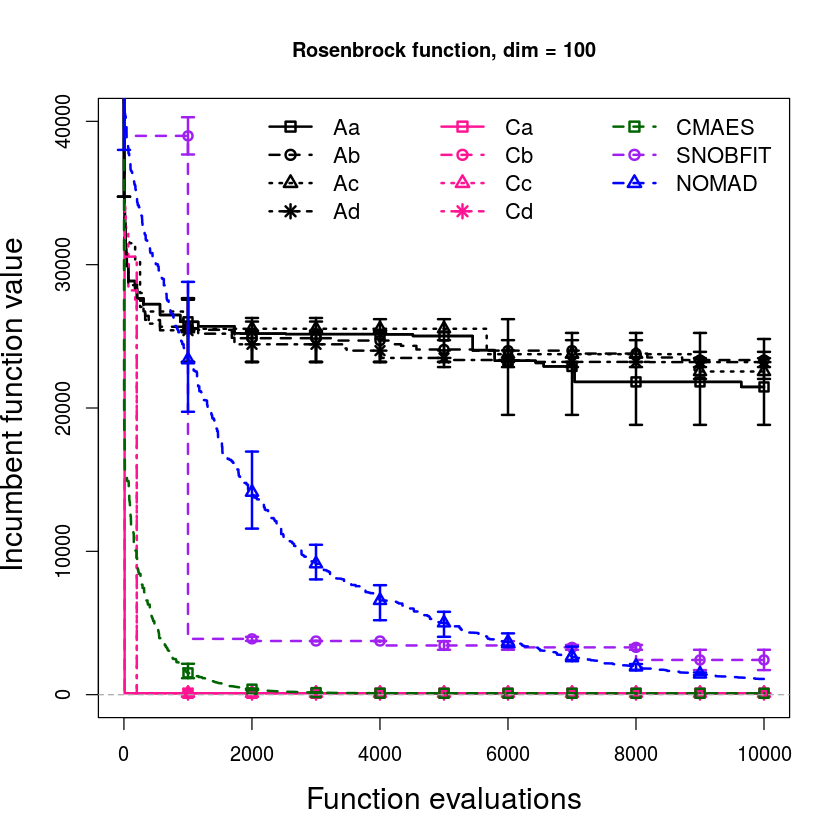}} 
        & {\includegraphics{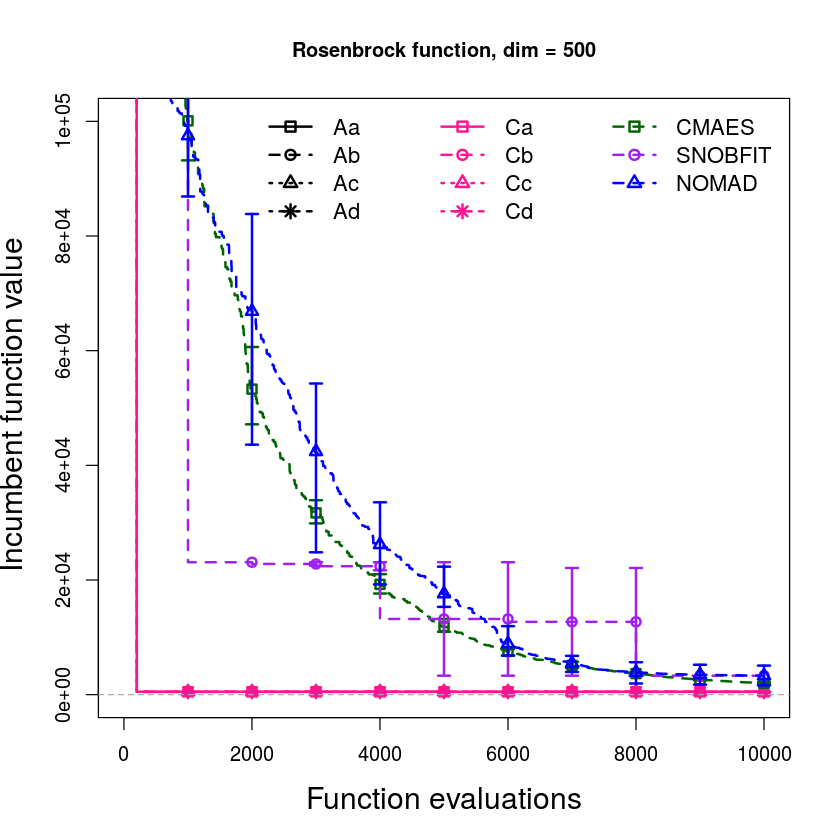}}
             & {\includegraphics{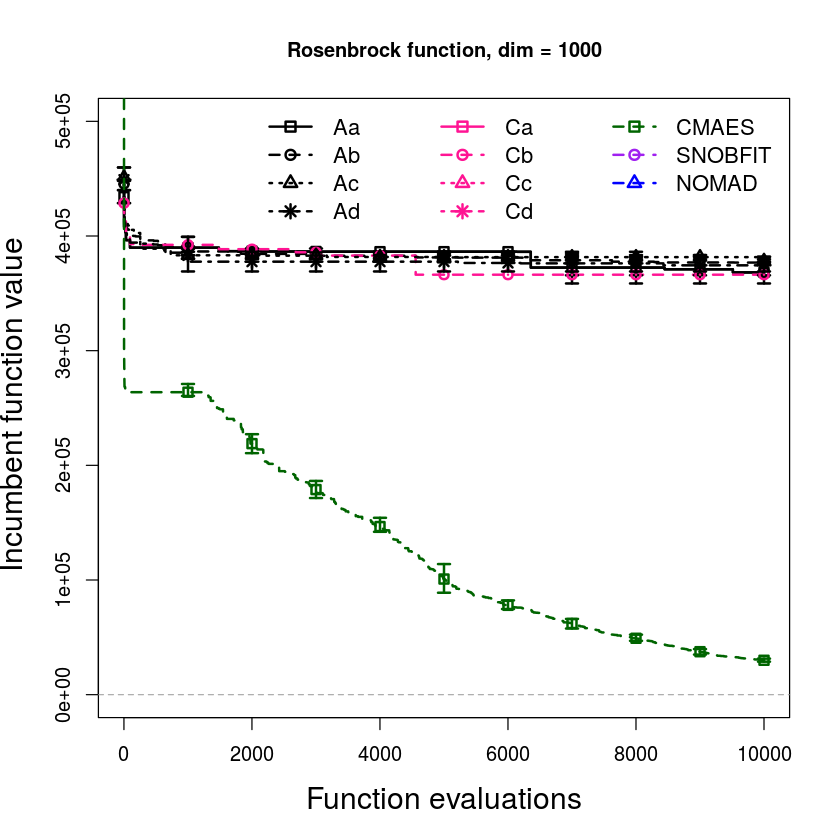}}
\\
\centering
Centered Sinusoidal &{\includegraphics{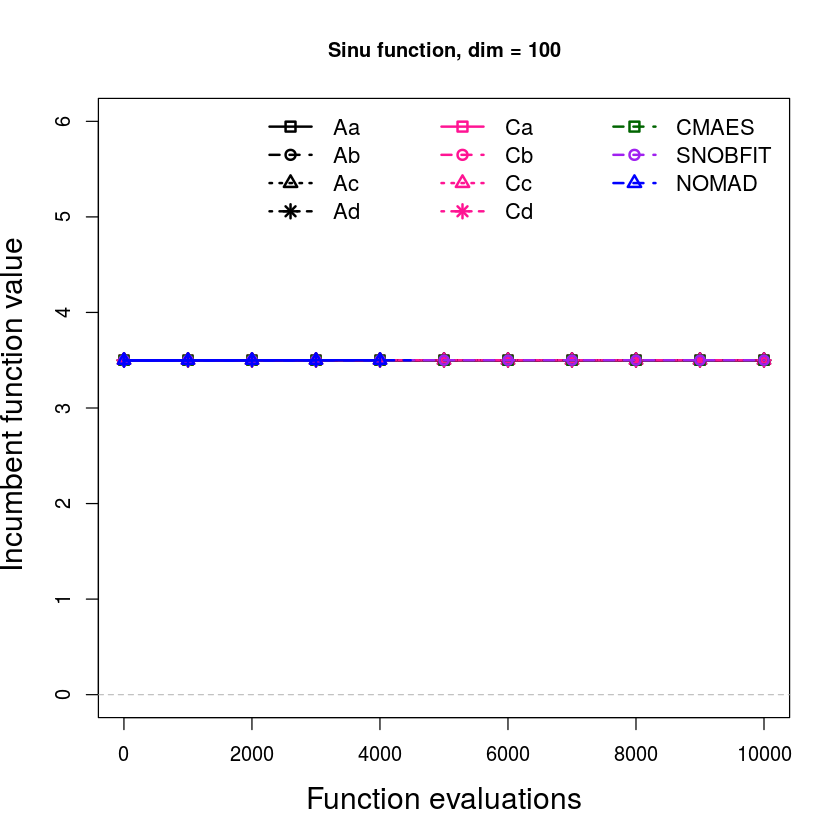}} 
        & {\includegraphics{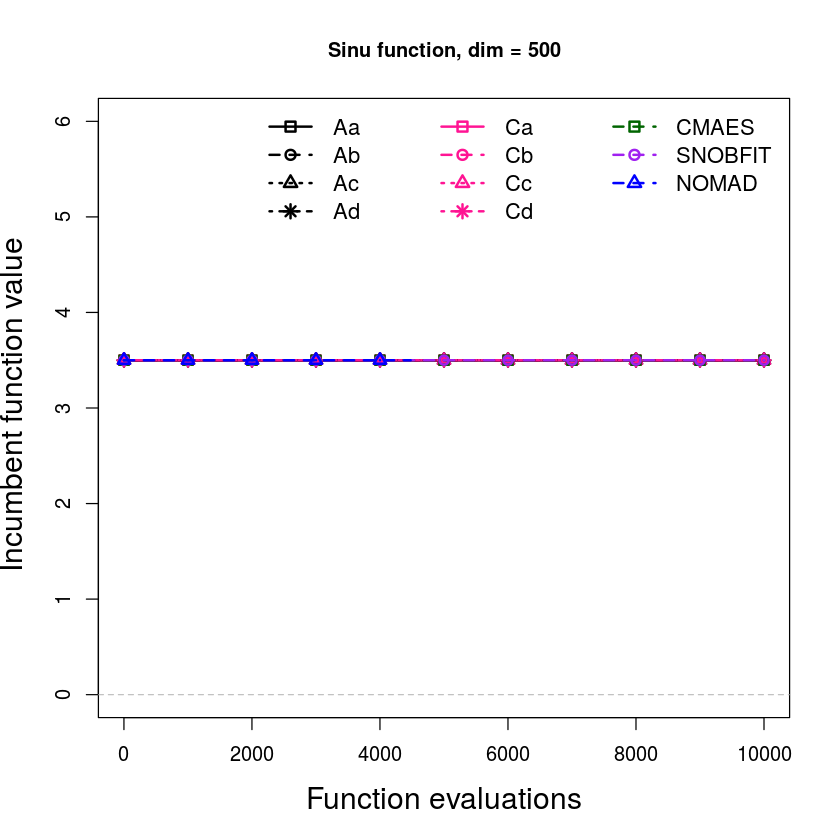}}
        & {\includegraphics{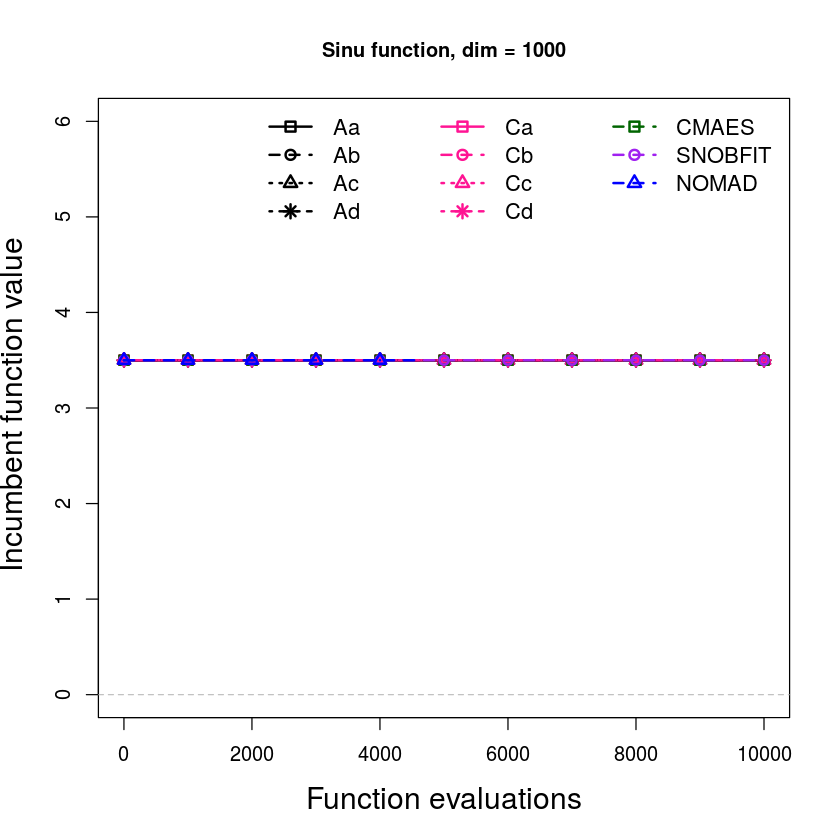}}
\\
\centering
Shifted Sinusoidal &{\includegraphics{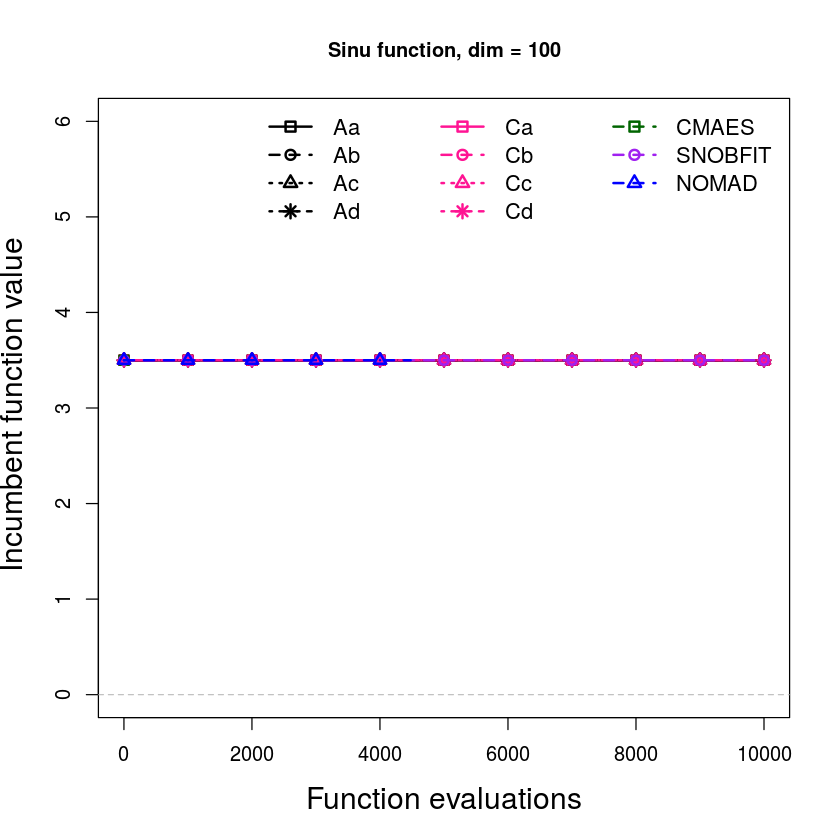}} 
        & {\includegraphics{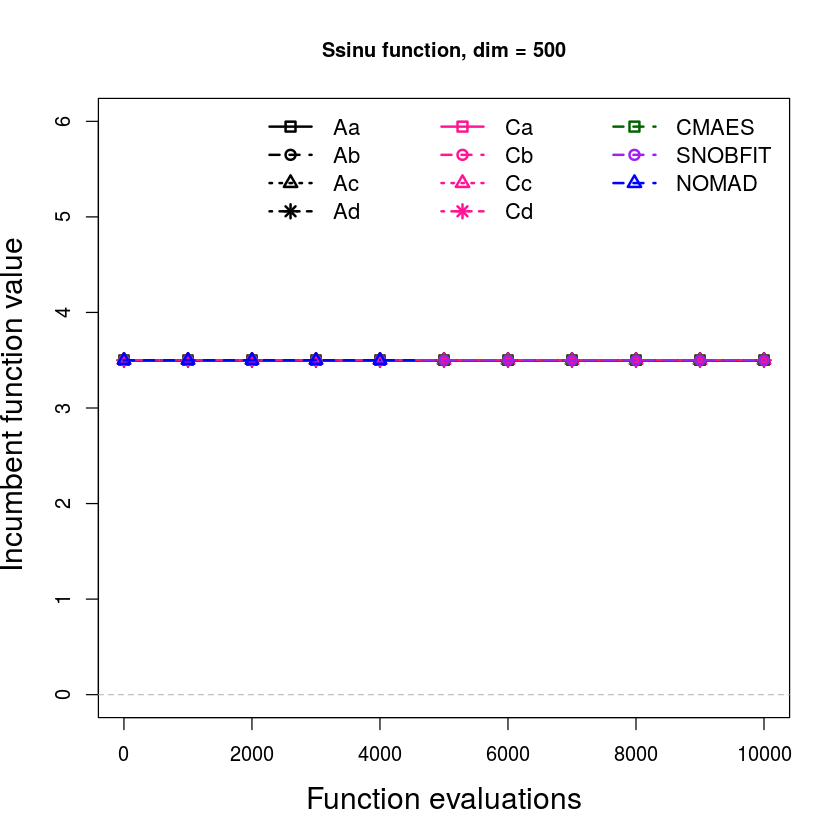}}
         & {\includegraphics{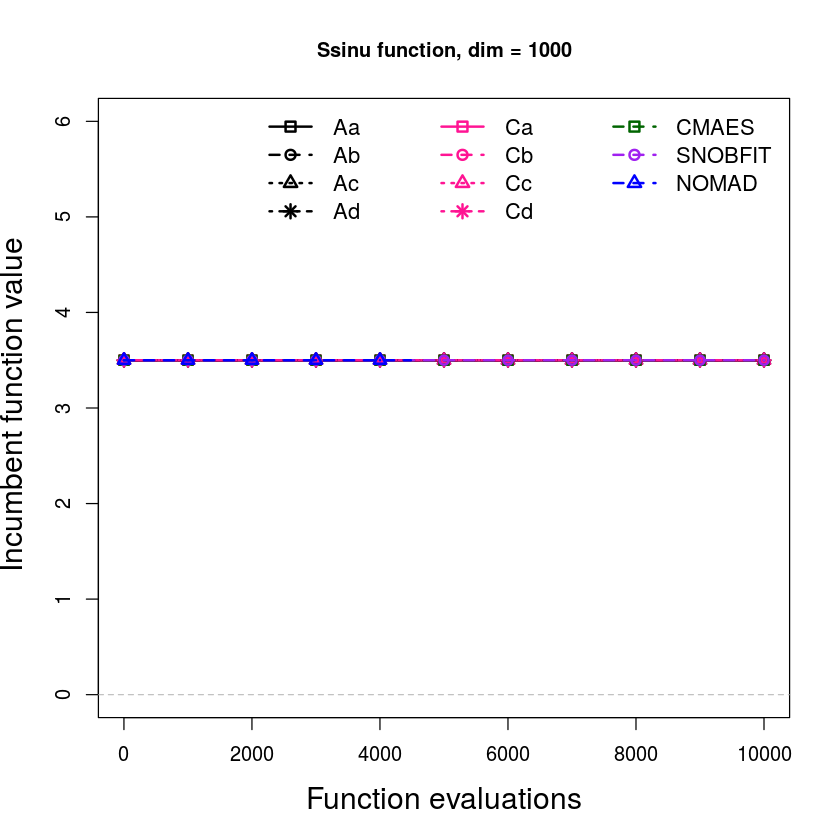}}
\\
    \bottomrule
    \end{tblr}
\caption{
%The incumbent function value vs. number of function evaluations.
Numerical results at dimensions 100, 500, and 1000.}
\label{fig:table_of_figures_exper_2_2}
\end{figure}
\end{appendices}
\backmatter

%%===========================================================================================%%
%% If you are submitting to one of the Nature Portfolio journals, using the eJP submission   %%
%% system, please include the references within the manuscript file itself. You may do this  %%
%% by copying the reference list from your .bbl file, paste it into the main manuscript .tex %%
%% file, and delete the associated \verb+\bibliography+ commands.                            %%
%%===========================================================================================%%

\newpage
\bibliography{BASSOREVISION}% common bib file
%% if required, the content of .bbl file can be included here once bbl is generated
%%\input sn-article.bbl
\end{document}